\documentclass[draft,preprint,12pt,numbers,sort&compress]{elsarticle}
\usepackage{mathrsfs}
\usepackage{amssymb}
\usepackage{amsthm}
\usepackage{mathrsfs}
\usepackage[centertags]{amsmath}
\usepackage{amsfonts}

\usepackage{color}

\input amssym.def
\input amssym.tex

\date{}

\newtheorem{Theorem}{Theorem}[section]

\newtheorem{Lemma}{Lemma}[section]

\newcommand\R{\mbox{\bf R}}

\newcommand\SR{\mbox{\scriptsize\bf R}}

\newcommand{\definition}{{\lower .5ex
  \hbox{$\>\>\stackrel{\triangle}{=}\>\>$} }}
\newcommand\supp{\mathop{\rm supp}}

\begin{document}

\baselineskip=22pt
\thispagestyle{empty}

\mbox{}
\bigskip

\begin{center}

{\Large \bf The Cauchy problem for higher-order modified}\\[1ex]
{\Large \bf Camassa-Holm equations on the circle}

{Wei YAN$^a$,\quad Yongsheng LI $^b$, \quad Xiaoping Zhai$^c$,
\quad Yimin Zhang$^d$}\\[2ex]

{$^a$School of Mathematics and Information Science, Henan Normal University,}\\
{Xinxiang, Henan 453007, P. R. China}\\[1ex]

{$^b$Department of Mathematics, South China University of Technology,}\\
{Guangzhou, Guangdong 510640, P. R. China}\\[1ex]
{$^c$Department of Mathematics, Sun Yat-sen University,
Guangzhou, 510275, China}
{$^d$Wuhan Institute of Physics and Mathematics, Chinese Academy of Sciences, Wuhan, Hubei 430071, P. R. China}\\[1ex]

\end{center}

\noindent{\bf Abstract.}
In this paper, we investigate the Cauchy problem
for the  shallow water type equation
\begin{eqnarray*}
      u_{t}+\partial_{x}^{2j+1}u
     + \frac{1}{2}\partial_{x}(u^{2})+
     \partial_{x}(1-\partial_{x}^{2})^{-1}\left[u^{2}+\frac{1}{2}u_{x}^{2}\right]=0
    \end{eqnarray*}
 with low regularity data in the periodic settings.  Himonas and  Misiolek (Communications in Partial Differential Equations,
23(1998), 123-139.) have proved that the problem is locally well-posed for small initial data in $H^{s}(\mathbf{T})$ with $s\geq-\frac{j}{2}+1,j\in N^{+}$ with the aid of the standard Fourier restriction norm method.
To the best of our knowledge, there is no result of well-posedness about the problem when $s<-\frac{j}{2}+1$.
 In this paper, firstly, we prove that the bilinear estimate related to the nonlinear term of the equation in standard Bourgain space is invalid with $s<-\frac{j}{2}+1.$
  Then we  prove that the Cauchy problem for the periodic shallow water-type
  equation is locally well-posed in $H^{s}(\mathbf{T})$ with $-j+\frac{3}{2}< s<-\frac{j}{2}+1,j\geq2$ for arbitrary initial data.
  The novelty is that
  we introduce   some new function spaces and give a useful relationship among new spaces.

\bigskip

{\large\bf 1. Introduction}
\bigskip

\setcounter{Theorem}{0} \setcounter{Lemma}{0}

\setcounter{section}{1}

In this paper, we consider the Cauchy problem for the periodic shallow water type equation
\begin{eqnarray}
&& u_{t}+\partial_{x}^{2j+1}u
     + \frac{1}{2}\partial_{x}(u^{2})+\partial_{x}(1-\partial_{x}^{2})^{-1}
     \left[u^{2}+\frac{1}{2}u_{x}^{2}\right]=0,\label{1.01}\\
    &&u(x,0)=u_{0}(x),\quad x\in \mathbf{T}=[0,2\pi\lambda),\lambda\geq1. \label{1.02}
\end{eqnarray}

Obviously, (\ref{1.01}) is the higher order modification of the Camassa-Holm equation
\begin{eqnarray}
u_{t}+ \frac{1}{2}\partial_{x}(u^{2})+\partial_{x}(1-\partial_{x}^{2})^{-1}
\left[u^{2}+\frac{1}{2}u_{x}^{2}\right]=0\label{1.03}
 \end{eqnarray}
in the nonlocal form. Equation (\ref{1.03}) has been investigated by many people, for instance, see \cite{BCARMA,CH,CHH,C,Con2000, C2000,C2001,C2006,CE,CECPAM,CS,EY,EYJFA,FF,HMPZ,LIMRN,XZCPAM,XZCPDE}.

Omitting the last term in (\ref{1.01})  yields the higher order Korteweg-de Vries  equation
\begin{eqnarray}
u_{t}+\partial_{x}^{2j+1}u+\frac{1}{2}\partial_{x}(u^{2})=0.\label{1.04}
\end{eqnarray}
Using the  Fourier restriction norm method, Hirayama \cite{H}  proved that (\ref{1.01}) is locally
well-posed in $H^{s}(\mathbf{T})$ with $s\geq-\frac{j}{2}.$
When $j=1$, (\ref{1.03}) reduces to the following KdV equation
which possesses the bi-Hamiltonian structure and completely integrable and infinite conserved laws.
Lots of  people have investigated the Cauchy problem for the KdV
equation, for instance, see \cite{Bourgain93, Bourgain97, KPV1996,KPV2001,L,CKSTT,H,Kis,T}.
Especially, Bourgain \cite{Bourgain93}
introduced the Fourier restriction norm method which is an effective tool in
solving the Cauchy problem for
dispersive equations in low regularity,  to establish the local well-posedness of
the Cauchy problem for the KdV. Kenig et al. \cite{KPV1996} proved that the Cauchy problem for the periodic KdV equation is
locally well-posed in $H^{s}(\mathbf{T})$ with $s\geq-\frac{1}{2}.$    Bourgain \cite{Bourgain97}
proved that the Cauchy problem for the periodic KdV  equation is ill-posed in $H^{s}(\mathbf{T})$
with $s< -\frac{1}{2}$ in the sense that the solution map is not $C^{3}.$ Colliander et al. \cite{CKSTT} proved that the Cauchy problem for the
periodic KdV  equation is globally well-posed in $H^{s}(\mathbf{T})$ with $s\geq -\frac{1}{2}$ with the aid of $I$ method.
Recently, by using the inverse scattering method,
Kappeler and  Topalov \cite{KT2006} proved that   the  Cauchy  problem for the KdV  equation
is globally
well-posed  in $H^{s}(\mathbf{T})$   with $s\geq -1$ in  $H^{s}(\mathbf{T})$.
Molinet \cite{Molinet} proved that the  Cauchy  problem for the KdV  equation  is  ill-posed in $H^{s}(\mathbf{T})$   with $s< -1$.
Many researchers have studied the non-periodic case of the KdV equation, for instance, see \cite{KPV1996, KPV2001,G,Kis}.

Many people have investigated the periodic case and nonperiodic case of (\ref{1.01}) \cite{HM,HM1998,HM2000,Go,O,WC,LYY,LYLH}.
Himonas and  Misiolek \cite{HM1998} have proved that the problem (\ref{1.01}) is locally well-posed for small initial data in $H^{s}(\mathbf{T})$ with $s\geq-\frac{j}{2}+1,j\in N^{+}$ with the aid of the standard Fourier restriction norm method. Himonas and  Misiolek \cite{HM2000} have proved that the problem (\ref{1.01}) with $j=1$ is locally well-posed for any initial data in $H^{s}(\mathbf{T})$ with $s\geq\frac{1}{2}$.
To the best of our knowledge, there is no result about the well-posedness of (\ref{1.01}) when initial data in $H^{s}(\mathbf{T})$ with $s<-\frac{j}{2}+1,j\in N^{+}$. The main difficulty is that the structure of (\ref{1.01})
is complicated. Recently, Yan et al. \cite{YLZZ} proved that the problem (\ref{1.01}) with $j=1$ is locally well-posed for small initial data in $H^{s}(\mathbf{T})$ with $\frac{1}{6}<s<\frac{1}{2}$ with the aid of the new spaces. The spaces of (\ref{1.01}) with  $j\geq2$,$j\in N$ are  different from  theirs of (\ref{1.01}) with  $j=1$ due to different structure.

In recent ten years, to obtain low regularity of dispersive equations, some resolution function spaces have been introduced by some researchers \cite{BT,IK,IKT,Kis,KJDE,KJDE2013,Kato,TKato}.
Choosing a suitable function space is useful  and difficult in dealing  with the low regularity of dispersive equations, for instance, see \cite{IK,IKT,G,BT,KJDE}.
In this paper, firstly, we prove that the bilinear estimate related to the nonlinear term of the equation in $W^{s}$ defined below  is invalid with $s<-\frac{j}{2}+1.$
Then, by introducing the  new function spaces and the Strichartz estimate which are
used to establish the bilinear estimates and using the fixed point Theorem, we prove that
the Cauchy problem for (\ref{1.01}) is locally well-posed in $H^{s}(\mathbf{T})$ with
$-j+\frac{3}{2}< s<1-\frac{j}{2}$ with $j\geq 2,j\in Z$ for arbitrary initial data.

We give some notations before presenting the main results. $A\sim B$ means that $|B|\leq |A|\leq 4|B|$.
 $A\gg B$ means that $|A|\geq 4|B|.$ $C$ is a positive constant
 which may vary from line to line.  $0<\epsilon\ll1$ means that
 $0<\epsilon<\frac{1}{100j^{5}}$.
 Throughout this paper,
 $\dot{Z}:=Z- \{ 0\}$.    Denote
 $dk$ by the normalized counting measure on $\dot{Z}$.
 $(dk)_{\lambda}$ the normalized counting measure on $\dot{Z_{\lambda}}=
 \frac{\dot{Z}}{\lambda}$:
 \begin{eqnarray*}
 \int a(k)(dk)_{\lambda}=\frac{1}{\lambda}\sum_{k\in  \dot{Z}_{\lambda}}a(k).
 \end{eqnarray*}
 Denote $\mathscr{F}_{x}f$ by the Fourier transformation of a function $f$
 defined on
  $[0,2\pi\lambda)$ with  respect to the space variable
 \begin{eqnarray*}
 \mathscr{F}_{x}f(k)=\frac{1}{\sqrt{2\pi}}\int_{0}^{2\pi\lambda}e^{- i kx}f(x)dx.
 \end{eqnarray*}
 and we have the Fourier inverse transformation formula
 \begin{eqnarray*}
 f(x)=\frac{1}{\sqrt{2\pi}}\int e^{ i kx} \mathscr{F}_{x}f(k)(dk)_{\lambda}=\frac{1}{\sqrt{2\pi}}\sum_{k \in
 \dot{Z}}e^{ i kx}\mathscr{F}_{x}f(k).
 \end{eqnarray*}
 Denote $\mathscr{F}_{t}f$ by the Fourier transformation of a function $f$
  with  respect to the time variable
 \begin{eqnarray*}
 \mathscr{F}_{t}f(\tau)=\frac{1}{\sqrt{2\pi}}\int_{\SR}e^{- i t\tau}f(t)dt
 \end{eqnarray*}
 and we have the Fourier inverse transformation formula
 \begin{eqnarray*}
 f(t)=\frac{1}{\sqrt{2\pi}}\int e^{ i t\tau} \mathscr{F}_{t}f(\tau)d\tau.
 \end{eqnarray*}
Let
\begin{eqnarray*}
S(t)\phi(x)=\frac{1}{\sqrt{2\pi}}\int e^{i kx }e^{i (-1)^{j+1}k^{2j+1}}\mathscr{F}_{x}\phi(k)(dk)_{\lambda}.
\end{eqnarray*}
We define the space-time  Fourier transform
$\mathscr{F}f(k,\tau)$ for $k\in \dot{Z}_{\lambda}$ and $\tau\in \R$ by
\begin{eqnarray*}
\mathscr{F}f(k,\tau)=\frac{1}{2\pi}\int\int_{0}^{2\pi\lambda}e^{-i kx}e^{-i \tau t}f(x,t)dxdt
\end{eqnarray*}
and
\begin{eqnarray*}
f(x,t)=\frac{1}{2\pi}\int\int e^{i kx}e^{i \tau t}\mathscr{F}f(k,\tau)(dk)_{\lambda}d\tau.
\end{eqnarray*}
Obviously, we have that
\begin{eqnarray*}
&&\|f\|_{L^{2}(\mathbf{T})}=\|\mathscr{F}_{x}f\|_{L^{2}((dk)_{\lambda})},\\
&&\int_{0}^{2\pi }f(x)\overline{g}(x)dx=\int \mathscr{F}_{x}f(k)\overline{\mathscr{F}_{x}
f}(k)(dk)_{\lambda},\\
&&\mathscr{F}_{x}(fg)=\mathscr{F}_{x}f*\mathscr{F}_{x}g=\int \mathscr{F}_{x}f(k-k_{1})
\mathscr{F}_{x}g(k_{1})(dk_{1})_{\lambda}.
\end{eqnarray*}
Let
\begin{eqnarray*}
&&P(k)=(-1)^{j+1}k^{2j+1},\sigma=\tau-P(k),\quad \sigma_{j}=\tau_{j}-P(k_{j}),\\
&&D_{1}=\left\{(\tau,k)\in \R\times \dot{Z}:|\tau-P(k)|\leq \frac{2j+1}{3}4^{-j}|k|^{2j},|k|\geq1\right\},\\
&&D_{2}=\left\{(\tau,k)\in \R\times \dot{Z}: \frac{2j+1}{3}4^{-j}|k|^{2j}< |\tau-P(k)|< \frac{2j+1}{3}4^{-j}|k|^{2j+1},|k|\geq1\right\},\\
&&D_{3}=\left\{(\tau,k)\in \R\times \dot{Z}:|\tau-P(k)|\geq \frac{2j+1}{3}4^{-j}|k|^{2j+1},|k|\geq1\right\},\\
&&D_{4}=\left\{(\tau,k)\in \R\times \dot{Z}:|\tau-P(k)|> \frac{2j+1}{3}4^{-j}|k|^{2j+1},\frac{1}{\lambda}\leq |k|\leq1\right\},\\
&&D_{5}=\left\{(\tau,k)\in \R\times \dot{Z}:|\tau-P(k)|\leq \frac{2j+1}{3}4^{-j}|k|^{2j+1},\frac{1}{\lambda}\leq |k|\leq1\right\},\\
&&\mathscr{F}(\Lambda^{-1})f=\langle \sigma \rangle ^{-1}\mathscr{F}f, \mathscr{F}J^{s}f=\langle k\rangle ^{s}\mathscr{F}f(k).
\end{eqnarray*}
The Sobolev space $H^{s}(\mathbf{T})$ is defined by the following  norm
\begin{eqnarray*}
\|f\|_{H^{s}(\mathbf{T})}=\|\langle k\rangle^{s}\mathscr{F}_{x}f(k)\|_{L^{2}(dk)_{\lambda}}
\end{eqnarray*}
and define the $X_{s,b}$ spaces for $2\pi$-periodic function via the norm
\begin{eqnarray*}
\|u\|_{X_{s,b}(\mathbf{T}\times \SR)}=\left\|\langle k\rangle^{s} \left\langle \sigma\right\rangle^{b}\mathscr{F}u(k,\tau)\right\|_{L^{2}((dk)_{\lambda}(d\tau))}.
\end{eqnarray*}
The $Z^{s}$ space is equipped with the following norm
\begin{eqnarray*}
\|u\|_{Z^{s}}=\|P_{D_{1}\cup D_{5}}u\|_{X_{s,\frac{2j-1}{2j}}}+\|P_{D_{2}}u\|_{X_{(1-2j)(s-1),s}}+
\|P_{D_{3}\cup D_{4}}u\|_{X_{-\frac{s-1}{j}-1,\frac{s-1}{j}+1}}+\left\|u\right\|_{Y^{s}},
\end{eqnarray*}
where $j\geq 2$ and $\|u\|_{Y^{s}}=\left\|\langle k\rangle^{s}
\mathscr{F}u(k,\tau)\right\|_{L^{2}(dk)_{\lambda}L^{1}(d\tau)}.$
Let
\begin{eqnarray*}
\|u\|_{W^{s}}=\|u\|_{X_{s,\frac{1}{2}}}+\|u\|_{Y^{s}}.
\end{eqnarray*}

The main result of this paper are as follows.

\begin{Theorem}\label{Thm1}
Let $s<\frac{2-j}{2}$,
\begin{eqnarray*}
F(u_{1},u_{2})=\frac{1}{2}\partial_{x}(u_{1}u_{2})+
\partial_{x}(1-\partial_{x}^{2})^{-1}\left[u_{1}u_{2}+\frac{1}{2}(\partial_{x}u_{1})(\partial_{x}u_{2})\right],\label{1.06}
\end{eqnarray*}
and $u_{j}(j=1,2)$ be $2\pi$-periodic
functions.
Then, we obtain that
\begin{eqnarray*}
\left\|\mathscr{F}^{-1}\left[\langle\tau+(-1)^{j}k^{2j+1}\rangle^{-1}\mathscr{F}F(u_{1},u_{2})\right]\right\|_{W^{s}}\leq C\prod_{j=1}^{2}\|u_{j}\|_{W^{s}}
\end{eqnarray*}
is untrue.
\end{Theorem}
\noindent {\bf Remark 1.}Theorem 1.1 implies that the standard Fourier restriction norm method is invalid when $s<1-\frac{j}{2}.$
Lack of the bilinear estimates in $W^{s}$ doesnot necessarily imply ill-posedness of problem. One can recover the bilinear estimates by changing new function spaces, for instance, see \cite{BT,G,IK,IKT}. Thus, by choosing suitable function spaces, we obtain the following  Theorem 1.2.

\begin{Theorem}\label{Thm2}
Let $ -j+\frac{3}{2}< s<-\frac{j}{2}+1$, $j\geq3$  and $u_{0}$ be
$2\pi\lambda$-periodic
function.
Then the Cauchy problems (\ref{1.01})(\ref{1.02})
are locally well-posed in $H^{s}(\mathbf{T})$. More  precisely,
for any $u_{0}\in H^{s}$ with  $\|u_{0}\|_{H^{s}}\leq r,$ there exists a solution $u\in C([-T,T];H^{s})$ to (\ref{1.01})(\ref{1.02}) with $T=T(r)>0.$ Moreover, the solution is uniquely derived in $Z_{T}^{s}$ embedded continuously into $C([-T,T];H^{s})$ and the data-to-solution map from
$\left\{u_{0}\in H^{s}| \quad \|u_{0}\|_{H^{s}}\leq r\right\}$ to $Z_{T}^{s}$ is Lipschitz.
\end{Theorem}

\noindent {\bf Remark 2.}
 The optimal regularity indices of the Cauchy problem for (\ref{1.01}) is unknown. We will  pursue the optimal regularity indices of the Cauchy problem for (\ref{1.01}).
 From Lemmas 2.1, 2.4 and the structure  of (\ref{1.01}), we choose the space $X_{s,\frac{2j-1}{2j}}$ related to $D_{1},D_{5}$. Since we consider the case $s<1-\frac{j}{2}$,  we choose the space
 $X_{-\frac{s-1}{j}-1,\frac{s-1}{j}+1}$ related to $D_{3},D_{4}$ in view of high $\times$ high $\rightarrow$ low interaction. By a direct computation, we know that $X_{(1-2j)(s-1),s}$ related to $D_{2}$ is suitable.

The rest of the paper is arranged as follows. In Section 2,  we present some
preliminaries. In Section 3, we present  some  bilinear estimates. In Section 4,
we present the proof of Theorem 1.1.
In Section 5, we present the proof of Theorem 1.2.

\bigskip
\bigskip
\bigskip

 \noindent{\large\bf 2. Preliminaries }

\setcounter{equation}{0}

\setcounter{Theorem}{0}

\setcounter{Lemma}{0}

\setcounter{section}{2}

In this section, we give some preliminaries which are crucial in establishing Lemmas 3.1, 3.2 and Theorems 1.1,1.2.

\begin{Lemma}\label{Lemma2.1}
Let $u(x,t),v(x,t)$ be  $2\pi \lambda$-periodic functions and $a+b\geq \frac{j+1}{2j+1}$
and ${\rm min}\{a,b\}>\frac{1}{2(2j+1)}$. Then, we have that
\begin{eqnarray}
     && \left\|uv\right\|_{L_{xt}^{2}}\leq C\|u\|_{X_{0,a}(\mathbf{T} \times \SR)}
     \|v\|_{X_{0,b}(\mathbf{T}  \times \SR)},
       \label{2.01}\\
       &&\left\|uv\right\|_{X_{0,-a}}\leq C\|u\|_{X_{0,b}(\mathbf{T}  \times \SR)}
       \|v\|_{L_{xt}^{2}}.
       \label{2.02}
\end{eqnarray}
\end{Lemma}

Lemma 2.1 can be proved similarly to  Lemma 2.3 of \cite{YJLH}.

\begin{Lemma}\label{Lemma2.2}
Assume that  $-j+\frac{3}{2}+j\epsilon\leq s\leq 1-\frac{j}{2}-j\epsilon\in \R$ and  $T>0$. Then, we have that
\begin{eqnarray*}
\left\|\eta(t)S(t)\phi\right\|_{Z^{s}}\leq C\|\phi\|_{H^{s}(\mathbf{T})}.
\end{eqnarray*}
\end{Lemma}
{\bf Proof.} Combining the definition of $Z^{s}$ with Lemma 2.3, we have that $X^{\frac{2j-1}{2j}}
\hookrightarrow Z^{s}\hookrightarrow C([0,T]:H^{s}(\mathbf{T})).$

We have completed the proof of Lemma 2.2.

\begin{Lemma}\label{Lemma2.3}Let $-j+\frac{3}{2}+j\epsilon\leq s\leq1-\frac{j}{2}-j\epsilon$ and  $T>0$. Then, we have that
\begin{eqnarray*}
\left\|\eta(t) \int_{0}^{t}S(t-\tau)\partial_{x}(uv)d\tau\right\|_{Z^{s}}
\leq C\|\partial_{x}\Lambda^{-1}(uv)\|_{Z^{s}}.
\end{eqnarray*}
\end{Lemma}

For the proof of Lemma 2.3, we refer the readers to Lemma 2.3 of \cite{KJDE}.

\begin{Lemma}\label{Lemma2.4}Let $-j+\frac{3}{2}+j\epsilon\leq 1-\frac{j}{2}-j\epsilon$ and $j\geq 2,j\in Z$. Then, we have that
\begin{eqnarray}
&&\|u\|_{X_{s,\frac{1}{2j}}}\leq C\|u\|_{Z^{s}}\leq C\|u\|_{X_{s,\frac{2j-1}{2j}}},\label{2.03}\\
&&\|u\|_{X_{s,\frac{1}{2}}(D_{1}\bigcup D_{2})}\leq C\|u\|_{Z^{s}(D_{1}\bigcup D_{2})}.\label{2.04}
\end{eqnarray}
\end{Lemma}
\noindent {\bf Proof.} We firstly prove that (\ref{2.03}). When $\supp \mathscr{F}u\subset D_{1}$,
since $\frac{2j-1}{2j}\geq \frac{1}{2j},$ we have that
$\|u\|_{X_{s,\frac{2j-1}{2j}}}\geq \|u\|_{X_{s,\frac{1}{2j}}}$. When $\supp \mathscr{F}u\subset D_{2}$,
 since $-j+\frac{3}{2}+j\epsilon\leq s\leq 1-\frac{j}{2}-j\epsilon$,  we have that
$\langle \sigma\rangle ^{s-\frac{1}{2j}}\geq C\langle k\rangle^{2js-2j+1}$ which yields
that $\langle k\rangle^{s}\langle \sigma\rangle^{\frac{1}{2j}}\leq C\langle k\rangle^{(1-2j)(s-1)}\langle \sigma\rangle^{s}$,
thus, we have that $\|u\|_{X_{(1-2j)(s-1),s}}\geq \|u\|_{X_{s,\frac{1}{2j}}}$. When $\supp \mathscr{F}u\subset D_{3}$,
since $-j+\frac{3}{2}+j\epsilon\leq s\leq1-\frac{j}{2}-j\epsilon$,  we have that
$\langle \sigma\rangle ^{\frac{s-1}{j}+\frac{2j-1}{2j}}\geq C\langle k\rangle^{s+1+\frac{s-1}{j}}$ which yields
that $|k|^{s}\langle \sigma\rangle^{\frac{1}{2j}}\leq C\langle k\rangle^{-\frac{s-1}{j}-1}\langle \sigma\rangle^{\frac{s-1}{j}+1}$,
thus, we have that $\|u\|_{X_{-\frac{s-1}{j}-1,\frac{s-1}{j}+1}}\geq \|u\|_{X_{s,\frac{1}{2j}}}$.
Consequently, we have that $\|u\|_{Z^{s}}\geq C\|u\|_{X_{s,\frac{1}{2j}}}.$
When $\supp \mathscr{F}u\subset D_{2}$,
 since $-j+\frac{3}{2}+j\epsilon\leq s\leq 1-\frac{j}{2}-j\epsilon$,  we have that
$\langle \sigma\rangle ^{-s+\frac{2j-1}{2j}}\geq C\langle k\rangle^{2js+1}$ which yields
that $\langle k\rangle^{(1-2j)(s-1)}\langle \sigma\rangle^{s}\leq C\langle k\rangle^{s}
\langle \sigma\rangle^{\frac{2j-1}{2j}}$,
thus, we have that $\|u\|_{X_{(1-2j)(s-1),s}}\leq C\|u\|_{X_{s,\frac{2j-1}{2j}}}$.
When $\supp \mathscr{F}u\subset D_{3}$,
since $-j+\frac{3}{2}+j\epsilon\leq s\leq 1-\frac{j}{2}-j\epsilon$,  we have that
$\langle \sigma\rangle ^{-\frac{s-1}{j}-\frac{1}{2j}}\geq C\langle k\rangle^{-s-1-\frac{s-1}{j}}$ which yields
that $\langle k\rangle^{-\frac{s-1}{j}-1}\langle \sigma\rangle^{\frac{s-1}{j}+1}
\leq C\langle k\rangle^{s}\langle \sigma\rangle^{\frac{2j-1}{2j}}$,
thus, we have that $\|u\|_{X_{-\frac{s-1}{j}-1,\frac{s-1}{j}+1}}\leq C \|u\|_{X_{s,\frac{2j-1}{2j}}}$.
Consequently, we have that $\|u\|_{Z^{s}}\leq C\|u\|_{X_{s,\frac{2j-1}{2j}}}.$
By Cauchy-Schwartz inequality with respect to $\tau$, we have that
$\|\langle k\rangle^{s}\mathscr{F}u\|_{l_{k}^{2}l_{\tau}^{1}}\leq C
\|u\|_{X_{s,\frac{2j-1}{2j}}},$  consequently, we have that $\|u\|_{Z^{s}}
\leq C\|u\|_{X_{s,\frac{2j-1}{2j}}}.$  Now we prove (\ref{2.04}). When
$\supp \mathscr{F}u\subset D_{1}$, since $\frac{2j-1}{2j}\geq \frac{1}{2},$ we have that
$\|u\|_{X_{s,\frac{2j-1}{2j}}}\geq \|u\|_{X_{s,\frac{1}{2}}}$. When
$\supp \mathscr{F}u\subset D_{2}$, since $s\geq-j+\frac{3}{2}+j\epsilon$, we have that
$\langle k\rangle^{s}\langle\sigma\rangle^{1/2}\leq C\langle k\rangle
^{(1-2j)(s-1)}\langle\sigma\rangle^{s},$ consequently, we have that
$\|u\|_{X_{(1-2j)(s-1),\>s}}\geq \|u\|_{X_{s,\frac{1}{2}}}$.

We have completed the proof of Lemma 2.4.

\begin{Lemma}\label{Lemma2.5}
Let $k=k_{1}+k_{2}$, $\tau=\tau_{1}+\tau_{2}$ and
\begin{eqnarray*}
\sigma=\tau-k^{2j+1},\sigma_{1}= \tau_{1}-k_{1}^{2j+1},\sigma_{2}= \tau_{2}-k_{2}^{2j+1}.
\end{eqnarray*}
 Then, we have that
\begin{eqnarray*}
3{\rm max}\left\{|\sigma|,|\sigma_{1}|,|\sigma_{2}|\right\}\geq|\sigma-\sigma_{1}-\sigma_{2}|=
\left|k^{2j+1}-k_{1}^{2j+1}-k_{2}^{2j+1}\right|\geq \frac{2j+1}{4^{j}}|k_{min}||k_{max}|^{2j}.
\end{eqnarray*}
where
\begin{eqnarray*}
|k_{min}|={\rm min}\left\{|k|,|k_{1}|,|k_{2}|\right\},|k_{max}|={\rm max}\left\{|k|,|k_{1}|,|k_{2}|\right\}.
\end{eqnarray*}
Moreover, we have that one of three following cases must occur:
\begin{eqnarray}
&& (a): |\sigma|={\rm max}\left\{|\sigma|,|\sigma_{1}|,|\sigma_{2}|\right\}\geq \frac{2j+1}{3}4^{-j}|k_{min}||k_{max}|^{2j},\label{2.05}\\
&& (b): |\sigma_{1}|={\rm max}\left\{|\sigma|,|\sigma_{1}|,|\sigma_{2}|\right\}\geq \frac{2j+1}{3}4^{-j}|k_{min}||k_{max}|^{2j},\label{2.06}\\
&& (c): |\sigma_{2}|={\rm max}\left\{|\sigma|,|\sigma_{1}|,|\sigma_{2}|\right\}\geq \frac{2j+1}{3}4^{-j}|k_{min}||k_{max}|^{2j}.\label{2.07}
\end{eqnarray}

\end{Lemma}

Proof  of Lemma 2.5 can be seen in  Lemma 2.4 of \cite{LYY}.

\begin{Lemma}\label{Lemma2.6}
Let $Z=\R/2\pi\lambda$, $\lambda>0$. Let $s\in \R$ and $\mathscr{X}^{s}$ be a Banach space of functions on $\R_{t}\times Z$ with the following properties:
(i) $\mathcal {S}(\R\times Z)$ is dense in $\mathscr{X}^{s}$,
(ii)  $X^{s,b}(\R\times Z)\hookrightarrow \mathscr{X}^{s} \hookrightarrow C_{t}(\R; H^{s}(Z))$ for some $b>\frac{1}{2}$,
(iii) $X^{s^{'},b^{'}}(\R\times Z)\hookrightarrow \mathscr{X}^{s}$ for some $s^{'}\in \R$ and $\frac{1}{2}\leq b^{'}<1$.
Assume that $u\in \mathscr{X}^{s}$ satisfies $u(\cdot,0)=0$ in $H^{s}(Z)$. Then, we have
\begin{eqnarray}
\lim\limits_{T\rightarrow +0}\|u\|_{\mathscr{X}_{T}^{s}}=0.
\end{eqnarray}
\end{Lemma}

For the proof of Lemma 2.6, we refer the readers to  Proposition 2.6 of \cite{KJDE}.

\noindent {\bf Remark 3.} From Lemma 2.4, we have that $X_{s,\frac{2j-1}{2j}}\hookrightarrow Z^{s}\hookrightarrow Y^{s}$.
It is easily checked that $Y^{s}\hookrightarrow C_{t}(\R; H^{s}(Z))$. Consequently, $Z^{s}$ satisfies (i)-(iii).

\bigskip
\bigskip

\noindent{\large\bf 3. Bilinear estimates }

\setcounter{equation}{0}

 \setcounter{Theorem}{0}

\setcounter{Lemma}{0}

 \setcounter{section}{3}
 In this section, we present some crucial bilinear estimates. We always assume that $j\geq2,j\in N^{+}$.

 \begin{Lemma}\label{Lemma3.1}
Let $j\geq 2$ and $-j+\frac{3}{2}+j\epsilon\leq s\leq1-\frac{j}{2}-j\epsilon$. Then, we have that
\begin{eqnarray}
      \left\|\Lambda^{-1}\partial_{x}(1-\partial_{x}^{2})^{-1}\prod_{j=1}^{2}(\partial_{x}u_{j})\right\|_{X^{s}}\leq C\prod\limits_{j=1}^{2}\|u_{j}\|_{Z^{s}},
        \label{3.01}
\end{eqnarray}
here $C>0$,  which  is  independent of  $\lambda$, $\left\|\cdot\right\|_{X^{s}}$ is the norm removing $\left\|\cdot\right\|_{Y^{s}}$ from $\left\|\cdot\right\|_{Z^{s}}.$
\end{Lemma}
{\bf Proof.} Obviously, $\left(\R\times\dot{Z}_{\lambda}\right)^{2}\subset \bigcup\limits_{j=1}^{8}\Omega_{j},$
where
\begin{eqnarray*}
&&\Omega_{1}=\left\{(\tau_{1},k_{1},\tau,k)\in \left(\R\times\dot{Z_{\lambda}}\right)^{2}:
{\rm max}\left\{|k_{1}|, |k|\right\}\leq1\right\},\\
&&\Omega_{2}=\left\{(\tau_{1},k_{1},\tau,k)\in \left(\R\times\dot{Z_{\lambda}}\right)^{2}\cap
\Omega_{1}^{c}:|k_{1}|\sim |k_{2}|\gg |k|\geq1\right\},\\
&&\Omega_{3}=\left\{(\tau_{1},k_{1},\tau,k)\in \left(\R\times\dot{Z_{\lambda}}\right)^{2}\cap
\Omega_{1}^{c}:|k_{1}|\sim |k_{2}|\gg |k|,1\geq|k|\geq\frac{1}{\lambda}\right\},\\
&&\Omega_{4}=\left\{(\tau_{1},k_{1},\tau,k)\in \left(\R\times\dot{Z_{\lambda}}\right)^{2}\cap
\Omega_{1}^{c}:|k|\sim |k_{2}|\gg |k_{1}|\geq 1\right\},\\
&&\Omega_{5}=\left\{(\tau_{1},k_{1},\tau,k)\in \left(\R\times\dot{Z_{\lambda}}\right)^{2}\cap
\Omega_{1}^{c}:|k|\sim |k_{2}|\gg |k_{1}|,1\geq|k_{1}|\geq \frac{1}{\lambda}\right\},\\
&&\Omega_{6}=\left\{(\tau_{1},k_{1},\tau,k)\in \left(\R\times\dot{Z_{\lambda}}\right)^{2}\cap
\Omega_{1}^{c}:|k|\sim |k_{1}|\gg |k_{2}|\geq 1\right\},\\
&&\Omega_{7}=\left\{(\tau_{1},k_{1},\tau,k)\in \left(\R\times\dot{Z_{\lambda}}\right)^{2}\cap
\Omega_{1}^{c}:|k|\sim |k_{1}|\gg |k_{2}|,1\geq|k_{2}|\geq \frac{1}{\lambda}\right\},\\
&&\Omega_{8}=\left\{(\tau_{1},k_{1},\tau,k)\in \left(\R\times\dot{Z_{\lambda}}\right)^{2}\cap
\Omega_{1}^{c}:|k|\sim |k_{1}|\sim |k_{2}|\geq 1\right\}.
\end{eqnarray*}
(1) In region $\Omega_{1}$.
By using Lemma 2.5 and the Young inequality, since ${\rm max}\left\{|k_{1}|, |k|\right\}\leq1$ and  $-j+\frac{3}{2}+j\epsilon\leq s\leq1-\frac{j}{2}-j\epsilon$, from the definition of $Z^{s},$  we have that
\begin{eqnarray*}
&&\left\|\Lambda^{-1}\partial_{x}(1-\partial_{x}^{2})^{-1}\prod_{j=1}^{2}(\partial_{x}u_{j})\right\|_{X^{s}}\leq C\left\|\Lambda^{-1}\partial_{x}(1-\partial_{x}^{2})^{-1}\prod_{j=1}^{2}(\partial_{x}u_{j})\right\|_{X_{s,\frac{2j-1}{2j}}}\nonumber\\
&&\leq C\left\||k|\langle \sigma\rangle^{-\frac{1}{2j}}\left(\mathscr{F}u_{1}*\mathscr{F}u_{2}\right)
\right\|_{l_{k}^{2}L_{\tau}^{2}}\nonumber\\
&&\leq C\||k|\|_{l_{k}^{2}}\left\|\mathscr{F}u_{1}*\mathscr{F}u_{2}\right\|_{l_{k}^{\infty}L_{\tau}^{2}}\leq C\|\mathscr{F}u_{1}\|_{l_{k}^{2}L_{\tau}^{2}}\|\mathscr{F}u_{2}\|_{l_{k}^{2}L_{\tau}^{1}}\leq C\|u_{1}\|_{X_{s,\frac{1}{2j}}}\|u_{2}\|_{Y^{s}}\nonumber\\&&\leq C\prod_{j=1}^{2}\|u_{j}\|_{Z^{s}}.
\end{eqnarray*}
(2) In region $\Omega_{2}$. In this  region, we consider (a)-(c) of Lemma 2.5, respectively.

\noindent(a) Case $|\sigma|={\rm max}\left\{|\sigma|,|\sigma_{1}|,|\sigma_{2}|\right\}.$ In this case, we have that
$\supp \left[\mathscr{F}u_{1}*\mathscr{F}u_{2}\right]\subset D_{3}.$

\noindent
When $\supp \mathscr{F}u_{j} \subset \Omega_{1}\cup \Omega_{2}$ with $j=1,2$,
 by using Lemmas 2.5, 2.3, 2.1,  since $-j+\frac{3}{2}+j\epsilon\leq s\leq1-\frac{j}{2}-j\epsilon$,  we have that
\begin{eqnarray*}
&&\left\|\Lambda^{-1}\partial_{x}(1-\partial_{x}^{2})^{-1}\prod_{j=1}^{2}(\partial_{x}u_{j})\right\|_{X^{s}}\leq C\left\|\langle k\rangle^{-\frac{s-1}{j}}\langle\sigma\rangle^{\frac{s-1}{j}}\left[(\langle k\rangle\mathscr{F}u_{1})*(\langle k\rangle\mathscr{F}u_{2})\right]\right\|_{l_{k}^{2}L_{\tau}^{2}}\nonumber\\
&&\leq C\|(J^{s}u_{1})(J^{s}u_{2})\|_{L_{xt}^{2}}\leq C\|u_{1}\|_{X_{s,\frac{1}{2}}}\|u_{2}\|_{X_{s,\frac{1}{2(2j+1)}+\epsilon}}\nonumber\\
&&\leq C\|u_{1}\|_{X_{s,\frac{1}{2}}}\|u_{2}\|_{X_{s,\frac{1}{2j}}}\leq C\prod_{j=1}^{2}\|u_{j}\|_{Z^{s}}.
\end{eqnarray*}
When $\supp \mathscr{F}u_{1} \subset \Omega_{3}$,
 by using  Plancherel identity and the H\"older inequality as well as Lemma 2.5, since  $-j+\frac{3}{2}+j\epsilon\leq s\leq1-\frac{j}{2}-j\epsilon$,  we have that
\begin{eqnarray*}
&&\left\|\Lambda^{-1}\partial_{x}(1-\partial_{x}^{2})^{-1}\prod_{j=1}^{2}
(\partial_{x}u_{j})\right\|_{X^{s}}\leq C\left\|\langle k\rangle^{-\frac{s-1}{j}}\langle\sigma\rangle^{\frac{s-1}{j}}\left[\langle k\rangle\mathscr{F}u_{1}*(\langle k\rangle\mathscr{F}u_{2})\right]\right\|_{l_{k}^{2}L_{\tau}^{2}}\nonumber\\
&&\leq C\left\|\mathscr{F}u_{1}*\left[\langle k\rangle^{2s}\mathscr{F}u_{2}\right]\right\|_{l_{k}^{2}L_{\tau}^{2}}\nonumber\\
&&\leq C\left\|\left[\langle k\rangle^{-\frac{s-1}{j}-1}\langle\sigma\rangle^{\frac{s-1}{j}+1}\mathscr{F}u_{1}\right]*\left[\langle k\rangle^{-2j}\mathscr{F}u_{2}\right]\right\|_{l_{k}^{2}L_{\tau}^{2}}
\nonumber\\
&&\leq C\|u_{1}\|_{X_{-\frac{s-1}{j}-1,\frac{s-1}{j}+1}}\|\langle k\rangle^{-2j+2}\mathscr{F}u_{2}\|_{l_{k}^{1}L_{\tau}^{1}}\nonumber\\
&&\leq C\|u_{1}\|_{X_{-\frac{s-1}{j}-1,\frac{s-1}{j}+1}}\|\langle k\rangle^{-2j+2-s}\langle k\rangle^{s}\mathscr{F}u_{2}\|_{l_{k}^{1}L_{\tau}^{1}}\nonumber\\
&&\leq C\|u_{1}\|_{X_{-\frac{s-1}{j}-1,\frac{s-1}{j}+1}}\|\langle k\rangle^{s}\mathscr{F}u_{2}\|_{l_{k}^{2}L_{\tau}^{1}}\leq
C\prod_{j=1}^{2}\|u_{j}\|_{Z^{s}}.
\end{eqnarray*}
When $\supp \mathscr{F}u_{2} \subset \Omega_{3}$,
this case can be proved similarly to $\supp \mathscr{F}u_{1} \subset \Omega_{3}$.

\noindent (b) Case $|\sigma_{1}|={\rm max}\left\{|\sigma|,|\sigma_{1}|,|\sigma_{2}|\right\},$
 we consider the following cases:
\begin{eqnarray*}
(i): |\sigma_{1}|>4{\rm max}\left\{|\sigma|,|\sigma_{2}|\right\},
(ii):|\sigma_{1}|\leq4{\rm max}\left\{|\sigma|,|\sigma_{2}|\right\},
\end{eqnarray*}
respectively.

\noindent
When (i) occurs: we consider $\supp \mathscr{F}u_{1}\subset D_{1}$, $\supp \mathscr{F}u_{1}\subset D_{2}$, respectively.

\noindent When $\supp \mathscr{F}u_{1}\subset D_{1}$ which yields that $|k|\leq C$, by using Lemmas 2.3, 2.5, 2.1,
$-j+\frac{3}{2}+j\epsilon\leq s\leq1-\frac{j}{2}-j\epsilon$,
 we have that
\begin{eqnarray*}
&&\left\|\Lambda^{-1}\partial_{x}(1-\partial_{x}^{2})^{-1}\prod_{j=1}^{2}
(\partial_{x}u_{j})\right\|_{X^{s}}
\leq C\left\|\langle k\rangle ^{s-1}\langle \sigma \rangle ^{-\frac{1}{2j}}
(\langle k\rangle\mathscr{F}u_{1})*(\langle k\rangle\mathscr{F}u_{2})\right\|_{l_{k}^{2}L_{\tau}^{2}}\nonumber\\&&
\leq C\left\|\langle \sigma \rangle ^{-\frac{1}{2j}}\left[(\langle k\rangle^{s}\langle\sigma\rangle^{\frac{2j-1}{2j}}
\mathscr{F}u_{1})*(\langle k\rangle ^{-s-2j+3}\mathscr{F}u_{2})\right]\right\|_{l_{k}^{2}L_{\tau}^{2}}\nonumber\\
&&\leq C\left\|\left(J^{s}\Lambda ^{\frac{2j-1}{2j}}u_{1}\right)
\left(J^{-s-2j+3}u_{2}\right)\right\|_{X_{0,-\frac{1}{2j}}}\nonumber\\
&&\leq C\|u_{1}\|_{X_{s,\frac{2j-1}{2j}}}\|u_{2}\|_{X_{-s-2j+3,\frac{1}{2}}}\leq C\|u_{1}\|_{X_{s,\frac{2j-1}{2j}}}\|u_{2}\|_{X_{s,\frac{1}{2}}}\leq
C\prod_{j=1}^{2}\|u_{j}\|_{Z^{s}}.
\end{eqnarray*}
When $\supp \mathscr{F} u_{1}\subset D_{2},$  since $-j+\frac{3}{2}+j\epsilon\leq s\leq1-\frac{j}{2}-j\epsilon$, by using Lemmas 2.3, 2.5, 2.1,
 we have that
\begin{eqnarray*}
&&\left\|\Lambda^{-1}\partial_{x}(1-\partial_{x}^{2})^{-1}\prod_{j=1}^{2}(\partial_{x}u_{j})\right\|_{X^{s}}
\leq C\left\|\langle k\rangle ^{s-1}\langle \sigma \rangle ^{-\frac{1}{2j}}
\left[\langle k\rangle\mathscr{F}u_{1}*\langle k\rangle\mathscr{F}u_{2}\right]\right\|_{l_{k}^{2}L_{\tau}^{2}}\nonumber\\&&
\leq C\left\|\left(J^{(1-2j)(s-1)}\Lambda ^{s}u_{1}\right)\left(J^{-s-2j+3}u_{2}\right)
\right\|_{X_{0,-\frac{1}{2j}}}\nonumber\\
&&\leq C\|u_{1}\|_{X_{(1-2j)(s-1),s}}\|u_{2}\|_{X_{-s-2j+3,\frac{1}{2}}}\nonumber\\
&&\leq C\|u_{1}\|_{X_{(1-2j)(s-1),s}}\|u_{2}\|_{X_{s,\frac{1}{2}}}\leq
C\prod_{j=1}^{2}\|u_{j}\|_{Z^{s}}.
\end{eqnarray*}
When (ii) occurs: we have that $|\sigma_{1}|\sim |\sigma|$ or $|\sigma_{1}|\sim |\sigma_{2}|$.

\noindent
When $|\sigma_{1}|\sim |\sigma|$  is valid, this case can be proved similarly to
$|\sigma|={\rm max}\left\{|\sigma|,|\sigma_{1}|,|\sigma_{2}|\right\}.$

\noindent
When $|\sigma_{1}|\sim |\sigma_{2}|$ is valid, we consider $\supp u_{1}\subset D_{1}$, $\supp u_{1}\subset D_{2}$, $\supp u_{1}\subset D_{3}$,
respectively.

 \noindent When  $\supp u_{1}\subset D_{1}$ which yields that $|k|\leq C$,
since $-j+\frac{3}{2}+j\epsilon\leq s\leq1-\frac{j}{2}-j\epsilon$, by using Lemmas 2.3, 2.5, 2.1,
 we have that
\begin{eqnarray*}
&&\left\|\Lambda^{-1}\partial_{x}(1-\partial_{x}^{2})^{-1}\prod_{j=1}^{2}
(\partial_{x}u_{j})\right\|_{X^{s}}\leq C\left\|\langle k\rangle ^{s-1}\langle \sigma \rangle ^{-\frac{1}{2j}}
(\langle k\rangle\mathscr{F}u_{1}*\langle k\rangle\mathscr{F}u_{2})\right\|_{l_{k}^{2}L_{\tau}^{2}}\nonumber\\&&
\leq C\left\|(\langle k\rangle^{s}\langle\sigma\rangle^{\frac{2j-1}{2j}}\mathscr{F}u_{1})*(\langle k\rangle ^{-s-2j+3}\mathscr{F}u_{2})\right\|_{X_{0,-\frac{1}{2j}}}\nonumber\\
&&\leq C\left\|\left(J^{s}\Lambda ^{\frac{2j-1}{2j}}u_{1}\right)
\left(J^{-s-2j+3}u_{2}\right)\right\|_{X_{0,-\frac{1}{2j}}}\nonumber\\
&&\leq C\|u_{1}\|_{X_{s,\frac{2j-1}{2j}}}\|u_{2}\|_{X_{-s-2j+3,\frac{1}{2}}}\leq C\|u_{1}\|_{X_{s,\frac{2j-1}{2j}}}\|u_{2}\|_{X_{s,\frac{1}{2}}}\leq
C\prod_{j=1}^{2}\|u_{j}\|_{Z^{s}}.
\end{eqnarray*}
When $\supp u_{1}\subset D_{2}$, since $-j+\frac{3}{2}+j\epsilon\leq s\leq1-\frac{j}{2}-j\epsilon$, by using Lemmas 2.3, 2.5, 2.1,
 we have that
\begin{eqnarray*}
&&\left\|\Lambda^{-1}\partial_{x}(1-\partial_{x}^{2})^{-1}\prod_{j=1}^{2}(\partial_{x}u_{j})\right\|_{X^{s}}\leq C\left\|\langle k\rangle ^{s-1}\langle \sigma \rangle ^{-\frac{1}{2j}}
\left[\langle k\rangle\mathscr{F}u_{1}*(\langle k\rangle\mathscr{F}u_{2})\right]\right\|_{l_{k}^{2}L_{\tau}^{2}}\nonumber\\&&
\leq C\left\|\left(J^{(1-2j)(s-1)}\Lambda ^{s}u_{1}\right)\left(J^{-s-2j+3}u_{2}\right)
\right\|_{X_{0,-\frac{1}{2j}}}\nonumber\\
&&\leq C\|u_{1}\|_{X_{(1-2j)(s-1),s}}\|u_{2}\|_{X_{-s-2j+3,\frac{1}{2}}}\leq C\|u_{1}\|_{X_{(1-2j)(s-1),s}}\|u_{2}\|_{X_{s,\frac{1}{2}}}\leq
C\prod_{j=1}^{2}\|u_{j}\|_{Z^{s}}.
\end{eqnarray*}
When $\supp\mathscr{F} u_{1}\subset D_{3}$
which yields that  $\mathscr{F}u_{2}\subset D_{3}$, we consider $\supp\left[\mathscr{F}u_{1}*\mathscr{F}u_{2}\right]\subset D_{1}$,
$\supp\left[\mathscr{F}u_{1}*\mathscr{F}u_{2}\right]\subset D_{2}$, $\supp\left[\mathscr{F}u_{1}*\mathscr{F}u_{2}\right]\subset D_{3}$,
respectively.

\noindent When $\supp\left[\mathscr{F}u_{1}*\mathscr{F}u_{2}\right]\subset D_{1}$, by using the H\"older  inequality  and the  Young inequality and Lemma 2.5, since $-j+\frac{3}{2}+j\epsilon\leq s\leq1-\frac{j}{2}-j\epsilon$,
 we  have that
\begin{eqnarray}
&&\left\|\Lambda^{-1}\partial_{x}(1-\partial_{x}^{2})^{-1}\prod_{j=1}^{2}
(\partial_{x}u_{j})\right\|_{X^{s}}\leq C
\left\|\langle k\rangle ^{s-1}\langle \sigma\rangle ^{-\frac{1}{2j}}
\left[|k|\mathscr{F}u_{1}*|k|\mathscr{F}u_{2}\right]\right\|_{l_{k}^{2}L_{\tau}^{2}}\nonumber\\
&&\leq C\left\|\langle k\rangle ^{s-\frac{1}{2}+\epsilon}\langle
 \sigma \rangle^{-\frac{1}{2j}+\frac{1}{2}+\epsilon}
\left[(|k|\mathscr{F}u_{1})*(|k|\mathscr{F}u_{2})\right]
\right\|_{l_{k}^{\infty}L_{\tau}^{\infty}}\nonumber\\
&&\leq C\left\|\left(\langle k\rangle ^{s+j-\frac{3}{2}+(2j+1)\epsilon}\right)(\langle k\rangle\mathscr{F}u_{1})
*(\langle k\rangle\mathscr{F}u_{2})\right\|_{l_{k}^{\infty}l_{\tau}^{\infty}}\nonumber\\&&
\leq C\left\|\left(\langle k_{1}\rangle ^{s+j+\frac{1}{2}+(2j+1)\epsilon}\mathscr{F}u_{1}
\right)*\mathscr{F}u_{2}\right\|_{l_{k}^{\infty}l_{\tau}^{\infty}}\nonumber\\&&
\leq C
\left\|\langle k\rangle ^{-3s-3j+\frac{9}{2}+(2j+1)\epsilon}\right\|_{l_{k}^{\infty}}
\prod_{j=1}^{2}\|u_{j}\|_{X_{-\frac{s-1}{j}-1,\frac{s-1}{j}+1}}\nonumber\\&&\leq C
\prod_{j=1}^{2}\|u_{j}\|_{X_{-\frac{s-1}{j}-1,\frac{s-1}{j}+1}}\leq C\prod_{j=1}^{2}\|u_{j}\|_{Z^{s}}.
\end{eqnarray}
\noindent When $\supp\left[\mathscr{F}u_{1}*\mathscr{F}u_{2}\right]\subset D_{2}$, by using the H\"older  inequality  and the  Young inequality and Lemma 2.5, since $-j+\frac{3}{2}+j\epsilon\leq s\leq1-\frac{j}{2}-j\epsilon$,
 we  have that
\begin{eqnarray}
&&\left\|\Lambda^{-1}\partial_{x}(1-\partial_{x}^{2})^{-1}\prod_{j=1}^{2}
(\partial_{x}u_{j})\right\|_{X^{s}}\leq C
\left\|\langle k\rangle ^{(1-2j)(s-1)+1}\langle \sigma\rangle ^{s-1}
\left[|k|\mathscr{F}u_{1}*|k|\mathscr{F}u_{2}\right]\right\|_{l_{k}^{2}L_{\tau}^{2}}\nonumber\\
&&\leq C\left\|\langle k\rangle ^{(1-2j)(s-1)+\frac{3}{2}+\epsilon}\langle
 \sigma \rangle^{s-\frac{1}{2}+\epsilon}
\left[(|k|\mathscr{F}u_{1})*(|k|\mathscr{F}u_{2})\right]
\right\|_{l_{k}^{\infty}L_{\tau}^{\infty}}\nonumber\\
&&\leq C\left\|\left(\langle k\rangle ^{s+j-\frac{3}{2}+(2j+1)\epsilon}\right)(\langle k\rangle\mathscr{F}u_{1})
*(\langle k\rangle\mathscr{F}u_{2})\right\|_{l_{k}^{\infty}l_{\tau}^{\infty}}\nonumber\\&&
\leq C\left\|\left(\langle k_{1}\rangle ^{s+j+\frac{1}{2}+(2j+1)\epsilon}\mathscr{F}u_{1}
\right)*\mathscr{F}u_{2}\right\|_{l_{k}^{\infty}l_{\tau}^{\infty}}\nonumber\\&&
\leq C
\left\|\langle k\rangle ^{-3s-3j+\frac{9}{2}+(2j+1)\epsilon}\right\|_{l_{k}^{\infty}}
\prod_{j=1}^{2}\|u_{j}\|_{X_{-\frac{s-1}{j}-1,\frac{s-1}{j}+1}}\nonumber\\&&\leq C
\prod_{j=1}^{2}\|u_{j}\|_{X_{-\frac{s-1}{j}-1,\frac{s-1}{j}+1}}\leq C\prod_{j=1}^{2}\|u_{j}\|_{Z^{s}}.
\end{eqnarray}
\noindent When $\supp\left[\mathscr{F}u_{1}*\mathscr{F}u_{2}\right]\subset D_{3}$, by using the H\"older  inequality  and the  Young inequality and Lemma 2.5, since $-j+\frac{3}{2}+j\epsilon\leq s\leq1-\frac{j}{2}-j\epsilon$,
 we  have that
\begin{eqnarray}
&&\left\|\Lambda^{-1}\partial_{x}(1-\partial_{x}^{2})^{-1}\prod_{j=1}^{2}
(\partial_{x}u_{j})\right\|_{X^{s}}\leq C
\left\|\langle k\rangle ^{-\frac{s-1}{j}}\langle \sigma\rangle ^{\frac{s-1}{j}}
\left[|k|\mathscr{F}u_{1}*|k|\mathscr{F}u_{2}\right]\right\|_{l_{k}^{2}L_{\tau}^{2}}\nonumber\\
&&\leq C\left\|\langle k\rangle ^{-\frac{s-1}{j}+\frac{1}{2}+\epsilon}\langle
 \sigma \rangle^{\frac{s-1}{j}+\frac{1}{2}+\epsilon}
\left[(|k|\mathscr{F}u_{1})*(|k|\mathscr{F}u_{2})\right]
\right\|_{l_{k}^{\infty}L_{\tau}^{\infty}}\nonumber\\
&&\leq C\left\|\left(\langle k\rangle ^{s+j-\frac{3}{2}+(2j+1)\epsilon}\right)(\langle k\rangle\mathscr{F}u_{1})
*(\langle k\rangle\mathscr{F}u_{2})\right\|_{l_{k}^{\infty}l_{\tau}^{\infty}}\nonumber\\&&
\leq C\left\|\left(\langle k_{1}\rangle ^{s+j+\frac{1}{2}+(2j+1)\epsilon}\mathscr{F}u_{1}
\right)*\mathscr{F}u_{2}\right\|_{l_{k}^{\infty}l_{\tau}^{\infty}}\nonumber\\&&
\leq C
\left\|\langle k\rangle ^{-3s-3j+\frac{9}{2}+(2j+1)\epsilon}\right\|_{l_{k}^{\infty}}
\prod_{j=1}^{2}\|u_{j}\|_{X_{-\frac{s-1}{j}-1,\frac{s-1}{j}+1}}\nonumber\\&&\leq C
\prod_{j=1}^{2}\|u_{j}\|_{X_{-\frac{s-1}{j}-1,\frac{s-1}{j}+1}}\leq C\prod_{j=1}^{2}\|u_{j}\|_{Z^{s}}.
\end{eqnarray}

When (c) occurs:  this case can be proved similarly to case (b).

\noindent (3) Region $\Omega_{3}$.
We  consider $|k|\leq |k_{1}|^{-2j}$  and  $|k_{1}|^{-2j}<|k|\leq 1,$
respectively.

\noindent
When $|k|\leq |k_{1}|^{-2j}$, by using Lemma 2.3, since $-j+\frac{3}{2}+j\epsilon\leq s\leq1-\frac{j}{2}-j\epsilon$, we have that
\begin{eqnarray}
&&\left\|\Lambda^{-1}\partial_{x}(1-\partial_{x}^{2})^{-1}\prod_{j=1}^{2}(\partial_{x}u_{j})\right\|_{X^{s}}\leq C
\left\|\langle k\rangle ^{s-1}\langle \sigma\rangle ^{-\frac{1}{2j}}
\left[|k|\mathscr{F}u_{1}*|k|\mathscr{F}u_{2}\right]\right\|_{l_{k}^{2}L_{\tau}^{2}}\nonumber\\
&&\leq C\left\||k_{1}|^{-3j+2}\langle \sigma\rangle ^{-\frac{1}{2j}}
\left[\mathscr{F}u_{1}*\mathscr{F}u_{2}\right]\right\|_{l_{k}^{\infty}L_{\tau}^{2}}\nonumber\\
&&\leq C\left\|\left[\langle k\rangle^{-(2j-1)}\mathscr{F}u_{1}*(\langle k\rangle^{-(j-1)}\mathscr{F}u_{2})\right]
\right\|_{l_{k}^{\infty}L_{\tau}^{2}}\nonumber\\
&&\leq C\|u_{1}\|_{X_{1-2j,0}}\|u_{2}\|_{Y^{1-j}}\leq C\|u_{1}\|_{X_{1-2j,0}}\|u_{2}\|_{Y^{s}}\nonumber\\
&&\leq C\|u_{1}\|_{X_{s,\frac{1}{2j}}}\|u_{2}\|_{Y^{s}}\leq C\prod_{j=1}^{2}\|u_{j}\|_{Z^{s}}.
\end{eqnarray}
Now we consider the case $|k_{1}|^{-2j}\leq|k|\leq 1.$ In this case,
we consider  (a)-(c) of Lemma 2.5, respectively.

\noindent When   (a) occurs:   in this case
$\supp \left[\mathscr{F}u_{1}*\mathscr{F}u_{2}\right]\subset D_{4}$,
by using the H\"older  inequality and the Young inequality,
since $|k|\leq 1$ and $1+\frac{s-1}{j}\geq 0,$  by using Lemma 2.5, since $-j+\frac{3}{2}+j\epsilon\leq s\leq1-\frac{j}{2}-j\epsilon$, we have that
\begin{eqnarray*}
&&\left\|\Lambda^{-1}\partial_{x}(1-\partial_{x}^{2})^{-1}
\prod_{j=1}^{2}(\partial_{x}u_{j})\right\|_{X^{s}}\leq C\left\||k|\langle k\rangle ^{-\frac{s-1}{j}-3}\langle \sigma \rangle ^{\frac{s-1}{j}}\left[(|k|\mathscr{F}u_{1})*(|k|\mathscr{F}u_{2})\right]
\right\|_{l_{k}^{2}L_{\tau}^{2}}\nonumber\\
&&\leq C\left\||k| \left[(|k|^{s}\mathscr{F}u_{1})*(|k|^{s}\mathscr{F}u_{2})\right]\right\|
_{l_{k}^{2}L_{\tau}^{2}}\nonumber\\
&&\leq C\left\|\left[(|k|^{s}\mathscr{F}u_{1})*(|k|^{s}\mathscr{F}u_{2})\right]
\right\|_{l_{k}^{\infty}L_{\tau}^{2}}\leq C\|u_{1}\|_{X_{s,0}}\|u_{2}\|_{Y^{s}}\nonumber\\
&&\leq C\|u_{1}\|_{X_{s,\frac{1}{2j}}}\|u_{2}\|_{Y^{s}}\leq C\prod_{j=1}^{2}\|u_{j}\|_{Z^{s}}.
\end{eqnarray*}
When (b)  occurs:  we   consider the case $|\sigma_{1}|>4{\rm max}\left\{|\sigma|,|\sigma_{2}|\right\}$  and
$|\sigma_{1}|\leq4{\rm max}\left\{|\sigma|,|\sigma_{2}|\right\}$, respectively.

\noindent When $|\sigma_{1}|>4{\rm max}\left\{|\sigma|,|\sigma_{2}|\right\}$
is valid, we consider $\supp \left[\mathscr{F}u_{1}*\mathscr{F}u_{2}\right]\subset D_{4}$ and $\supp \left[\mathscr{F}u_{1}*\mathscr{F}u_{2}\right]\subset D_{5}$, respectively.

\noindent When $\supp \left[\mathscr{F}u_{1}*\mathscr{F}u_{2}\right]\subset D_{4}$, by
 using the H\"older  inequality and the Young inequality and  Lemma 2.5,
 since $|k|\leq 1$, since $-j+\frac{3}{2}+j\epsilon\leq s\leq1-\frac{j}{2}-j\epsilon$, we have that
\begin{eqnarray*}
&&\left\|\Lambda^{-1}\partial_{x}(1-\partial_{x}^{2})^{-1}\prod_{j=1}^{2}
(\partial_{x}u_{j})\right\|_{X^{s}}\leq C\left\||k|\langle k\rangle
^{-\frac{s-1}{j}-3}\langle \sigma \rangle ^{\frac{s-1}{j}}\left[(|k|\mathscr{F}u_{1})*(|k|\mathscr{F}u_{2})\right]
\right\|_{l_{k}^{2}L_{\tau}^{2}}\nonumber\\
&&\leq C\left\||k|^{\frac{1}{2j}} \left[\langle k\rangle^{s}
\langle\sigma\rangle^{\frac{2j-1}{2j}}  \mathscr{F}u_{1}*(\langle k\rangle ^{-s-2j+3}\mathscr{F}u_{2})\right]\right\|_{l_{k}^{2}L_{\tau}^{2}}\nonumber\\
&&\leq C\left\| \left[(\langle k\rangle^{s}\langle\sigma\rangle^{\frac{2j-1}{2j}}  \mathscr{F}u_{1})*(\langle k\rangle ^{-s-2j+3}\mathscr{F}u_{2})\right]\right\|_{l_{k}^{\infty}L_{\tau}^{2}}\nonumber\\
&&\leq C\|u_{1}\|_{X_{s,\frac{2j-1}{2j}}}\|u_{2}\|_{Y^{s}}\leq C\prod_{j=1}^{2}\|u_{j}\|_{X_{s,\frac{2j-1}{2j}}}\leq C\prod_{j=1}^{2}\|u_{j}\|_{Z^{s}}.
\end{eqnarray*}
When $\supp \left[\mathscr{F}u_{1}*\mathscr{F}u_{2}\right]\subset D_{5}$, by using Lemmas
 2.3,  2.5, 2.1,  since $-j+\frac{3}{2}+j\epsilon\leq s\leq1-\frac{j}{2}-j\epsilon$, we have that
\begin{eqnarray*}
&&\left\|\Lambda^{-1}\partial_{x}(1-\partial_{x}^{2})^{-1}\prod_{j=1}^{2}
(\partial_{x}u_{j})\right\|_{X^{s}}
\leq C\left\||k|\langle k\rangle ^{s-2}\langle \sigma \rangle^{-\frac{1}{2j}}
\left[(|k|\mathscr{F}u_{1})*(|k|\mathscr{F}u_{2})\right]
\right\|_{l_{k}^{2}L_{\tau}^{2}}\nonumber\\
&&\leq C\left\|\langle \sigma \rangle^{-\frac{1}{2j}}
\left[
\left(\langle k\rangle ^{s}\langle \sigma \rangle
^{\frac{2j-1}{2j}}\mathscr{F}u_{1}\right)*
\left(\langle k\rangle^{-s-2j+3}\mathscr{F}u_{2}\right)
\right]\right\|_{l_{k}^{2}L_{\tau}^{2}}\nonumber\\
&&\leq C\left\|\left(J^{s}\Lambda^{\frac{2j-1}{2j}}u_{1}
\right)\left(J^{-s-2j+3}u_{2}\right)\right\|_{X_{0,-\frac{1}{2j}}}\nonumber\\
&&\leq C\|J^{s}\Lambda^{\frac{2j-1}{2j}}u_{1}\|_{L_{xt}^{2}}
\|J^{-s-2j+3}u_{2}\|_{X_{0,\frac{1}{2}}}\nonumber\\
&&\leq C\|u_{1}\|_{X_{s,\frac{2j-1}{2j}}}\|u_{2}\|_{X_{s,\frac{1}{2}}}\leq C\prod_{j=1}^{2}\|u_{j}\|_{X_{s,\frac{2j-1}{2j}}}\leq C\prod_{j=1}^{2}\|u_{j}\|_{Z^{s}}.
\end{eqnarray*}
When $|\sigma_{1}|\leq4{\rm max}\left\{|\sigma|,|\sigma_{2}|\right\}$ is valid, we have that
 $|\sigma_{1}|\sim |\sigma|$ or  $|\sigma_{1}|\sim |\sigma_{2}|$.

 \noindent
 When $|\sigma_{1}|\sim|\sigma|$, then this case an be proved similarly to
 case $|\sigma|={\rm max}\left\{|\sigma|,|\sigma_{1}|,|\sigma_{2}|\right\}.$

 \noindent
  When $|\sigma_{1}|\sim|\sigma_{2}|$, we consider $\supp \mathscr{F}u_{1}\subset D_{1},$ $\supp \mathscr{F}u_{1}\subset D_{2},$
  $\supp \mathscr{F}u_{1}\subset D_{3},$ respectively.

 \noindent When $\supp \mathscr{F}u_{1}\subset D_{1},$
  since $-j+\frac{3}{2}+j\epsilon\leq s\leq1-\frac{j}{2}-j\epsilon$, by using Lemma 2.1, 2.5,
 we have that
\begin{eqnarray*}
&&\left\|\Lambda^{-1}\partial_{x}(1-\partial_{x}^{2})^{-1}
\prod_{j=1}^{2}(\partial_{x}u_{j})\right\|_{X^{s}}
\leq C\left\||k|\langle k\rangle^{s-2}\langle \sigma \rangle^{-\frac{1}{2j}}\left[(|k|\mathscr{F}u_{1})*(|k|\mathscr{F}u_{2})\right]
\right\|_{l_{k}^{2}L_{\tau}^{2}}\nonumber\\
&&\leq C\left\|\left(J^{s}\Lambda^{\frac{2j-1}{2j}}u_{1}\right)\left(J^{-s-2j+3}u_{2}\right)
\right\|_{X_{0,-\frac{1}{2j}}}\nonumber\\
&&\leq C\|u_{1}\|_{X_{s,\frac{2j-1}{2j}}}\|u_{2}\|_{X_{-s-2j+3,\frac{1}{2}}}\leq C\|u_{1}\|_{X_{s,\frac{2j-1}{2j}}}\|u_{2}\|_{X_{s,\frac{1}{2}}}\leq C\prod_{j=1}^{2}\|u_{j}\|_{Z^{s}}.
\end{eqnarray*}
When $\supp \mathscr{F}u_{1}\subset D_{2},$
by using Lemma 2.3, since $-j+\frac{3}{2}+j\epsilon\leq s\leq1-\frac{j}{2}-j\epsilon$ and $|\sigma|\leq C|k_{1}|^{2j+1},$ we have that
\begin{eqnarray*}
&&\left\|\Lambda^{-1}(1-\partial_{x}^{2})^{-1}
\partial_{x}\prod_{j=1}^{2}(\partial_{x}u_{j})\right\|_{X^{s}}\leq C\left\||k|\langle k\rangle^{s-2}\langle \sigma \rangle^{-\frac{1}{2j}}\left[(|k|\mathscr{F}u_{1})*(|k|\mathscr{F}u_{2})\right]
\right\|_{l_{k}^{2}L_{\tau}^{2}}\nonumber\\
&&\leq C\left\||k|\langle k\rangle^{s-\frac{3}{2}+\epsilon}\langle \sigma \rangle^{-\frac{1}{2j}+\frac{1}{2}+\epsilon}
\left[(|k|\mathscr{F}u_{1})*(|k|\mathscr{F}u_{2})\right]
\right\|_{l_{k}^{\infty}L_{\tau}^{\infty}}\nonumber\\
&&\leq C\left\|\langle k\rangle^{s-\frac{1}{2}+\epsilon}
\left[(|k|^{j+\frac{1}{2}-\frac{1}{2j}+(2j+1)\epsilon}\mathscr{F}u_{1})*(|k|\mathscr{F}u_{2})\right]
\right\|_{l_{k}^{\infty}L_{\tau}^{\infty}}\nonumber\\
&&\leq C\left\|\left[(\langle k\rangle^{s+1-\frac{1}{2j}+j+(2j+2)\epsilon}\mathscr{F}u_{1})*(\mathscr{F}u_{2})\right]
\right\|_{l_{k}^{\infty}L_{\tau}^{\infty}}\nonumber\\
&&\leq C\left\|\langle k\rangle^{-s-3j+3-\frac{1}{2j}+(2j+2)\epsilon}\right\|_{l_{k}^{\infty}}\prod_{j=1}^{2}\left\|\left(J^{(1-2j)(s-1)}\Lambda^{s}u_{j}\right)
\right\|_{l_{k}^{2}L_{\tau}^{2}}\nonumber\\
&&\leq C\prod_{j=1}^{2}\|u_{j}\|_{X_{(1-2j)(s-1),s}}\leq C\prod_{j=1}^{2}\|u_{j}\|_{Z^{s}};
\end{eqnarray*}
When $\supp \mathscr{F}u_{1}\subset D_{3},$
we consider $\supp\left[\mathscr{F}u_{1}*\mathscr{F}u_{2}\right]\subset D_{1}$, $\supp\left[\mathscr{F}u_{1}*\mathscr{F}u_{2}\right]\subset D_{2}$,
$\supp\left[\mathscr{F}u_{1}*\mathscr{F}u_{2}\right]\subset D_{3}$, respectively.

\noindent When $\supp\left[\mathscr{F}u_{1}*\mathscr{F}u_{2}\right]\subset D_{1}$,
since $-j+\frac{3}{2}+j\epsilon\leq s\leq1-\frac{j}{2}-j\epsilon$,
then, we have that
\begin{eqnarray*}
&&\left\|\Lambda^{-1}(1-\partial_{x}^{2})^{-1}\partial_{x}
\prod_{j=1}^{2}(\partial_{x}u_{j})\right\|_{X^{s}}
\leq C\left\||k|\langle k\rangle^{s-2}\langle
\sigma \rangle^{-\frac{1}{2j}}\left[(|k|\mathscr{F}u_{1})
*(|k|\mathscr{F}u_{2})\right]\right\|_{l_{k}^{2}L_{\tau}^{2}}\nonumber\\
&&\leq C\left\|\langle k\rangle^{s-\frac{1}{2}+\epsilon}\langle \sigma \rangle^{-\frac{1}{2j}+\frac{1}{2}+\epsilon}
\left[(|k|\mathscr{F}u_{1})*(|k|\mathscr{F}u_{2})\right]
\right\|_{l_{k}^{\infty}L_{\tau}^{\infty}}\nonumber\\
&&\leq C\left\|\langle k\rangle^{s+j-\frac{3}{2}+(2j+1)\epsilon}
\left[(\langle k\rangle\mathscr{F}u_{1})*(\langle k\rangle \mathscr{F}u_{2})\right]
\right\|_{l_{k}^{\infty}L_{\tau}^{\infty}}\nonumber\\
&&\leq C\left\|\langle k\rangle^{-4s-4j+6}\right\|_{l_{k}^{\infty}}\prod_{j=1}^{2}\|u_{j}\|_{X_{-\frac{s-1}{j}-1,\frac{s-1}{j}+1}}\nonumber\\
&&\leq C\prod_{j=1}^{2}\|u_{j}\|_{X_{-\frac{s-1}{j}-1,\frac{s-1}{j}+1}}\leq C\prod_{j=1}^{2}\|u_{j}\|_{Z^{s}}.
\end{eqnarray*}
When $\supp\left[\mathscr{F}u_{1}*\mathscr{F}u_{2}\right]\subset D_{2}$,
since $-j+\frac{3}{2}+j\epsilon\leq s\leq1-\frac{j}{2}-j\epsilon$,
 we have that
\begin{eqnarray*}
&&\left\|\Lambda^{-1}(1-\partial_{x}^{2})^{-1}\partial_{x}
\prod_{j=1}^{2}(\partial_{x}u_{j})\right\|_{X^{s}}
\leq C\left\|\langle k\rangle^{(1-2j)(s-1)-1}\langle
\sigma \rangle^{s-1}\left[(|k|\mathscr{F}u_{1})
*(|k|\mathscr{F}u_{2})\right]\right\|_{l_{k}^{2}L_{\tau}^{2}}\nonumber\\
&&\leq C\left\|\langle k\rangle^{(1-2j)(s-1)-\frac{1}{2}+\epsilon}\langle \sigma \rangle^{s-\frac{1}{2}+\epsilon}
\left[(|k|\mathscr{F}u_{1})*(|k|\mathscr{F}u_{2})\right]
\right\|_{l_{k}^{\infty}L_{\tau}^{\infty}}\nonumber\\
&&\leq C\left\|\langle k\rangle^{s+j-\frac{3}{2}+(2j+1)\epsilon}
\left[(\langle k\rangle\mathscr{F}u_{1})*(\langle k\rangle \mathscr{F}u_{2})\right]
\right\|_{l_{k}^{\infty}L_{\tau}^{\infty}}\nonumber\\
&&\leq C\left\|\langle k\rangle^{-4s-4j+6}\right\|_{l_{k}^{\infty}}\prod_{j=1}^{2}\|u_{j}\|_{X_{-\frac{s-1}{j}-1,\frac{s-1}{j}+1}}\nonumber\\
&&\leq C\prod_{j=1}^{2}\|u_{j}\|_{X_{-\frac{s-1}{j}-1,\frac{s-1}{j}+1}}\leq C\prod_{j=1}^{2}\|u_{j}\|_{Z^{s}}.
\end{eqnarray*}
When $\supp\left[\mathscr{F}u_{1}*\mathscr{F}u_{2}\right]\subset D_{3}$,
since $-j+\frac{3}{2}+j\epsilon\leq s\leq1-\frac{j}{2}-j\epsilon$,
 we have that
\begin{eqnarray*}
&&\left\|\Lambda^{-1}(1-\partial_{x}^{2})^{-1}\partial_{x}
\prod_{j=1}^{2}(\partial_{x}u_{j})\right\|_{X^{s}}
\leq C\left\|\langle k\rangle^{-\frac{s-1}{j}-2}\langle
\sigma \rangle^{\frac{s-1}{j}}\left[(|k|\mathscr{F}u_{1})
*(|k|\mathscr{F}u_{2})\right]\right\|_{l_{k}^{2}L_{\tau}^{2}}\nonumber\\
&&\leq C\left\|\langle k\rangle^{-\frac{s-1}{j}-\frac{3}{2}+\epsilon}\langle \sigma \rangle^{\frac{s-1}{j}+\frac{1}{2}+\epsilon}
\left[(|k|\mathscr{F}u_{1})*(|k|\mathscr{F}u_{2})\right]
\right\|_{l_{k}^{\infty}L_{\tau}^{\infty}}\nonumber\\
&&\leq C\left\|\langle k\rangle^{2s+j-3+(2j+2)\epsilon}
\left[(\langle k\rangle\mathscr{F}u_{1})*(\langle k\rangle \mathscr{F}u_{2})\right]
\right\|_{l_{k}^{\infty}L_{\tau}^{\infty}}\nonumber\\
&&\leq C\left\|\langle k\rangle^{-4s-4j+6}\right\|_{l_{k}^{\infty}}\prod_{j=1}^{2}\|u_{j}\|_{X_{-\frac{s-1}{j}-1,\frac{s-1}{j}+1}}\nonumber\\
&&\leq C\prod_{j=1}^{2}\|u_{j}\|_{X_{-\frac{s-1}{j}-1,\frac{s-1}{j}+1}}\leq C\prod_{j=1}^{2}\|u_{j}\|_{Z^{s}}.
\end{eqnarray*}
When (c) occurs: this case can be proved similarly to case (b).

\noindent
(4) Region $\Omega_{4}$.
We consider  (a)-(c) of Lemma 2.5,  respectively.

\noindent When (a) occurs: we consider $|\sigma|>4{\rm max}\left\{|\sigma_{1}|,|\sigma_{2}|\right\}$, $|\sigma|\leq4{\rm max}\left\{|\sigma_{1}|,|\sigma_{2}|\right\}$, respectively.

\noindent When $|\sigma|>4{\rm max}\left\{|\sigma_{1}|,|\sigma_{2}|\right\}$,
 then $\supp \mathscr{F}u_{2}\subset D_{1}\cup D_{2}$ and  $
 \supp \left[\mathscr{F}u_{1}*\mathscr{F}u_{2}\right]\subset D_{2}.$
In this case, by using Lemmas 2.5,  2.3, since $-j+\frac{3}{2}+j\epsilon\leq s\leq1-\frac{j}{2}-j\epsilon$, we have that
\begin{eqnarray*}
&&\left\|\Lambda^{-1}\partial_{x}(1-\partial_{x}^{2})^{-1}
\prod_{j=1}^{2}(\partial_{x}u_{j})\right\|_{X^{s}}\leq C
\left\|\langle k\rangle^{-(2j-1)(s-1)-1}\langle \sigma \rangle^{s-1}\left[(|k|\mathscr{F}u_{1})*(|k|\mathscr{F}u_{2})\right]
\right\|_{l_{k}^{2}L_{\tau}^{2}}\nonumber\\
&&\leq C\left\|(J^{s}u_{1})(J^{s}u_{2})\right\|_{L_{xt}^{2}}\nonumber\\&&
\leq C\|u_{1}\|_{X_{s,\frac{1}{2j}}}\|u_{2}\|_{X_{s,\frac{1}{2}}}
\leq C\prod_{j=1}^{2}\|u_{j}\|_{Z^{s}}.
\end{eqnarray*}
When $|\sigma|\leq4{\rm max}\left\{|\sigma_{1}|,|\sigma_{2}|\right\}$, we have
 that $|\sigma|\sim |\sigma_{1}|$ or $|\sigma|\sim |\sigma_{2}|$.

\noindent
 When  $|\sigma|\sim |\sigma_{1}|$, we consider
 $\supp \left[\mathscr{F}u_{1}*\mathscr{F}u_{2}\right]\subset D_{1},$ $\supp \left[\mathscr{F}u_{1}*\mathscr{F}u_{2}\right]\subset D_{2},$
 $\supp \left[\mathscr{F}u_{1}*\mathscr{F}u_{2}\right]\subset D_{3},$ respectively.

\noindent When
 $\supp \left[\mathscr{F}u_{1}*\mathscr{F}u_{2}\right]\subset D_{1},$
 by using Lemmas 2.3, 2.1, since $-j+\frac{3}{2}+j\epsilon\leq s\leq1-\frac{j}{2}-j\epsilon$, we have that
\begin{eqnarray*}
&&\left\|\Lambda^{-1}(1-\partial_{x}^{2})^{-1}\partial_{x}\prod_{j=1}^{2}
(\partial_{x}u_{j})\right\|_{X^{s}}\leq C\left\|\langle k\rangle^{s-1}\langle \sigma \rangle^{-\frac{1}{2j}}
\left[(|k|\mathscr{F}u_{1})*(|k|\mathscr{F}u_{2})\right]\right\|_{l_{k}^{2}L_{\tau}^{2}}
 \nonumber\\&&\leq\left\|(J^{s}u_{1})(J^{s}u_{2})\right\|_{L_{xt}^{2}}
\leq C\|u_{1}\|_{X_{s,\frac{1}{2j}}}\|u_{2}\|_{X_{s,\frac{1}{2}}}\leq C\prod_{j=1}^{2}\|u_{j}\|_{Z^{s}}.
\end{eqnarray*}
When $\supp \left[\mathscr{F}u_{1}*\mathscr{F}u_{2}\right]\subset D_{2},$
 by using Lemmas 2.5, 2.1, since $-j+\frac{3}{2}+j\epsilon\leq s\leq1-\frac{j}{2}-j\epsilon$, we have that
\begin{eqnarray*}
&&\left\|\Lambda^{-1}(1-\partial_{x}^{2})^{-1}\partial_{x}(\prod_{j=1}^{2}u_{j})
\right\|_{X^{s}}\leq C\left\|\langle k\rangle^{-(2j-1)(s-1)-1}\langle \sigma \rangle^{s-1}
\left[(|k|\mathscr{F}u_{1})*(|k|\mathscr{F}u_{2})\right]\right\|_{l_{k}^{2}L_{\tau}^{2}} \nonumber\\&&\leq\left\|(J^{s}u_{1})(J^{s}u_{2})\right\|_{L_{xt}^{2}}
\leq C\|u_{1}\|_{X_{s,\frac{1}{2j}}}\|u_{2}\|_{X_{s,\frac{1}{2}}}\leq C\prod_{j=1}^{2}\|u_{j}\|_{Z^{s}}.
\end{eqnarray*}
When $\supp \left[\mathscr{F}u_{1}*\mathscr{F}u_{2}\right]\subset D_{3},$
 by using Lemma 2.5, the Young inequality, since $-j+\frac{3}{2}+j\epsilon\leq s\leq1-\frac{j}{2}-j\epsilon$, we have that
\begin{eqnarray*}
&&\left\|\Lambda^{-1}\partial_{x}(1-\partial_{x}^{2})^{-1}\prod_{j=1}^{2}(\partial_{x}u_{j})\right\|_{X^{s}}\leq C\left\|\langle k\rangle^{-\frac{s-1}{j}-2}\langle \sigma \rangle^{\frac{s-1}{j}}
\left[|k|\mathscr{F}u_{1}*|k|\mathscr{F}u_{2}\right]\right\|_{l_{k}^{2}L_{\tau}^{2}} \nonumber\\&&\leq\left\|(J^{-\frac{s-1}{j}-1}\Lambda^{\frac{s-1}{j}+1}u_{1})(J^{-2j}u_{2})\right\|_{L_{xt}^{2}}\nonumber\\&&\leq C\|u_{1}\|_{X_{-\frac{s-1}{j}-1,\frac{s-1}{j}+1}}\|J^{-2j}u_{2}\|_{l_{k}^{1}L_{\tau}^{1}}\leq C\|u_{1}\|_{X_{-\frac{s-1}{j}-1,\frac{s-1}{j}+1}}\|u_{2}\|_{Y^{s}}\leq C\prod_{j=1}^{2}\|u_{j}\|_{Z^{s}}.
\end{eqnarray*}
When $|\sigma|\sim |\sigma_{2}|$, we consider $\supp \left[\mathscr{F}u_{1}*\mathscr{F}u_{2}\right]\subset D_{1},$ $\supp \left[\mathscr{F}u_{1}*\mathscr{F}u_{2}\right]\subset D_{2},$
$\supp \left[\mathscr{F}u_{1}*\mathscr{F}u_{2}\right]\subset D_{3},$ respectively.

\noindent When $\supp \left[\mathscr{F}u_{1}*\mathscr{F}u_{2}\right]\subset D_{1},$
 $1\leq |k_{1}|\leq C,$ by using Lemmas 2.3, 2.1,  since $-j+\frac{3}{2}+j\epsilon\leq s\leq1-\frac{j}{2}-j\epsilon$,  we have that
\begin{eqnarray*}
&&\left\|\Lambda^{-1}(1-\partial_{x}^{2})^{-1}\partial_{x}\prod_{j=1}^{2}
(\partial_{x}u_{j})\right\|_{X^{s}}\leq C\left\|\langle k\rangle^{s-1}\langle \sigma \rangle^{-\frac{1}{2j}}
\left[(|k|\mathscr{F}u_{1})*(|k|\mathscr{F}u_{2})\right]\right\|_{l_{k}^{2}L_{\tau}^{2}}
 \nonumber\\&&\leq\left\|(J^{s}u_{1})(J^{s}u_{2})\right\|_{L_{xt}^{2}}
\leq C\|u_{1}\|_{X_{s,\frac{1}{2j}}}\|u_{2}\|_{X_{s,\frac{1}{2}}}\leq C\prod_{j=1}^{2}\|u_{j}\|_{Z^{s}}.
\end{eqnarray*}
When $\supp \left[\mathscr{F}u_{1}*\mathscr{F}u_{2}\right]\subset D_{2},$
 by using Lemma 2.5  and  the Young inequality, since $-j+\frac{3}{2}+j\epsilon\leq s\leq1-\frac{j}{2}-j\epsilon$,  we have that
\begin{eqnarray*}
&&\left\|\Lambda^{-1}(1-\partial_{x}^{2})^{-1}\partial_{x}(\prod_{j=1}^{2}u_{j})
\right\|_{X^{s}}\leq C\left\|\langle k\rangle^{-(2j-1)(s-1)-1}\langle \sigma \rangle^{s-1}
\left[(|k|\mathscr{F}u_{1})*(|k|\mathscr{F}u_{2})\right]\right\|_{l_{k}^{2}L_{\tau}^{2}} \nonumber\\&&\leq\left\|(J^{s}u_{1})(J^{s}u_{2})\right\|_{L_{xt}^{2}}
\leq C\|u_{1}\|_{X_{s,\frac{1}{2j}}}\|u_{2}\|_{X_{s,\frac{1}{2}}}\leq C\prod_{j=1}^{2}\|u_{j}\|_{Z^{s}}.
\end{eqnarray*}
When $\supp \left[\mathscr{F}u_{1}*\mathscr{F}u_{2}\right]\subset D_{3},$
 by using  using Lemma 2.5  and the Young inequality, since $-j+\frac{3}{2}+j\epsilon\leq s\leq1-\frac{j}{2}-j\epsilon$, we have that
\begin{eqnarray*}
&&\left\|\Lambda^{-1}\partial_{x}(1-\partial_{x}^{2})^{-1}
\prod_{j=1}^{2}(\partial_{x}u_{j})\right\|_{X^{s}}\leq C\left\|\langle k\rangle^{-\frac{s-1}{j}-2}\langle \sigma \rangle^{\frac{s-1}{j}}
\left[(|k|\mathscr{F}u_{1})*(|k|\mathscr{F}u_{2})\right]
\right\|_{l_{k}^{2}L_{\tau}^{2}} \nonumber\\&&\leq\left\|(J^{-2j}u_{1})(J^{-\frac{s-1}{j}-1}
\Lambda^{\frac{s-1}{j}+1}u_{2})\right\|_{L_{xt}^{2}}\nonumber\\&&\leq C\|J^{-2j}u_{1}\|_{l_{k}^{1}L_{\tau}^{1}}\|u_{2}\|_{X_{-\frac{s-1}{j}-1,\frac{s-1}{j}+1}}\leq C\|u_{1}\|_{Y^{s}}\|u_{2}\|_{X_{-\frac{s-1}{j}-1,\frac{s-1}{j}+1}}\leq C\prod_{j=1}^{2}\|u_{j}\|_{Z^{s}}.
\end{eqnarray*}
(b): $|\sigma_{1}|={\rm max}\left\{|\sigma|,|\sigma_{1}|,|\sigma_{2}|\right\}.$
If $|\sigma_{1}|>4{\rm max}\left\{|\sigma|,|\sigma_{2}|\right\}$,
 then $\supp \mathscr{F}u_{1}\subset D_{3}$.

   \noindent By using Lemmas 2.3, 2.5,  2.1,   since $-j+\frac{3}{2}+j\epsilon\leq s\leq1-\frac{j}{2}-j\epsilon$,  we have that
\begin{eqnarray*}
&&\left\|\Lambda^{-1}\partial_{x}(1-\partial_{x}^{2})^{-1}\prod_{j=1}^{2}(\partial_{x}u_{j})\right\|_{X^{s}}\leq C\left\|\langle k\rangle^{s}\langle \sigma \rangle^{-\frac{1}{2j}}
\left[(|k_{1}|\mathscr{F}u_{1})*\mathscr{F}u_{2}\right]\right\|_{l_{k}^{2}L_{\tau}^{2}}
\nonumber\\
&&\leq C\left\|(J^{-\frac{s-1}{j}-1}\Lambda^{\frac{s-1}{j}+1}u_{1})
(J^{-s-(2j-3)}u_{2})\right\|_{X_{0,-\frac{1}{2j}}}
\nonumber\\&&\leq C\|u_{1}\|_{X_{-\frac{s-1}{j}-1,\frac{s-1}{j}+1}}\|u_{2}\|_{X_{s,\frac{1}{2}}}
\leq C\prod_{j=1}^{2}\|u_{j}\|_{Z^{s}}.
\end{eqnarray*}
 When $|\sigma_{1}|\leq4{\rm max}\left\{|\sigma|,|\sigma_{2}|\right\}$
we have that $|\sigma_{1}|\sim |\sigma|$ or $|\sigma_{1}|\sim |\sigma_{2}|$.

\noindent
 When  $|\sigma_{1}|\sim |\sigma|$, this case can be proved similarly
 to case $|\sigma|={\rm max}\left\{|\sigma|,|\sigma_{1}|,|\sigma_{2}|\right\}.$

\noindent When $|\sigma_{1}|\sim |\sigma_{2}|$,
we consider $\supp \mathscr{F}u_{2}\subset D_{2}$ and  $\supp \mathscr{F}u_{2}\subset D_{3},$
respectively.

\noindent When $\supp \mathscr{F}u_{2}\subset D_{2},$
 by using Lemmas 2.3, 2.5,  2.1, since $-j+\frac{3}{2}+j\epsilon\leq s\leq1-\frac{j}{2}-j\epsilon$,  we have that
\begin{eqnarray*}
&&\left\|\Lambda^{-1}\partial_{x}(1-\partial_{x}^{2})^{-1}\prod_{j=1}^{2}
(\partial_{x}u_{j})\right\|_{X^{s}}\leq C\left\|\langle k\rangle^{s}\langle \sigma \rangle^{-\frac{1}{2j}}
\left[(|k|\mathscr{F}u_{1})*\mathscr{F}u_{2}\right]\right\|_{l_{k}^{2}L_{\tau}^{2}}
 \nonumber\\&&\leq\left\|(J^{-s-(2j-2)}u_{1})(J^{(1-2j)(s-1)}\Lambda^{s}u_{2})\right\|_{X_{0,-\frac{1}{2j}}}
\nonumber\\&&\leq C\|u_{1}\|_{Z^{s}{(D_{2}\cup D_{3})}}\|u_{2}\|_{X_{(1-2j)(s-1),s}}\leq C\prod_{j=1}^{2}\|u_{j}\|_{Z^{s}}.
\end{eqnarray*}
When $\supp \mathscr{F}u_{2}\subset D_{3},$  we consider $\supp\left[\mathscr{F}u_{1}*\mathscr{F}u_{2}\right]\subset D_{1}$, $\supp\left[\mathscr{F}u_{1}*\mathscr{F}u_{2}\right]\subset D_{2}$,
$\supp\left[\mathscr{F}u_{1}*\mathscr{F}u_{2}\right]\subset D_{3}$, respectively.

\noindent
When $\supp\left[\mathscr{F}u_{1}*\mathscr{F}u_{2}\right]\subset D_{1}$, since $-j+\frac{3}{2}+j\epsilon\leq s\leq1-\frac{j}{2}-j\epsilon$, by using Lemmas 2.3, 2.5,  2.1, we have that
\begin{eqnarray*}
&&\left\|\Lambda^{-1}\partial_{x}(1-\partial_{x}^{2})^{-1}\prod_{j=1}^{2}(\partial_{x}u_{j})\right\|_{X^{s}}\leq C\left\|\langle k\rangle^{s}\langle \sigma \rangle^{-\frac{1}{2j}}
\left[(|k_{1}|\mathscr{F}u_{1})*\mathscr{F}u_{2}\right]\right\|_{l_{k}^{2}L_{\tau}^{2}}
\nonumber\\
&&\leq C\left\|\langle k\rangle^{s+\frac{1}{2}+\epsilon}\langle \sigma \rangle^{-\frac{1}{2j}+\frac{1}{2}+\epsilon}
\left[(|k_{1}|\mathscr{F}u_{1})*\mathscr{F}u_{2}\right]\right\|_{l_{k}^{\infty}L_{\tau}^{\infty}}
\nonumber\\
&&\leq C\left\|\langle k\rangle^{s+j-\frac{1}{2}+(2j+1)\epsilon}
\left[(|k_{1}|\mathscr{F}u_{1})*\mathscr{F}u_{2}\right]\right\|_{l_{k}^{\infty}L_{\tau}^{\infty}}
\nonumber\\
&&\leq C\left\|
\left[(\langle k\rangle^{-2s-2j+3}\mathscr{F}u_{1})*(\langle k\rangle^{-s-j+\frac{3}{2}+(2j+1)\epsilon}\langle k\rangle^{-\frac{s-1}{j}-1}\langle \sigma\rangle^{\frac{s-1}{j}+1}\mathscr{F}u_{2})\right]\right\|_{l_{k}^{\infty}L_{\tau}^{\infty}}
\nonumber\\
&&\leq  C\left\|
\left[(\langle k\rangle^{-\frac{s-1}{j}-1}\langle \sigma\rangle^{\frac{s-1}{j}+1}\mathscr{F}u_{1})*(\langle k\rangle^{-\frac{s-1}{j}-1}\langle \sigma\rangle^{\frac{s-1}{j}+1}\mathscr{F}u_{2})\right]\right\|_{l_{k}^{\infty}L_{\tau}^{\infty}}
\nonumber\\
&&\leq C\|u_{1}\|_{X_{-\frac{s-1}{j}-1,\frac{s-1}{j}+1}}\|u_{2}\|_{X_{-\frac{s-1}{j}-1,\frac{s-1}{j}+1}}
\leq C\prod_{j=1}^{2}\|u_{j}\|_{Z^{s}}.
\end{eqnarray*}
When $\supp\left[\mathscr{F}u_{1}*\mathscr{F}u_{2}\right]\subset D_{2}$, by using Lemmas 2.3, 2.5,  2.1, since $-j+\frac{3}{2}+j\epsilon\leq s\leq1-\frac{j}{2}-j\epsilon$, we have that
\begin{eqnarray*}
&&\left\|\Lambda^{-1}\partial_{x}(1-\partial_{x}^{2})^{-1}\prod_{j=1}^{2}(\partial_{x}u_{j})\right\|_{X^{s}}\leq C\left\|\langle k\rangle^{(1-2j)(s-1)-1}\langle \sigma \rangle^{s-1}
\left[(|k_{1}|\mathscr{F}u_{1})*\mathscr{F}u_{2}\right]\right\|_{l_{k}^{2}L_{\tau}^{2}}
\nonumber\\
&&\leq C\left\|\langle k\rangle^{(1-2j)(s-1)-\frac{1}{2}+\epsilon}\langle \sigma \rangle^{s-\frac{1}{2}+\epsilon}
\left[(|k_{1}|\mathscr{F}u_{1})*\mathscr{F}u_{2}\right]\right\|_{l_{k}^{\infty}L_{\tau}^{\infty}}
\nonumber\\
&&\leq C\left\|\langle k\rangle^{s+j-\frac{3}{2}+(2j+1)\epsilon}
\left[(|k_{1}|\mathscr{F}u_{1})*\mathscr{F}u_{2}\right]\right\|_{l_{k}^{\infty}L_{\tau}^{\infty}}
\nonumber\\
&&\leq C\left\|
\left[(\langle k\rangle^{-2s-2j+3}\mathscr{F}u_{1})*(\langle k\rangle^{-s-j+\frac{1}{2}+(2j+1)\epsilon}\langle k\rangle^{-\frac{s-1}{j}-1}\langle \sigma\rangle^{\frac{s-1}{j}+1}\mathscr{F}u_{2})\right]\right\|_{l_{k}^{\infty}L_{\tau}^{\infty}}
\nonumber\\
&&\leq  C\left\|
\left[(\langle k\rangle^{-\frac{s-1}{j}-1}\langle \sigma\rangle^{\frac{s-1}{j}+1}\mathscr{F}u_{1})*(\langle k\rangle^{-\frac{s-1}{j}-1}\langle \sigma\rangle^{\frac{s-1}{j}+1}\mathscr{F}u_{2})\right]\right\|_{l_{k}^{\infty}L_{\tau}^{\infty}}
\nonumber\\
&&\leq C\|u_{1}\|_{X_{-\frac{s-1}{j}-1,\frac{s-1}{j}+1}}\|u_{2}\|_{X_{-\frac{s-1}{j}-1,\frac{s-1}{j}+1}}
\leq C\prod_{j=1}^{2}\|u_{j}\|_{Z^{s}}.
\end{eqnarray*}
When $\supp\left[\mathscr{F}u_{1}*\mathscr{F}u_{2}\right]\subset D_{3}$, by using Lemmas 2.3, 2.5,  2.1, since $-j+\frac{3}{2}+j\epsilon\leq s\leq1-\frac{j}{2}-j\epsilon$, we have that
\begin{eqnarray*}
&&\left\|\Lambda^{-1}\partial_{x}(1-\partial_{x}^{2})^{-1}\prod_{j=1}^{2}(\partial_{x}u_{j})\right\|_{X^{s}}\leq C\left\|\langle k\rangle^{-\frac{s-1}{j}-2}\langle \sigma \rangle^{\frac{s-1}{j}}
\left[(|k_{1}|\mathscr{F}u_{1})*\mathscr{F}u_{2}\right]\right\|_{l_{k}^{2}L_{\tau}^{2}}
\nonumber\\
&&\leq C\left\|\langle k\rangle^{-\frac{s-1}{j}-\frac{3}{2}+\epsilon}\langle \sigma \rangle^{\frac{s-1}{j}+\frac{1}{2}+\epsilon}
\left[(|k_{1}|\mathscr{F}u_{1})*\mathscr{F}u_{2}\right]\right\|_{l_{k}^{\infty}L_{\tau}^{\infty}}
\nonumber\\
&&\leq C\left\|\langle k\rangle^{2s+j-3+(2j+2)\epsilon}
\left[(|k_{1}|\mathscr{F}u_{1})*\mathscr{F}u_{2}\right]\right\|_{l_{k}^{\infty}L_{\tau}^{\infty}}
\nonumber\\
&&\leq C\left\|
\left[(\langle k\rangle^{-2s-2j+3}\mathscr{F}u_{1})*(\langle k\rangle^{-j-1+(2j+2)\epsilon}\langle k\rangle^{-\frac{s-1}{j}-1}\langle \sigma\rangle^{\frac{s-1}{j}+1}\mathscr{F}u_{2})\right]\right\|_{l_{k}^{\infty}L_{\tau}^{\infty}}
\nonumber\\
&&\leq  C\left\|
\left[(\langle k\rangle^{-\frac{s-1}{j}-1}\langle \sigma\rangle^{\frac{s-1}{j}+1}\mathscr{F}u_{1})*(\langle k\rangle^{-\frac{s-1}{j}-1}\langle \sigma\rangle^{\frac{s-1}{j}+1}\mathscr{F}u_{2})\right]\right\|_{l_{k}^{\infty}L_{\tau}^{\infty}}
\nonumber\\
&&\leq C\|u_{1}\|_{X_{-\frac{s-1}{j}-1,\frac{s-1}{j}+1}}\|u_{2}\|_{X_{-\frac{s-1}{j}-1,\frac{s-1}{j}+1}}
\leq C\prod_{j=1}^{2}\|u_{j}\|_{Z^{s}}.
\end{eqnarray*}
When (c) occurs: we consider case $|\sigma_{2}|>4{\rm max}\left\{
|\sigma|,|\sigma_{1}|\right\},$ case $|\sigma_{2}|\leq4{\rm max}\left\{
|\sigma|,|\sigma_{1}|\right\},$ respectively.

\noindent When $|\sigma_{2}|>4{\rm max}\left\{
|\sigma|,|\sigma_{1}|\right\},$ obviously, $\supp \mathscr{F}u_{2}\subset D_{2}\cup D_{3}$.

\noindent When $\supp \mathscr{F}u_{2}\subset D_{2}$, by using Lemmas 2.3, 2.5, 2.1, since $-j+\frac{3}{2}+j\epsilon\leq s\leq1-\frac{j}{2}-j\epsilon$,  we have that
\begin{eqnarray*}
&&\left\|\Lambda^{-1}\partial_{x}(1-\partial_{x}^{2})^{-1}\prod_{j=1}^{2}
(\partial_{x}u_{j})\right\|_{X^{s}}\leq C\left\|\langle k\rangle^{s}\langle \sigma \rangle^{-\frac{1}{2j}}
\left[(|k|\mathscr{F}u_{1})*\mathscr{F}u_{2}\right]\right\|_{l_{k}^{2}L_{\tau}^{2}}
 \nonumber\\&&\leq\left\|(J^{-s-(2j-2)}u_{1})(J^{(1-2j)(s-1)}\Lambda^{s}u_{2})\right\|_{X_{0,-\frac{1}{2j}}}
\nonumber\\&&\leq C\|u_{1}\|_{Z^{s}}\|u_{2}\|_{X_{(1-2j)(s-1),s}}\leq C\prod_{j=1}^{2}\|u_{j}\|_{Z^{s}}.
\end{eqnarray*}
When $\supp \mathscr{F}u_{2}\subset D_{3}$ and $\supp\mathscr{F}u_{1}\subset D_{1}\cup D_{2}$, by using Lemmas 2.3, 2.5, 2.1, since $-j+\frac{3}{2}+j\epsilon\leq s\leq1-\frac{j}{2}-j\epsilon$,
we have that
\begin{eqnarray*}
&&\left\|\Lambda^{-1}\partial_{x}(1-\partial_{x}^{2})^{-1}\prod_{j=1}^{2}(\partial_{x}u_{j})\right\|_{X^{s}}\leq C\left\|\langle k\rangle^{s}\langle \sigma \rangle^{-\frac{1}{2j}}
\left[(|k_{1}|\mathscr{F}u_{1})*\mathscr{F}u_{2}\right]\right\|_{l_{k}^{2}L_{\tau}^{2}}
\nonumber\\
&&\leq C\left\|(J^{-s-(2j-3)}u_{1})(J^{-\frac{s-1}{j}-1}\Lambda^{\frac{s-1}{j}+1}u_{2})
\right\|_{X_{0,-\frac{1}{2j}}}
\nonumber\\&&\leq C\|u_{1}\|_{X_{s,\frac{1}{2}}}\|u_{2}\|_{X_{-\frac{s-1}{j}-1,\frac{s-1}{j}+1}}
\leq C\prod_{j=1}^{2}\|u_{j}\|_{Z^{s}}.
\end{eqnarray*}
When $\supp \mathscr{F}u_{2}\subset D_{3}$ and $\supp\mathscr{F}u_{1}\subset D_{3}$, this case can be proved similarly to case $\supp \mathscr{F}u_{2}\subset D_{3}$ of $|\sigma_{1}|\sim|\sigma_{2}|$ in case (b).

\noindent Case $|\sigma_{2}|\leq4{\rm max}\left\{
|\sigma|,|\sigma_{1}|\right\}$ can be proved similarly to $|\sigma_{1}|\leq4{\rm max}\left\{
|\sigma|,|\sigma_{2}|\right\}$.

\noindent(5) In region $\Omega_{5}$. In this region, we consider
cases $|k_{1}|\leq |k|^{-2j}$ and $|k|^{-2j}< |k_{1}|\leq 1,$ respectively.

\noindent
When $|k_{1}|\leq |k|^{-2j}$, by using the Cauchy-Schwartz inequality and Young inequality as well as Lemma 2.3, since $-j+\frac{3}{2}+j\epsilon\leq s\leq1-\frac{j}{2}-j\epsilon$, we have that
\begin{eqnarray}
&&\left\|\Lambda^{-1}\partial_{x}(1-\partial_{x}^{2})^{-1}
\prod_{j=1}^{2}(\partial_{x}u_{j})\right\|_{X^{s}}\leq C\left\|\langle k\rangle^{s}\langle \sigma\rangle ^{-\frac{1}{2j}}
\left[(|k|\mathscr{F}u_{1})*\mathscr{F}u_{2}\right]\right\|_{l_{k}^{2}L_{\tau}^{2}}\nonumber\\
&&\leq C\left\|\langle k\rangle ^{-2j}\left[\mathscr{F}u_{1}*(\langle k\rangle^{s}\mathscr{F}u_{2})\right]\right\|_{l_{k}^{2}L_{\tau}^{2}}\nonumber\\
&&\leq C\left\|\left[\mathscr{F}u_{1}*(\langle k\rangle^{s}\mathscr{F}u_{2})\right]\right\|_{l_{k}^{\infty}L_{\tau}^{2}}\nonumber\\
&&\leq C\left\|\mathscr{F}u_{1}\right\|_{l_{k}^{2}l_{\tau}^{2}}\|u_{2}\|_{Y^{s}}
\leq C\|\mathscr{F}u_{1}\|_{l_{k}^{2}L_{\tau}^{2}}\|u_{2}\|_{Y^{s}}
\leq C\prod_{j=1}^{2}\|u_{j}\|_{Z^{s}}.
\end{eqnarray}
When  $|k|^{-2j}\leq |k_{1}|\leq 1,$
we consider  (a)-(c) of Lemma 2.5,  respectively.

\noindent
When (a) occurs:
 by using the Young inequality and Cauchy-Schwartz inequality, since $-j+\frac{3}{2}+j\epsilon\leq s\leq1-\frac{j}{2}-j\epsilon$, we have that
\begin{eqnarray}
&&\left\|\Lambda^{-1}\partial_{x}(1-\partial_{x}^{2})^{-1}
\prod_{j=1}^{2}(\partial_{x}u_{j})\right\|_{X^{s}}\leq C
\left\|\langle k\rangle^{s}\langle \sigma\rangle ^{-\frac{1}{2j}}
\left[(|k|\mathscr{F}u_{1})*\mathscr{F}u_{2}\right]
\right\|_{l_{k}^{2}L_{\tau}^{2}}\nonumber\\
&&\leq C\left\|\left[(|k|^{-\frac{1}{2j}}
\mathscr{F}u_{1})*(\langle k\rangle ^{s-1}\mathscr{F}u_{2})\right]
\right\|_{l_{k}^{2}L_{\tau}^{2}}\nonumber\\
&&\leq C\left\||k|^{-\frac{1}{2j}}\mathscr{F}u_{1}\right\|_{l_{k}^{1}l_{\tau}^{2}}\|u_{2}\|_{Y^{s}}
\leq C\|\mathscr{F}u_{1}\|_{l_{k}^{2}L_{\tau}^{2}}\|u_{2}\|_{Y^{s}}
\leq C\prod_{j=1}^{2}\|u_{j}\|_{Z^{s}}.
\end{eqnarray}
When (b) occurs:
 by using the Young inequality and Cauchy-Schwartz
 inequality, since $-j+\frac{3}{2}+j\epsilon\leq s\leq1-\frac{j}{2}-j\epsilon$, we have that
\begin{eqnarray}
&&\left\|\Lambda^{-1}\partial_{x}(1-\partial_{x}^{2})^{-1}\prod_{j=1}^{2}(\partial_{x}u_{j})\right\|_{X^{s}}\leq C\left\|\langle k\rangle ^{s}\langle \sigma\rangle ^{-\frac{1}{2j}}
\left[(|k|\mathscr{F}u_{1})*\mathscr{F}u_{2}\right]\right\|_{l_{k}^{2}L_{\tau}^{2}}\nonumber\\
&&\leq C\left\|\left[(|k|^{1-\frac{1}{2j}}\langle \sigma \rangle ^{\frac{1}{2j}}\mathscr{F}u_{1})*(\langle k\rangle ^{s-1}\mathscr{F}u_{2})\right]
\right\|_{l_{k}^{2}L_{\tau}^{2}}\nonumber\\
&&\leq C\left\||k|^{1-\frac{1}{2j}}\langle \sigma \rangle ^{\frac{1}{2j}}\mathscr{F}u_{1}\right\|_{l_{k}^{1}l_{\tau}^{2}}\|u_{2}\|_{Y^{s}}\nonumber\\&&
\leq C\|\langle \sigma \rangle ^{\frac{1}{2j}}\mathscr{F}u_{1}\|_{l_{k}^{2}L_{\tau}^{2}}\|u_{2}\|_{Y^{s}}\nonumber\\&&
\leq C\|u_{1}\|_{X_{0,\frac{1}{2j}}}\|u_{2}\|_{Y^{s}}
\leq C\prod_{j=1}^{2}\|u_{j}\|_{Z^{s}}.
\end{eqnarray}
When (c)  occurs:
by using the Young inequality and Cauchy-Schwartz inequality, since $-j+\frac{3}{2}+j\epsilon\leq s\leq1-\frac{j}{2}-j\epsilon$, we have that
\begin{eqnarray}
&&\left\|\Lambda^{-1}\partial_{x}(1-\partial_{x}^{2})^{-1}
\prod_{j=1}^{2}(\partial_{x}u_{j})\right\|_{X^{s}}\leq C\left\|\langle k\rangle ^{s}\langle \sigma\rangle ^{-\frac{1}{2j}}
\left[\mathscr{F}u_{1}*\mathscr{F}u_{2}\right]\right\|_{l_{k}^{2}L_{\tau}^{2}}\nonumber\\
&&\leq C\left\|\left[(|k|^{1-\frac{1}{2j}}\mathscr{F}u_{1})*(\langle k\rangle ^{s}\langle \sigma \rangle ^{\frac{1}{2j}}\mathscr{F}u_{2})\right]
\right\|_{l_{k}^{2}L_{\tau}^{2}}\nonumber\\
&&\leq C\left\||k|^{1-\frac{1}{2j}}\mathscr{F}u_{1}\right\|_{l_{k}^{1}l_{\tau}^{1}}
\|u_{2}\|_{X_{s,\frac{1}{2j}}}\nonumber\\&&
\leq C\|\mathscr{F}u_{1}\|_{l_{k}^{2}L_{\tau}^{1}}\|u_{2}\|_{X_{s,\frac{1}{2j}}}\nonumber\\&&
\leq C\|u_{1}\|_{Y^{s}}\|u_{2}\|_{X_{s,\frac{1}{2j}}}
\leq C\prod_{j=1}^{2}\|u_{j}\|_{Z^{s}}.
\end{eqnarray}
(6)In region $\Omega_{6}$.
This case can be proved similarly to $\Omega_{4}$.

\noindent (7)In region $\Omega_{7}$.
This case can be proved similarly to $\Omega_{5}$.

\noindent (8)In region $\Omega_{8}$.
We consider  (a)-(c) of Lemma 2.5,  respectively.

\noindent
When (a) occurs: we have that $\supp\left( \mathscr{F}u_{1}*\mathscr{F}u_{2}\right) \subset D_{3}$.

\noindent
If $|\sigma|>4{\rm max}\left\{|\sigma_{1}|,|\sigma_{2}|\right\}$ and $\supp \mathscr{F}u_{1}\subset D_{1}\cup D_{2}$.
In this case, by using Lemma 2.5, 2.1,  since $-j+\frac{3}{2}+j\epsilon\leq s\leq1-\frac{j}{2}-j\epsilon$,  we have that
\begin{eqnarray*}
&&\left\|\Lambda^{-1}\partial_{x}(1-\partial_{x}^{2})^{-1}\prod_{j=1}^{2}(\partial_{x}u_{j})\right\|_{X^{s}}\leq C\left\|\langle k\rangle ^{-\frac{s-1}{j}}\langle \sigma \rangle ^{\frac{s-1}{j}}
(\mathscr{F}u_{1}*\mathscr{F}u_{2})\right\|_{l_{k}^{2}L_{\tau}^{2}}\nonumber\\&&
\leq C\left\|(J^{s}u_{1})(J^{s}u_{2})\right\|_{L_{xt}^{2}}\nonumber\\&&\leq
C\|u_{1}\|_{X_{s,\frac{1}{2}}}\|u_{2}\|_{X_{s,\frac{1}{2(2j+1)}+\epsilon}}\nonumber\\&&
\leq C\|u_{1}\|_{X_{s,\frac{1}{2}}}\|u_{2}\|_{X_{s,\frac{1}{2j}}}\leq C\prod_{j=1}^{2}\|u_{j}\|_{Z^{s}}.
\end{eqnarray*}
If $|\sigma|>4{\rm max}\left\{|\sigma_{1}|,|\sigma_{2}|\right\}$
 and  $\supp \mathscr{F}u_{1}\subset  D_{3}$.
In this case, by using Lemma 2.1 and the Young  inequality,  since $-j+\frac{3}{2}+j\epsilon\leq s\leq1-\frac{j}{2}-j\epsilon$,  we have that
\begin{eqnarray*}
&&\left\|\Lambda^{-1}\partial_{x}(1-\partial_{x}^{2})^{-1}\prod_{j=1}^{2}(\partial_{x}u_{j})\right\|_{X^{s}}\leq C\left\|\langle k\rangle ^{-\frac{s-1}{j}}\langle \sigma \rangle ^{\frac{s-1}{j}}
(\mathscr{F}u_{1}*\mathscr{F}u_{2})\right\|_{l_{k}^{2}L_{\tau}^{2}}\nonumber\\&&
\leq \left\|\left[J^{-\frac{s-1}{j}-1}\Lambda^{\frac{s-1}{j}+1}u_{1}\right]\left[J^{-2j}u_{2}\right]\right\|_{L_{xt}^{2}} \nonumber\\&&\leq \|u_{1}\|_{X_{-\frac{s-1}{j}-1,\frac{s-1}{j}+1}}\|\langle k\rangle ^{-2j}\mathscr{F}u_{2}\|_{l_{k}^{1}L_{\tau}^{1}}
\leq C\|u_{1}\|_{X_{-\frac{s-1}{j}-1,\frac{s-1}{j}+1}}\|u_{2}\|_{Y^{s}}\nonumber\\
&&\leq C\prod_{j=1}^{2}\|u_{j}\|_{Z^{s}}.
\end{eqnarray*}
If $|\sigma|\leq 4{\rm max}\left\{|\sigma_{1}|,|\sigma_{2}|\right\},$  then we have that $|\sigma|\sim |\sigma_{1}|$ or $|\sigma|\sim |\sigma_{2}|.$

\noindent
When $|\sigma|\sim |\sigma_{1}|$. In this case, we have that $\supp\left( \mathscr{F}u_{1}*\mathscr{F}u_{2}\right) \subset D_{3}.$
Since $-j+\frac{3}{2}+j\epsilon\leq s\leq1-\frac{j}{2}-j\epsilon$,  by using Lemma 2.5 and the Young inequality,
 we have that
\begin{eqnarray*}
&&\left\|\Lambda^{-1}\partial_{x}(1-\partial_{x}^{2})^{-1}\prod_{j=1}^{2}(\partial_{x}u_{j})\right\|_{X^{s}}\leq C\left\|\langle k\rangle ^{-\frac{s-1}{j}}\langle \sigma \rangle^{\frac{s-1}{j}}\left(\mathscr{F}u_{1}*\mathscr{F}u_{2}\right) \right\|_{l_{k}^{2}L_{\tau}^{2}}\nonumber\\
&&\leq C\left\|\mathscr{F}u_{1}*\left[\langle k\rangle^{2s-2}\mathscr{F}u_{2}\right]\right\|_{l_{k}^{2}L_{\tau}^{2}}\nonumber\\
&&\leq C\left\|\left[\langle k\rangle^{-\frac{s-1}{j}-1}\langle\sigma\rangle^{\frac{s-1}{j}+1}\mathscr{F}u_{1}\right]*\left[\langle k\rangle^{-2j}\mathscr{F}u_{2}\right]\right\|_{l_{k}^{2}L_{\tau}^{2}}
\nonumber\\
&&\leq C\|u_{1}\|_{X_{-\frac{s-1}{j}-1,\frac{s-1}{j}+1}}\|\langle k\rangle^{-2j}\mathscr{F}u_{2}\|_{l_{k}^{1}L_{\tau}^{1}}\nonumber\\
&&\leq C\|u_{1}\|_{X_{-\frac{s-1}{j}-1,\frac{s-1}{j}+1}}\|\langle k\rangle^{s}\mathscr{F}u_{2}\|_{l_{k}^{2}L_{\tau}^{1}}\leq
C\prod_{j=1}^{2}\|u_{j}\|_{Z^{s}}.
\end{eqnarray*}
When  $|\sigma|\sim |\sigma_{2}|$, this case can be proved similarly to case  $|\sigma|\sim |\sigma_{1}|$.

\noindent When (b)  occurs:  if  $|\sigma_{1}|>4{\rm max}\left\{|\sigma|,|\sigma_{2}|\right\}$
which yields $\supp\mathscr{F}u_{1}\subset D_{3} $.
In this case, we consider
$\mathscr{F}u_{2}\subset D_{2}\cup D_{3}$.
When $\mathscr{F}u_{2}\subset D_{2},$
 by using Lemma 2.3, 2.5, 2.1, since $-j+\frac{3}{2}+j\epsilon\leq s\leq1-\frac{j}{2}-j\epsilon$, we have that
\begin{eqnarray}
&&\left\|\Lambda^{-1}\partial_{x}(1-\partial_{x}^{2})^{-1}\prod_{j=1}^{2}(\partial_{x}u_{j})\right\|_{X^{s}}\leq C\left\|\langle k\rangle ^{s+1}\langle \sigma \rangle ^{-\frac{1}{2j}}(\mathscr{F}u_{1}*\mathscr{F}u_{2})\right\|_{l_{k}^{2}L_{\tau}^{2}}\nonumber\\
&&\leq C\left\|(J^{-\frac{s-1}{j}-1}\Lambda ^{\frac{s-1}{j}+1}u_{1})(J^{-s-(2j-3)}u_{2})\right\|_{X_{0,-\frac{1}{2j}}}\leq C\|u_{1}\|_{X_{-\frac{s-1}{j}-1,\frac{s-1}{j}+1}}\|u_{2}\|_{X_{s,\frac{1}{2}}}\nonumber\\
&&\leq
C\prod_{j=1}^{2}\|u_{j}\|_{Z^{s}}.
\end{eqnarray}
When $\mathscr{F}u_{2}\subset D_{3},$
since $-j+\frac{3}{2}+j\epsilon\leq s\leq1-\frac{j}{2}-j\epsilon$,  by using Lemmas  2.3,  2.5 and the Young inequality, we have that
\begin{eqnarray*}
&&\left\|\Lambda^{-1}\partial_{x}(1-\partial_{x}^{2})^{-1}
\prod_{j=1}^{2}(\partial_{x}u_{j})\right\|_{X^{s}}\leq C\left\|\langle k\rangle ^{s+1}\langle \sigma\rangle ^{-\frac{1}{2j}}
\left[\mathscr{F}u_{1}*\mathscr{F}u_{2}\right]\right\|_{l_{k}^{2}L_{\tau}^{2}}\nonumber\\
&&\leq C\left\|\langle k\rangle ^{s+\frac{3}{2}+\epsilon}\langle \sigma \rangle^{-\frac{1}{2j}+\frac{1}{2}+\epsilon}
\left[\mathscr{F}u_{1}*\mathscr{F}u_{2}\right]
\right\|_{l_{k}^{\infty}L_{\tau}^{\infty}}\nonumber\\
&&\leq C\left\|\left(\langle k\rangle ^{s+j+\frac{1}{2}+(2j+1)\epsilon}\mathscr{F}u_{1}
\right)*\mathscr{F}u_{2}\right\|_{l_{k}^{\infty}l_{\tau}^{\infty}}\nonumber\\&&\leq C
\left\|\langle k\rangle ^{-3s+\frac{9}{2}-3j+(2j+1)\epsilon}\right\|_{l_{k}^{\infty}}
\prod_{j=1}^{2}\|u_{j}\|_{X_{-\frac{s-1}{j}-1,\frac{s-1}{j}+1}}\nonumber\\&&\leq C
\prod_{j=1}^{2}\|u_{j}\|_{X_{-\frac{s-1}{j}-1,\frac{s-1}{j}+1}}\leq C\prod_{j=1}^{2}\|u_{j}\|_{Z^{s}};
\end{eqnarray*}
If $|\sigma_{1}|\leq 4{\rm max}\left\{|\sigma|,|\sigma_{2}|\right\},$ we have that $|\sigma_{1}|\sim |\sigma|$ or $|\sigma_{1}|\sim |\sigma_{2}|.$

\noindent
When $|\sigma_{1}|\sim |\sigma|$, this case can be proved similarly to $|\sigma|={\rm max}\left\{|\sigma|,|\sigma_{1}|,|\sigma_{2}|\right\}.$

\noindent
When $|\sigma_{1}|\sim |\sigma_{2}|,$ we consider $\supp\left( \mathscr{F}u_{1}*\mathscr{F}u_{2}\right) \subset D_{1},$
$\supp\left( \mathscr{F}u_{1}*\mathscr{F}u_{2}\right) \subset D_{2},$ $\supp\left( \mathscr{F}u_{1}*\mathscr{F}u_{2}\right) \subset D_{3},$
respectively.

\noindent When $\supp\left( \mathscr{F}u_{1}*\mathscr{F}u_{2}\right) \subset D_{1},$
by using Lemmas 2.3, 2.5,
 since $-j+\frac{3}{2}+j\epsilon\leq s\leq1-\frac{j}{2}-j\epsilon$,  we have that
\begin{eqnarray*}
&&\left\|\Lambda^{-1}\partial_{x}(1-\partial_{x}^{2})^{-1}
\prod_{j=1}^{2}(\partial_{x}u_{j})\right\|_{X^{s}}\leq C\left\|\langle k\rangle ^{s+1}\langle \sigma\rangle ^{-\frac{1}{2j}}
\left[\mathscr{F}u_{1}*\mathscr{F}u_{2}\right]\right\|_{l_{k}^{2}L_{\tau}^{2}}\nonumber\\
&&\leq C\left\|\langle k\rangle ^{s+\frac{3}{2}+\epsilon}\langle \sigma \rangle^{-\frac{1}{2j}+\frac{1}{2}+\epsilon}
\left[\mathscr{F}u_{1}*\mathscr{F}u_{2}\right]
\right\|_{l_{k}^{\infty}L_{\tau}^{\infty}}\nonumber\\
&&\leq C\left\|\left(\langle k\rangle ^{s+j+\frac{1}{2}+(2j+1)\epsilon}\mathscr{F}u_{1}
\right)*\mathscr{F}u_{2}\right\|_{l_{k}^{\infty}l_{\tau}^{\infty}}\nonumber\\&&\leq C
\left\|\langle k\rangle ^{-3s+\frac{9}{2}-3j+(2j+1)\epsilon}\right\|_{l_{k}^{\infty}}
\prod_{j=1}^{2}\|u_{j}\|_{X_{-\frac{s-1}{j}-1,\frac{s-1}{j}+1}}\nonumber\\&&\leq C
\prod_{j=1}^{2}\|u_{j}\|_{X_{-\frac{s-1}{j}-1,\frac{s-1}{j}+1}}\leq C\prod_{j=1}^{2}\|u_{j}\|_{Z^{s}}.
\end{eqnarray*}
\noindent When $\supp\left( \mathscr{F}u_{1}*\mathscr{F}u_{2}\right) \subset D_{2},$
by using Lemmas 2.3, 2.5,
 since $-j+\frac{3}{2}+j\epsilon\leq s\leq1-\frac{j}{2}-j\epsilon$,  we have that
\begin{eqnarray*}
&&\left\|\Lambda^{-1}\partial_{x}(1-\partial_{x}^{2})^{-1}
\prod_{j=1}^{2}(\partial_{x}u_{j})\right\|_{X^{s}}\leq C\left\|\langle k\rangle ^{(1-2j)(s-1)+1}\langle \sigma\rangle ^{s-1}
\left[\mathscr{F}u_{1}*\mathscr{F}u_{2}\right]\right\|_{l_{k}^{2}L_{\tau}^{2}}\nonumber\\
&&\leq C\left\|\langle k\rangle ^{(1-2j)(s-1)+\frac{3}{2}+\epsilon}\langle \sigma \rangle^{s-\frac{1}{2}+\epsilon}
\left[\mathscr{F}u_{1}*\mathscr{F}u_{2}\right]
\right\|_{l_{k}^{\infty}L_{\tau}^{\infty}}\nonumber\\
&&\leq C\left\|\left(\langle k\rangle ^{s+j+\frac{1}{2}+(2j+1)\epsilon}\mathscr{F}u_{1}
\right)*\mathscr{F}u_{2}\right\|_{l_{k}^{\infty}l_{\tau}^{\infty}}\nonumber\\&&\leq C
\left\|\langle k\rangle ^{-3s+\frac{9}{2}-3j+(2j+1)\epsilon}\right\|_{l_{k}^{\infty}}
\prod_{j=1}^{2}\|u_{j}\|_{X_{-\frac{s-1}{j}-1,\frac{s-1}{j}+1}}\nonumber\\&&\leq C
\prod_{j=1}^{2}\|u_{j}\|_{X_{-\frac{s-1}{j}-1,\frac{s-1}{j}+1}}\leq C\prod_{j=1}^{2}\|u_{j}\|_{Z^{s}}.
\end{eqnarray*}
When $\supp\left( \mathscr{F}u_{1}*\mathscr{F}u_{2}\right) \subset D_{3},$
by using Lemmas 2.3, 2.5,
 since $-j+\frac{3}{2}+j\epsilon\leq s\leq1-\frac{j}{2}-j\epsilon$,  we have that
\begin{eqnarray*}
&&\left\|\Lambda^{-1}\partial_{x}(1-\partial_{x}^{2})^{-1}
\prod_{j=1}^{2}(\partial_{x}u_{j})\right\|_{X^{s}}\leq C\left\|\langle k\rangle ^{-\frac{s-1}{j}}\langle \sigma\rangle ^{\frac{s-1}{j}}
\left[\mathscr{F}u_{1}*\mathscr{F}u_{2}\right]\right\|_{l_{k}^{2}L_{\tau}^{2}}\nonumber\\
&&\leq C\left\|\langle k\rangle ^{-\frac{s-1}{j}+\frac{1}{2}+\epsilon}\langle \sigma \rangle^{\frac{s-1}{j}+\frac{1}{2}+\epsilon}
\left[\mathscr{F}u_{1}*\mathscr{F}u_{2}\right]
\right\|_{l_{k}^{\infty}L_{\tau}^{\infty}}\nonumber\\
&&\leq C\left\|\left(\langle k\rangle ^{2s+2j-2+(2j+2)\epsilon}\mathscr{F}u_{1}
\right)*\mathscr{F}u_{2}\right\|_{l_{k}^{\infty}l_{\tau}^{\infty}}\nonumber\\&&\leq C
\left\|\langle k\rangle ^{-2s-2j+2+(2j+2)\epsilon}\right\|_{l_{k}^{\infty}}
\prod_{j=1}^{2}\|u_{j}\|_{X_{-\frac{s-1}{j}-1,\frac{s-1}{j}+1}}\nonumber\\&&\leq C
\prod_{j=1}^{2}\|u_{j}\|_{X_{-\frac{s-1}{j}-1,\frac{s-1}{j}+1}}\leq C\prod_{j=1}^{2}\|u_{j}\|_{Z^{s}}.
\end{eqnarray*}
When (c) occurs: this case can be proved similarly to case (b).

We have completed the proof of Lemma 3.1.

\noindent {\bf Remark 4.} Regions $\Omega_{2}$ determines the indices $-j+\frac{3}{2}+j\epsilon\leq s\leq1-\frac{j}{2}-j\epsilon$.

\begin{Lemma}\label{Lemma3.2}
Let $j\geq 2$ and $-j+\frac{3}{2}+j\epsilon\leq s\leq1-\frac{j}{2}-j\epsilon$. Then, we have that
\begin{eqnarray}
      \left\|\Lambda^{-1}\partial_{x}(1-\partial_{x}^{2})^{-1}\prod_{j=1}^{2}(\partial_{x}u_{j})\right\|_{Y^{s}}
      \leq C\prod\limits_{j=1}^{2}\|u_{j}\|_{Z^{s}}.
        \label{3.01}
\end{eqnarray}
\end{Lemma}
{\bf Proof.}
 Obviously, $\left(\R\times\dot{Z}_{\lambda}\right)^{2}\subset \bigcup\limits_{j=1}^{8}\Omega_{j},$
where $\Omega_{j}(1\leq j\leq8)$ are defined as Lemma 3.1.

(1) In region $\Omega_{1}$.
By using the Lemma 2.3 and the H\"older inequality as well
as the Cauchy-Schwartz inequality, since $-j+\frac{3}{2}+j\epsilon\leq s\leq1-\frac{j}{2}-j\epsilon$, we have that
\begin{eqnarray*}
&&\left\|\Lambda^{-1}\partial_{x}(1-\partial_{x}^{2})^{-1}
\prod_{j=1}^{2}(\partial_{x}u_{j})\right\|_{Y^{s}}\leq C\left\|\Lambda^{-1}\partial_{x}(1-\partial_{x}^{2})^{-1}
\prod_{j=1}^{2}(\partial_{x}u_{j})\right\|_{X_{s,\frac{2j-1}{2j}}}\nonumber\\
&&\leq C\left\||k|\langle \sigma\rangle^{-\frac{1}{2j}}
\left(\mathscr{F}u_{1}*\mathscr{F}u_{2}\right)
\right\|_{l_{k}^{2}L_{\tau}^{2}}\leq C\|k\|_{l_{k}^{2}}\left\|\mathscr{F}u_{1}*\mathscr{F}u_{2}
\right\|_{l_{k}^{\infty}L_{\tau}^{2}}\leq C\|\mathscr{F}u_{1}\|_{l_{k}^{2}L_{\tau}^{2}}
\|\mathscr{F}u_{2}\|_{l_{k}^{2}L_{\tau}^{1}}\nonumber\\
&&\leq
 C\|u_{1}\|_{X_{s,\frac{1}{2j}}}\|u_{2}\|_{Y^{s}}
\leq C\prod_{j=1}^{2}\|u_{j}\|_{Z^{s}}.
\end{eqnarray*}
(2) In region $\Omega_{2}$.
In this case, we consider case (a)-(c) of Lemma 2.5, respectively.

\noindent When (a) is valid,
since $-j+\frac{3}{2}+j\epsilon\leq s\leq1-\frac{j}{2}-j\epsilon$, by using Lemma 2.5 and the Young inequality,
  we have that
\begin{eqnarray*}
&&\left\|\Lambda^{-1}\partial_{x}(1-\partial_{x}^{2})^{-1}\prod_{j=1}^{2}(\partial_{x}u_{j})\right\|_{Y^{s}}\leq\left\|\langle k\rangle^{s-1}\langle\sigma\rangle^{-1}\left[(|k|\mathscr{F}u_{1})*(|k|\mathscr{F}u_{2})\right]\right\|_{l_{k}^{2}L_{\tau}^{1}}\nonumber\\
&&\leq C\left\|(\langle k\rangle ^{-j+1}\mathscr{F}u_{1})*(\langle k\rangle ^{-j+1}\mathscr{F}u_{2})\right\|_{l_{k}^{\infty}L_{\tau}^{1}}\nonumber\\
&&\leq \left\|(\langle k\rangle ^{s}\mathscr{F}u_{1})*(\langle k\rangle ^{s}\mathscr{F}u_{2})\right\|_{l_{k}^{\infty}L_{\tau}^{1}}\leq C\prod_{j=1}^{2}\|u_{j}\|_{Y^{s}}\leq C\prod_{j=1}^{2}\|u_{j}\|_{Z^{s}}.
\end{eqnarray*}
When (b) is valid,
 we consider the following cases:
\begin{eqnarray*}
(i): |\sigma_{1}|>4{\rm max}\left\{|\sigma|,|\sigma_{2}|\right\},\quad
(ii):|\sigma_{1}|\leq4{\rm max}\left\{|\sigma|,|\sigma_{2}|\right\},
\end{eqnarray*}
respectively.

\noindent When (i) occurs: we consider $\supp u_{1}\subset D_{1}$, $\supp u_{1}\subset D_{2}$, $\supp u_{1}\subset D_{3}$, respectively.

\noindent When $\supp u_{1}\subset D_{1}$ which yields that $|k|\leq C$, by using Lemmas 2.3, 2.1, 2.5,   since $-j+\frac{3}{2}+j\epsilon\leq s\leq1-\frac{j}{2}-j\epsilon$,
 we have that
\begin{eqnarray*}
&&\left\|\Lambda^{-1}\partial_{x}(1-\partial_{x}^{2})^{-1}
\prod_{j=1}^{2}(\partial_{x}u_{j})\right\|_{Y^{s}}\leq C\left\|\Lambda^{-1}\partial_{x}(1-\partial_{x}^{2})^{-1}
\prod_{j=1}^{2}(\partial_{x}u_{j})\right\|_{Z^{s}}\nonumber\\&&
\leq C\left\|\partial_{x}(1-\partial_{x}^{2})^{-1}
\prod_{j=1}^{2}(\partial_{x}u_{j})\right\|_{X_{s,-\frac{1}{2j}}}
\leq C\left\|\langle k\rangle ^{s-1}\langle \sigma \rangle ^{-\frac{1}{2j}}
(|k|\mathscr{F}u_{1})*(|k|\mathscr{F}u_{2})\right\|_{l_{k}^{2}L_{\tau}^{2}}\nonumber\\&&
\leq C\left\|\langle \sigma \rangle ^{-\frac{1}{2j}}(\langle k\rangle^{s}\langle\sigma\rangle^{\frac{2j-1}{2j}}\mathscr{F}u_{1})*(\langle k\rangle ^{-s-2j+3}\mathscr{F}u_{2})\right\|_{l_{k}^{2}L_{\tau}^{2}}\nonumber\\
&&\leq C\left\|\left(J^{s}\Lambda ^{\frac{2j-1}{2j}}u_{1}\right)
\left(J^{-s-2j+3}u_{2}\right)\right\|_{X_{0,-\frac{1}{2j}}}\leq C\|u_{1}\|_{X_{s,\frac{2j-1}{2j}}}\|u_{2}\|_{X_{s,\frac{1}{2}}}\leq
C\prod_{j=1}^{2}\|u_{j}\|_{Z^{s}}.
\end{eqnarray*}
When $\supp u_{1}\subset D_{2},$   by using Lemmas 2.3, 2.1, 2.5,
since $-j+\frac{3}{2}+j\epsilon\leq s\leq1-\frac{j}{2}-j\epsilon$,
 we have that
\begin{eqnarray*}
&&\left\|\Lambda^{-1}\partial_{x}(1-\partial_{x}^{2})^{-1}
\prod_{j=1}^{2}(\partial_{x}u_{j})\right\|_{Y^{s}}\leq C\left\|\Lambda^{-1}\partial_{x}(1-\partial_{x}^{2})^{-1}
\prod_{j=1}^{2}(\partial_{x}u_{j})\right\|_{Z^{s}}\nonumber\\&&
\leq C\left\|\partial_{x}(1-\partial_{x}^{2})^{-1}
\prod_{j=1}^{2}(\partial_{x}u_{j})\right\|_{X_{s,-\frac{1}{2j}}}
\nonumber\\&&\leq \left\|\langle k\rangle ^{s-1}\langle \sigma \rangle ^{-\frac{1}{2j}}
\left[(|k|\mathscr{F}u_{1})*(|k|\mathscr{F}u_{2})\right]\right\|_{l_{k}^{2}L_{\tau}^{2}}\nonumber\\&&
\leq C\left\|\left(J^{(1-2j)(s-1)}\Lambda ^{s}u_{1}\right)\left(J^{-s-2j+3}u_{2}\right)
\right\|_{X_{0,-\frac{1}{2j}}}\nonumber\\
&&\leq C\|u_{1}\|_{X_{(1-2j)(s-1),s}}\|u_{2}\|_{X_{s,\frac{1}{2}}}\leq
C\prod_{j=1}^{2}\|u_{j}\|_{Z^{s}}.
\end{eqnarray*}
When (ii) occurs: we have that $|\sigma_{1}|\sim |\sigma|$ or $|\sigma_{1}|\sim |\sigma_{2}|$.

\noindent
When $|\sigma_{1}|\sim |\sigma|$  is valid,
this case can be proved similarly to
$|\sigma|={\rm max}\left\{|\sigma|,|\sigma_{1}|,|\sigma_{2}|\right\}.$
When $|\sigma_{1}|\sim |\sigma_{2}|$, we consider $\supp u_{1}\subset D_{1}$, $\supp u_{1}\subset D_{2}$, $\supp u_{1}\subset D_{3}$, respectively.

\noindent When  $\supp u_{1}\subset D_{1}$ which yields that $|k|\leq C$,
by using Lemmas 2.3, 2.1, 2.5,   since $-j+\frac{3}{2}+j\epsilon\leq s\leq1-\frac{j}{2}-j\epsilon$,
 we have that
\begin{eqnarray*}
&&\left\|\Lambda^{-1}\partial_{x}(1-\partial_{x}^{2})^{-1}
\prod_{j=1}^{2}(\partial_{x}u_{j})\right\|
_{Y^{s}}\leq C\left\|\Lambda^{-1}\partial_{x}(1-\partial_{x}^{2})^{-1}
\prod_{j=1}^{2}(\partial_{x}u_{j})
\right\|_{Z^{s}}\nonumber\\&&
\leq C\left\|\partial_{x}(1-\partial_{x}^{2})^{-1}
\prod_{j=1}^{2}(\partial_{x}u_{j})\right\|
_{X_{s,-\frac{1}{2j}}}\leq
\left\|\langle k\rangle ^{s-1}\langle \sigma \rangle ^{-\frac{1}{2j}}
\left[|k|\mathscr{F}u_{1}*|k|\mathscr{F}u_{2}\right]
\right\|_{l_{k}^{2}L_{\tau}^{2}}\nonumber\\&&
\leq C
\left\|\langle \sigma \rangle ^{-\frac{1}{2j}}(\langle k\rangle^{s}\langle\sigma\rangle^{\frac{2j-1}{2j}}
\mathscr{F}u_{1})*(\langle k\rangle ^{-s-2j+3}\mathscr{F}u_{2})
\right\|_{l_{k}^{2}L_{\tau}^{2}}\nonumber\\
&&\leq C\left\|
\left(J^{s}\Lambda ^{\frac{2j-1}{2j}}u_{1}\right)
\left(J^{-s-2j+3}u_{2}\right)\right\|_{X_{0,-\frac{1}{2j}}}\leq  C
\|u_{1}\|_{X_{s,\frac{2j-1}{2j}}}\|u_{2}\|_{X_{s,\frac{1}{2}}}\leq
C\prod_{j=1}^{2}\|u_{j}\|_{Z^{s}}.
\end{eqnarray*}
When $\supp u_{1}\subset D_{2},$ we  can assume that $\supp u_{2}\subset D_{2}$ and  $|\sigma|\leq C|k_{1}|^{2j+1},$
by using $X_{s,\frac{1}{2}+\epsilon}\hookrightarrow Y^{s}$ and
the H\"older  inequality  as well as the  Young inequality,
since $-j+\frac{3}{2}+j\epsilon\leq s\leq1-\frac{j}{2}-j\epsilon$,
 we  have that
\begin{eqnarray}
&&\left\|\Lambda^{-1}\partial_{x}(1-\partial_{x}^{2})^{-1}
\prod_{j=1}^{2}(\partial_{x}u_{j})\right\|_{Y^{s}}\leq C\left\|\Lambda^{-1}\partial_{x}(1-\partial_{x}^{2})^{-1}
\prod_{j=1}^{2}(\partial_{x}u_{j})\right\|_{Z^{s}}\nonumber\\&&
\leq C\left\|\partial_{x}(1-\partial_{x}^{2})^{-1}
\prod_{j=1}^{2}(\partial_{x}u_{j})\right\|_{X_{s,-\frac{1}{2}+\epsilon}}
\nonumber\\&&\leq \left\|\langle k\rangle ^{s-1}\langle
\sigma \rangle ^{-\frac{1}{2}+\epsilon}
\left[(|k|\mathscr{F}u_{1})*(|k|\mathscr{F}u_{2})\right]
\right\|_{l_{k}^{2}L_{\tau}^{2}}\nonumber\\&&
\leq C\left\|\langle k\rangle ^{s-\frac{1}{2}+\epsilon}\langle
\sigma \rangle^{2\epsilon}
\left[(|k|\mathscr{F}u_{1})*(|k|\mathscr{F}u_{2})\right]
\right\|_{l_{k}^{\infty}L_{\tau}^{\infty}}\nonumber\\
&&\leq C\left\|\left(\langle k\rangle ^{2+(4j+2)\epsilon}\mathscr{F}u_{1}
\right)*\mathscr{F}u_{2}\right\|_{l_{k}^{\infty}l_{\tau}^{\infty}}\nonumber\\&&\leq C
\left\|\langle k\rangle ^{-4s+4-4j+(4j+2)\epsilon}\right\|_{l_{k}^{\infty}}
\prod_{j=1}^{2}\|u_{j}\|_{X_{(1-2j)(s-1),s}}\nonumber\\&&\leq C
\prod_{j=1}^{2}\|u_{j}\|_{X_{(1-2j)(s-1),s}}\leq C\prod_{j=1}^{2}\|u_{j}\|_{Z^{s}}.
\end{eqnarray}
(c) Case $|\sigma_{2}|={\rm max}\left\{|\sigma|,|\sigma_{1}|,|\sigma_{2}|\right\}.$
 This case can be proved similarly to case (b).

\noindent (3) Region $\Omega_{3}$.
We  consider $|k|\leq |k_{1}|^{-2j}$  and  $|k_{1}|^{-2j}\leq|k|\leq 1,$
respectively.

\noindent
When $|k|\leq |k_{1}|^{-2j}$, by using Lemma 2.3 and the Young inequality, since $-j+\frac{3}{2}+j\epsilon\leq s\leq1-\frac{j}{2}-j\epsilon$, we have that
\begin{eqnarray*}
&&\left\|\Lambda^{-1}\partial_{x}(1-\partial_{x}^{2})^{-1}\prod_{j=1}^{2}
(\partial_{x}u_{j})\right\|_{Y^{s}}\leq C
\left\||k|\langle \sigma\rangle ^{-\frac{1}{2j}}
\left[(|k|\mathscr{F}u_{1})*(|k|\mathscr{F}u_{2})\right]
\right\|_{l_{k}^{2}L_{\tau}^{2}}\nonumber\\
&&\leq C\left\|\left[(\langle k\rangle^{-(\frac{3j}{2}-1)}
\mathscr{F}u_{1})*(\langle k\rangle^{-(\frac{3j}{2}-1)}
\mathscr{F}u_{2})\right]
\right\|_{l_{k}^{\infty}L_{\tau}^{2}}\nonumber\\
&&\leq C\|u_{1}\|_{X_{1-\frac{3j}{2},0}}\|u_{2}\|_{Y^{1-2j}}\leq C\|u_{1}\|_{X_{1-\frac{3j}{2},0}}\|u_{2}\|_{Y^{s}}\leq C
\prod_{j=1}^{2}\|u_{j}\|_{Z^{s}}.
\end{eqnarray*}
When  $|k_{1}|^{-2j}\leq|k|\leq 1,$
we consider  (a)-(c) of Lemma 2.5,  respectively.

\noindent
When (a) occurs:  by using the H\"older inequality and the Young
inequality,  since $-j+\frac{3}{2}+j\epsilon\leq s\leq1-\frac{j}{2}-j\epsilon$, we have that
\begin{eqnarray*}
&&\left\|\Lambda^{-1}\partial_{x}(1-\partial_{x}^{2})^{-1}
\prod_{j=1}^{2}(\partial_{x}u_{j})\right\|_{Y^{s}}\leq
 C\left\||k|\langle k\rangle^{s-2}\langle \sigma \rangle^{-1}\left[(|k|\mathscr{F}u_{1})*(|k|\mathscr{F}u_{2})\right]
\right\|_{l_{k}^{2}L_{\tau}^{1}}\nonumber\\
&&\leq C\left\|\left[(\langle k\rangle ^{-j+1}\mathscr{F}u_{1})*(\langle k\rangle ^{-j+1}\mathscr{F}u_{2})\right]
\right\|_{l_{k}^{\infty}L_{\tau}^{1}}\nonumber\\
&&\leq C\prod_{j=1}^{2}\|\langle k\rangle^{1-j}\mathscr{F}u_{j}\|_{l_{k}^{2}L_{\tau}^{1}}\leq C\prod_{j=1}^{2}\|u_{j}\|_{Y^{s}}\leq C\prod_{j=1}^{2}\|u_{j}\|_{Z^{s}}.
\end{eqnarray*}
When  (b) occurs:  we consider  $|\sigma_{1}|>4{\rm max}\left\{|\sigma|,|\sigma_{2}|\right\}$  and
$|\sigma_{1}|\leq4{\rm max}\left\{|\sigma|,|\sigma_{2}|\right\}$, respectively.

\noindent
When $|\sigma_{1}|>4{\rm max}\left\{|\sigma|,|\sigma_{2}|\right\}$,
 $\supp \mathscr{F}u_{1}\subset D_{1}$, by using $X_{s,\frac{1}{2}+\epsilon}\hookrightarrow Y^{s}$,
 the H\"older inequality and the Young inequality, since $-j+\frac{3}{2}+j\epsilon\leq s\leq1-\frac{j}{2}-j\epsilon$, we have that
\begin{eqnarray*}
&&\left\|\Lambda^{-1}\partial_{x}(1-\partial_{x}^{2})^{-1}
\prod_{j=1}^{2}(\partial_{x}u_{j})\right\|_{Y^{s}}\leq C\left\||k|\langle \sigma \rangle ^{-\frac{1}{2}+\epsilon}\left[(|k|\mathscr{F}u_{1})*(|k|\mathscr{F}u_{2})\right]
\right\|_{l_{k}^{2}L_{\tau}^{2}}\nonumber\\
&&\leq C\left\||k|^{\frac{1}{2j}}\left(\langle k\rangle^{s}\langle \sigma
\rangle ^{\frac{2j-1}{2j}}\mathscr{F}u_{1}\right)*\left(\langle k\rangle ^{-s-2j+3}
\mathscr{F}u_{2}\right)\right\|_{l_{k}^{2}L_{\tau}^{2}}\nonumber\\
&&\leq C\left\|\left(\langle k\rangle^{s}\langle \sigma \rangle^{\frac{2j-1}{2j}}\mathscr{F}u_{1}\right)*\left(\langle k\rangle ^{-s-2j+3}\mathscr{F}u_{2}\right)\right\|_{l_{k}^{\infty}L_{\tau}^{2}}\nonumber\\
&&\leq C\|u_{1}\|_{X_{s,\frac{2j-1}{2j}}}\|u_{2}\|_{Y^{s}}\leq C\prod_{j=1}^{2}\|u_{j}\|_{Z^{s}}..
\end{eqnarray*}
When $|\sigma_{1}|\leq4{\rm max}\left\{|\sigma|,|\sigma_{2}|\right\}$, we have
that $|\sigma_{1}|\sim |\sigma|$ or $|\sigma_{1}|\sim |\sigma_{2}|.$

\noindent
When $|\sigma_{1}|\sim |\sigma|$, this case can be proved similarly to case
$|\sigma|={\rm max}\left\{|\sigma|,|\sigma_{1}|,|\sigma_{2}|\right\}.$
When $|\sigma_{1}|\sim |\sigma_{2}|$, we consider $\supp \mathscr{F}u_{j}\subset D_{1}$, $\supp \mathscr{F}u_{j}\subset D_{2}$, $\supp \mathscr{F}u_{j}\subset D_{3}$, respectively.

\noindent When $\supp \mathscr{F}u_{j}\subset D_{1}$
with $j=1,2,$ by using $X_{s,\frac{1}{2}+\epsilon}\hookrightarrow Y^{s}$,
 the H\"older inequality and the Young inequality, since $-j+\frac{3}{2}+j\epsilon\leq s\leq1-\frac{j}{2}-j\epsilon$, we have that
\begin{eqnarray*}
&&\left\|\Lambda^{-1}\partial_{x}(1-\partial_{x}^{2})^{-1}
\prod_{j=1}^{2}(\partial_{x}u_{j})\right\|_{Y^{s}}
\leq C\left\||k|\langle \sigma \rangle ^{-\frac{1}{2}+\epsilon}\left[(|k|\mathscr{F}u_{1})*(|k|\mathscr{F}u_{2})\right]
\right\|_{l_{k}^{2}L_{\tau}^{2}}\nonumber\\
&&\leq C\left\|\left(\langle k\rangle^{s}
\langle \sigma \rangle ^{\frac{2j-1}{2j}}\mathscr{F}u_{1}\right)*\left(\langle k\rangle ^{-s-2j+3}\mathscr{F}u_{2}\right)\right\|_{l_{k}^{2}L_{\tau}^{2}}\nonumber\\
&&\leq C\left\|\left(\langle k\rangle^{s}\langle \sigma \rangle^{\frac{2j-1}{2j}}
\mathscr{F}u_{1}\right)*\left(\langle k\rangle ^{-s-2j+3}\mathscr{F}u_{2}\right)
\right\|_{l_{k}^{\infty}L_{\tau}^{2}}\nonumber\\
&&\leq C\prod_{j=1}^{2}\|u_{j}\|_{X_{s,\frac{2j-1}{2j}}}\leq C
\prod_{j=1}^{2}\|u_{j}\|_{Z^{s}}.
\end{eqnarray*}
When $\supp \mathscr{F}u_{j}\subset D_{2}$ with $j=1,2,$ by using
$X_{s,\frac{1}{2}+\epsilon}\hookrightarrow Y^{s}$,
 the H\"older inequality and the Young inequality, since $-j+\frac{3}{2}+j\epsilon\leq s\leq1-\frac{j}{2}-j\epsilon$, we have that
\begin{eqnarray*}
&&\left\|\Lambda^{-1}\partial_{x}(1-\partial_{x}^{2})^{-1}
\prod_{j=1}^{2}(\partial_{x}u_{j})\right\|_{Y^{s}}\leq
 C\left\||k|\langle \sigma \rangle ^{-\frac{1}{2}+\epsilon}\left[(|k|\mathscr{F}u_{1})*(|k|\mathscr{F}u_{2})\right]
 \right\|_{l_{k}^{2}L_{\tau}^{2}}\nonumber\\
&&\leq C\left\|\left(\langle k\rangle^{(1-2j)(s-1)}\langle \sigma \rangle ^{s}
\mathscr{F}u_{1}\right)*\left(\langle k\rangle ^{-s-2j+3}\mathscr{F}u_{2}\right)\right\|_{X^{0,-\frac{1}{2}+\epsilon}}\nonumber\\
&&\leq C\left\|u_{1}\right\|_{X^{(1-2j)(s-1),s}}\|u\|_{X_{-s-2j+3,\frac{1}{2}}}\nonumber\\
&&\leq C\|u_{1}\|_{X_{(1-2j)(s-1),s}}\|u_{2}\|_{X_{s,\frac{1}{2}}}\leq C
\prod_{j=1}^{2}\|u_{j}\|_{Z^{s}}.
\end{eqnarray*}
When $\supp \mathscr{F}u_{j}\subset D_{3}$ with $j=1,2,$ we consider
$\supp\left( \mathscr{F}u_{1}*\mathscr{F}u_{2}\right) \subset D_{1}$, $\supp\left( \mathscr{F}u_{1}*\mathscr{F}u_{2}\right) \subset D_{2}$, $\supp\left( \mathscr{F}u_{1}*\mathscr{F}u_{2}\right) \subset D_{3}$, respectively.

\noindent When $\supp\left( \mathscr{F}u_{1}*\mathscr{F}u_{2}\right) \subset D_{1}\cup D_{2}$, by using  $X_{s,\frac{1}{2}+\epsilon}\hookrightarrow Y^{s}$,
 the H\"older inequality and the Young inequality, since $-j+\frac{3}{2}+j\epsilon\leq s\leq1-\frac{j}{2}-j\epsilon$, we have that
\begin{eqnarray*}
&&\left\|\Lambda^{-1}\partial_{x}(1-\partial_{x}^{2})^{-1}
\prod_{j=1}^{2}(\partial_{x}u_{j})\right\|_{Y^{s}}\leq C\left\|\Lambda^{-1}\partial_{x}(1-\partial_{x}^{2})^{-1}
\prod_{j=1}^{2}(\partial_{x}u_{j})\right\|_{X_{s,-\frac{1}{2}+\epsilon}}
\nonumber\\&&\leq C\left\|\langle k\rangle ^{s-1}\langle \sigma \rangle ^{-\frac{1}{2}+\epsilon}\left((|k|\mathscr{F}u_{1})*(|k|\mathscr{F}u_{2})\right)
\right\|_{l_{k}^{2}L_{\tau}^{2}}\nonumber\\
&&\leq C\left\|\langle k\rangle ^{s-\frac{1}{2}+\epsilon}\langle \sigma \rangle ^{2\epsilon}\left((|k|\mathscr{F}u_{1})*(|k|\mathscr{F}u_{2})\right)
\right\|_{l_{k}^{\infty}L_{\tau}^{\infty}}\nonumber\\
&&\leq C\left\|\langle k\rangle ^{s-\frac{1}{2}+(4j+1)\epsilon}\left((|k|\mathscr{F}u_{1})*(|k|\mathscr{F}u_{2})\right)
\right\|_{l_{k}^{\infty}L_{\tau}^{\infty}}\nonumber\\
&&
\leq C\left\|\langle k\rangle^{-4s-4j+6}\right\|_{l_{k}^{\infty}}\prod_{j=1}^{2}
\|u_{j}\|_{X_{-\frac{s-1}{j}-1,\frac{s-1}{j}+1}}\leq
C\prod_{j=1}^{2}\|u_{j}\|_{Z^{s}}.
\end{eqnarray*}
\noindent When $\supp\left( \mathscr{F}u_{1}*\mathscr{F}u_{2}\right) \subset D_{3}$, we consider $|\sigma|\leq C|k_{1}|^{2j+1}$ and
$|\sigma|> C|k_{1}|^{2j+1}$, respectively.

\noindent When $|\sigma|\leq C|k_{1}|^{2j+1}$,
by using  $X_{s,\frac{1}{2}+\epsilon}\hookrightarrow Y^{s}$,
 the H\"older inequality and the Young inequality, since $-j+\frac{3}{2}+j\epsilon\leq s\leq1-\frac{j}{2}-j\epsilon$, we have that
\begin{eqnarray*}
&&\left\|\Lambda^{-1}\partial_{x}(1-\partial_{x}^{2})^{-1}
\prod_{j=1}^{2}(\partial_{x}u_{j})\right\|_{Y^{s}}\leq C\left\|\Lambda^{-1}\partial_{x}(1-\partial_{x}^{2})^{-1}
\prod_{j=1}^{2}(\partial_{x}u_{j})\right\|_{X_{s,-\frac{1}{2}+\epsilon}}
\nonumber\\&&\leq C\left\|\langle k\rangle ^{s-1}\langle \sigma \rangle ^{-\frac{1}{2}+\epsilon}\left((|k|\mathscr{F}u_{1})*(|k|\mathscr{F}u_{2})\right)
\right\|_{l_{k}^{2}L_{\tau}^{2}}\nonumber\\
&&\leq C\left\|\langle k\rangle ^{s-\frac{1}{2}+\epsilon}\langle \sigma \rangle ^{2\epsilon}\left((|k|\mathscr{F}u_{1})*(|k|\mathscr{F}u_{2})\right)
\right\|_{l_{k}^{\infty}L_{\tau}^{\infty}}\nonumber\\
&&\leq C\left\|\langle k\rangle ^{s-\frac{1}{2}+\epsilon}\left((|k|^{1+4j\epsilon}\mathscr{F}u_{1})*(|k|\mathscr{F}u_{2})\right)
\right\|_{l_{k}^{\infty}L_{\tau}^{\infty}}\nonumber\\
&&
\leq C\left\|\langle k\rangle^{-4s-4j+4+4j\epsilon}\right\|_{l_{k}^{\infty}}\prod_{j=1}^{2}
\|u_{j}\|_{X_{-\frac{s-1}{j}-1,\frac{s-1}{j}+1}}\leq
C\prod_{j=1}^{2}\|u_{j}\|_{Z^{s}}.
\end{eqnarray*}
 When $|\sigma|> C|k_{1}|^{2j+1}$,
by using  $X_{s,\frac{1}{2}+\epsilon}\hookrightarrow Y^{s}$,
 the H\"older inequality and the Young inequality, since $-j+\frac{3}{2}+j\epsilon\leq s\leq1-\frac{j}{2}-j\epsilon$, we have that
\begin{eqnarray*}
&&\left\|\Lambda^{-1}\partial_{x}(1-\partial_{x}^{2})^{-1}
\prod_{j=1}^{2}(\partial_{x}u_{j})\right\|_{Y^{s}}\leq C\left\|\langle k\rangle ^{s-1}\langle\sigma\rangle ^{-1}\left((|k|\mathscr{F}u_{1})*(|k|\mathscr{F}u_{2})\right)\right\|_{l_{k}^{2}L_{\tau}^{1}}
\nonumber\\&&\leq C\left\|\left((|k|^{-j+\frac{1}{2}}\mathscr{F}u_{1})*(|k|^{-j+\frac{1}{2}}\mathscr{F}u_{2})\right)\right\|_{l_{k}^{2}L_{\tau}^{1}}\nonumber\\
&&\leq C\left\|\left((|k|^{-j+\frac{1}{2}}\mathscr{F}u_{1}\right)\right\|_{l_{k}^{2}L_{\tau}^{1}}\left\|\left(|k|^{-j+\frac{1}{2}}\mathscr{F}u_{2}\right)\right\|_{l_{k}^{1}L_{\tau}^{1}}\nonumber\\
&&\leq C\left\|\left((|k|^{-j+\frac{1}{2}}\mathscr{F}u_{1}\right)\right\|_{l_{k}^{2}L_{\tau}^{1}}\left\|\langle k\rangle^{s}\left(|k|^{-j+\frac{1}{2}-s}\mathscr{F}u_{2}\right)\right\|_{l_{k}^{1}L_{\tau}^{1}}\nonumber\\
&&\leq C\prod_{j=1}^{2}\left\|\left((|k|^{-j+\frac{1}{2}}\mathscr{F}u_{j}\right)\right\|_{l_{k}^{2}L_{\tau}^{1}}\leq
C\prod_{j=1}^{2}\|u_{j}\|_{Z^{s}}.
\end{eqnarray*}

\noindent (c) Case $|\sigma_{2}|={\rm max}\left\{|\sigma|,|\sigma_{1}|,|\sigma_{2}|\right\}.$
 This case can be proved similarly to case (b).

 \noindent
(4) Region $\Omega_{4}$.
We consider  (a)-(c) of Lemma 2.5,  respectively.

\noindent
When (a) occurs:  by using Lemma 2.5, since $-j+\frac{3}{2}+j\epsilon\leq s\leq1-\frac{j}{2}-j\epsilon$,
we have that
\begin{eqnarray*}
&&\left\|\Lambda^{-1}\partial_{x}(1-\partial_{x}^{2})^{-1}
\prod_{j=1}^{2}(\partial_{x}u_{j})\right\|_{Y^{s}}\leq C\left\|\langle k\rangle^{s}\langle \sigma \rangle^{-1}\left[(|k|\mathscr{F}u_{1})* \mathscr{F}u_{2}\right]\right\|_{l_{k}^{2}L_{\tau}^{1}}\nonumber\\
&&\leq C\left\|\langle k\rangle ^{-2j}\left[\mathscr{F}u_{1}*(\langle k\rangle ^{s}\mathscr{F}u_{2})\right]\right\|_{l_{k}^{2}L_{\tau}^{1}}\nonumber\\
&&\leq C\left\|\left[(\langle k\rangle ^{-2j}\mathscr{F}u_{1})*(\langle k\rangle ^{s}\mathscr{F}u_{2})\right]\right\|_{l_{k}^{2}L_{\tau}^{1}}\nonumber\\
&&\leq C\|\langle k\rangle ^{-2j}u_{1}\|_{l_{k}^{1}L_{\tau}^{1}}\|u_{2}\|_{Y^{s}}\leq C\prod_{j=1}^{2}\|u_{j}\|_{Y^{s}}\leq C\prod_{j=1}^{2}\|u_{j}\|_{Z^{s}}.
\end{eqnarray*}
(b): $|\sigma_{1}|={\rm max}\left\{|\sigma|,|\sigma_{1}|,|\sigma_{2}|\right\}.$
In this case, we consider $|\sigma_{1}|>4{\rm max}\left\{|\sigma|,|\sigma_{2}|\right\}$
and $|\sigma_{1}|\leq4{\rm max}\left\{|\sigma|,|\sigma_{2}|\right\}$, respectively.

\noindent
If $|\sigma_{1}|>4{\rm max}\left\{|\sigma|,|\sigma_{2}|\right\}$,
 then $\supp \mathscr{F}u_{1}\subset D_{3}$ and  $\supp \mathscr{F}u_{2}\subset D_{1}\cup D_{2}$,
 by using Lemma 2.3, 2.5, 2.1,   since $-j+\frac{3}{2}+j\epsilon\leq s\leq1-\frac{j}{2}-j\epsilon$,  we have that
\begin{eqnarray*}
&&\left\|\Lambda^{-1}\partial_{x}(1-\partial_{x}^{2})^{-1}\prod_{j=1}^{2}(\partial_{x}u_{j})\right\|_{Y^{s}}\nonumber\\&&\leq C\left\|\Lambda^{-1}\partial_{x}(1-\partial_{x}^{2})^{-1}\prod_{j=1}^{2}(\partial_{x}u_{j})\right\|_{Z^{s}}
\leq C\left\|\partial_{x}(1-\partial_{x}^{2})^{-1}\prod_{j=1}^{2}(\partial_{x}u_{j})\right\|_{X_{s,-\frac{1}{2j}}}\nonumber\\&&\leq C\left\|\langle k\rangle^{s+1}\langle \sigma \rangle^{-\frac{1}{2j}}
\left[\mathscr{F}u_{1}*\mathscr{F}u_{2}\right]\right\|_{l_{k}^{2}L_{\tau}^{2}}
\nonumber\\
&&\leq C\left\|(J^{-\frac{s-1}{j}-1}\langle \sigma \rangle ^{\frac{s-1}{j}+1}u_{1})(J^{-s-(2j-3)}u_{2})\right\|_{X_{0,-\frac{1}{2j}}}
\nonumber\\&&\leq C\|u_{1}\|_{X_{-\frac{s-1}{j}-1,\frac{s-1}{j}+1}}\|u_{2}\|_{X_{s,\frac{1}{2}}}
\leq C\prod_{j=1}^{2}\|u_{j}\|_{Z^{s}}.
\end{eqnarray*}
When $|\sigma_{1}|\leq4{\rm max}\left\{|\sigma|,|\sigma_{2}|\right\}$,
we have that $|\sigma_{1}|\sim |\sigma|$ or $|\sigma_{1}|\sim |\sigma_{2}|$.

\noindent
 When  $|\sigma_{1}|\sim |\sigma|$, this case can be proved similarly
 to case $|\sigma|={\rm max}\left\{|\sigma|,|\sigma_{1}|,|\sigma_{2}|\right\}.$

\noindent When $|\sigma_{1}|\sim |\sigma_{2}|$, we have that $\supp \mathscr{F}u_{1}\subset D_{3}$.
In this case, we consider
$
\supp \mathscr{F}u_{2}\subset D_{2},\supp \mathscr{F}u_{2}\subset D_{3},
$
respectively.

\noindent When  $\supp \mathscr{F}u_{2}\subset D_{2}$, by using Lemma 2.3 and $X_{s,\frac{1}{2}+\epsilon}\hookrightarrow Y^{s}$, since $-j+\frac{3}{2}+j\epsilon\leq s\leq1-\frac{j}{2}-j\epsilon$,
 we have that
\begin{eqnarray*}
&&\left\|\Lambda^{-1}\partial_{x}(1-\partial_{x}^{2})^{-1}
\prod_{j=1}^{2}(\partial_{x}u_{j})\right\|_{Y^{s}}\nonumber\\&&\leq C
\left\|\Lambda^{-1}\partial_{x}(1-\partial_{x}^{2})^{-1}
\prod_{j=1}^{2}(\partial_{x}u_{j})\right\|_{Z^{s}}
\leq C\left\|\partial_{x}(1-\partial_{x}^{2})^{-1}
\prod_{j=1}^{2}(\partial_{x}u_{j})\right\|_{X_{s,-\frac{1}{2}+\epsilon}}\nonumber\\&&\leq C\left\|\langle \sigma \rangle^{-\frac{1}{2}+\epsilon}
\left[(|k|^{-2j+2}\mathscr{F}u_{1})*(\langle k\rangle ^{(1-2j)(s-1)}\langle\sigma\rangle^{s}\mathscr{F}u_{2})\right]
\right\|_{l_{k}^{2}L_{\tau}^{2}}\nonumber\\
&&\leq C\left\|\left(J^{-2j+2}u_{1}
\right)*\left(J^{(1-2j)(s-1)}\Lambda ^{s}u_{2}\right)\right\|_{X_{0,-\frac{1}{2}+\epsilon}}\nonumber\\&&\leq C
\|u_{1}\|_{X_{-2j+2,\frac{1}{2j}}}\|u_{2}\|_{X_{(1-2j)(s-1),s}}\nonumber\\&&\leq C
\|u_{1}\|_{X_{-\frac{s-1}{j}-1,\frac{s-1}{j}+1}}\|u_{2}\|_{X_{(1-2j)(s-1),s}}\leq C\prod_{j=1}^{2}\|u_{j}\|_{Z^{s}}.
\end{eqnarray*}
When    $\supp \mathscr{F}u_{j}\subset D_{3}$ with $j=1,2$, without loss of generality, we can assume that $|\sigma|\leq C|k|^{2j+1}$ since $|\sigma|> C|k|^{2j+1}$ can be easily proved.

 \noindent By using the Young inequality, by using Lemma 2.3 and  $X_{s,\frac{1}{2}+\epsilon}\hookrightarrow Y^{s}$, since $-j+\frac{3}{2}+j\epsilon\leq s\leq1-\frac{j}{2}-j\epsilon$, we have that
\begin{eqnarray}
&&\left\|\Lambda^{-1}(1-\partial_{x}^{2})^{-1}\prod_{j=1}^{2}(\partial_{x}u_{j})\right\|_{Y^{s}}\leq C
\left\|\Lambda^{-1}(1-\partial_{x}^{2})^{-1}\prod_{j=1}^{2}(\partial_{x}u_{j})\right\|_{Z^{s}}
\nonumber\\&&\leq C\left\|\langle k\rangle ^{s+1}\langle \sigma\rangle ^{-\frac{1}{2}+\epsilon}
\left[\mathscr{F}u_{1}*\mathscr{F}u_{2}\right]\right\|_{l_{k}^{2}L_{\tau}^{2}}\nonumber\\
&&\leq C\left\|\langle k\rangle ^{s+\frac{3}{2}+\epsilon}\langle \sigma \rangle^{2\epsilon}
\left[\mathscr{F}u_{1}*\mathscr{F}u_{2}\right]
\right\|_{l_{k}^{\infty}L_{\tau}^{\infty}}\nonumber\\
&&\leq C\left\|\left(\mathscr{F}u_{1}
\right)*(\langle k\rangle ^{s+\frac{3}{2}+(4j+3)\epsilon}\mathscr{F}u_{2})\right\|_{l_{k}^{\infty}l_{\tau}^{\infty}}\nonumber\\&&\leq C
\left\|\langle k\rangle ^{-3s+\frac{11}{2}-4j+(4j+3)\epsilon}\right\|_{l_{k}^{\infty}}
\prod_{j=1}^{2}\|u_{j}\|_{X_{-\frac{s-1}{j}-1,\frac{s-1}{j}+1}}\nonumber\\&&\leq C
\prod_{j=1}^{2}\|u_{j}\|_{X_{-\frac{s-1}{j}-1,\frac{s-1}{j}+1}}\leq C\prod_{j=1}^{2}\|u_{j}\|_{Z^{s}}.
\end{eqnarray}
When case (c) occurs:  by using Lemma 2.3, we have that
\begin{eqnarray*}
&&\left\|\Lambda^{-1}(1-\partial_{x}^{2})^{-1}\partial_{x}
\prod_{j=1}^{2}(\partial_{x}u_{j})
\right\|_{Y^{s}}\leq C\left\|\langle k\rangle ^{s}\langle
\sigma\rangle ^{-\frac{1}{2j}}
\left[|k|\mathscr{F}u_{1}*\mathscr{F}u_{2}\right]\right
\|_{l_{k}^{2}L_{\tau}^{2}}.
\end{eqnarray*}
By using a proof similar to case (c)  of region $\Omega_{4}$ of   Lemma 3.1, we can obtain that
\begin{eqnarray*}
&&\left\|(1-\partial_{x}^{2})^{-1}\partial_{x}
\prod_{j=1}^{2}(\partial_{x}u_{j})
\right\|_{Y^{s}}\leq C\prod_{j=1}^{2}\|u_{j}\|_{Z^{s}}.
\end{eqnarray*}
(5) In region $\Omega_{5}$.
 In this region, we consider  $|k_{1}|\leq |k|^{-2j}$ and $|k|^{-2j}< |k_{1}|\leq 1,$
 respectively.

 \noindent
When $|k_{1}|\leq |k|^{-2j}$, by using Lemma 2.3 and the Young inequality
as well as Cauchy-Schwartz inequality,
since $-j+\frac{3}{2}+j\epsilon\leq s\leq1-\frac{j}{2}-j\epsilon$, we have that
\begin{eqnarray*}
&&\left\|\Lambda^{-1}\partial_{x}(1-\partial_{x}^{2})^{-1}
\prod_{j=1}^{2}(\partial_{x}u_{j})\right\|_{Y^{s}}\leq
C\left\|\langle k\rangle^{s}\langle \sigma\rangle ^{-\frac{1}{2j}}
\left[(|k_{1}|\mathscr{F}u_{1})*\mathscr{F}u_{2}\right]\right\|_{l_{k}^{2}L_{\tau}^{2}}\nonumber\\
&&\leq C\left\|\langle k\rangle ^{-2j}\left[\mathscr{F}u_{1}*(\langle k\rangle^{s}\mathscr{F}u_{2})\right]\right\|_{l_{k}^{2}L_{\tau}^{2}}\nonumber\\
&&\leq C\left\|\left[\mathscr{F}u_{1}*(\langle k\rangle^{s}\mathscr{F}u_{2})\right]\right\|_{l_{k}^{\infty}L_{\tau}^{2}}
\leq C\|\mathscr{F}u_{1}\|_{l_{k}^{2}L_{\tau}^{2}}\|u_{2}\|_{Y^{s}}
\leq C\prod_{j=1}^{2}\|u_{j}\|_{Z^{s}}.
\end{eqnarray*}
When  $|k|^{-2j}\leq |k_{1}|\leq 1,$
we consider  (a)-(c) of Lemma 2.5,  respectively.

\noindent
When (a) occurs:
 by using Lemma 2.3
 and the Young inequality as well as Cauchy-Schwartz inequality,
 since $-j+\frac{3}{2}+j\epsilon\leq s\leq1-\frac{j}{2}-j\epsilon$,
 we have that
\begin{eqnarray*}
&&\left\|\Lambda^{-1}\partial_{x}(1-\partial_{x}^{2})^{-1}
\prod_{j=1}^{2}(\partial_{x}u_{j})\right\|_{Y^{s}}\leq
C\left\|\langle k\rangle^{s}\langle \sigma\rangle ^{-\frac{1}{2j}}
\left[(|k|\mathscr{F}u_{1})*\mathscr{F}u_{2}\right]
\right\|_{l_{k}^{2}L_{\tau}^{2}}\nonumber\\
&&\leq C\left\|\left[(|k|^{1-\frac{1}{2j}}\mathscr{F}u_{1})
*(\langle k\rangle ^{s-1}\mathscr{F}u_{2})\right]
\right\|_{l_{k}^{2}L_{\tau}^{2}}\nonumber\\
&&\leq C\left\||k|^{1-\frac{1}{2j}}\mathscr{F}u_{1}
\right\|_{l_{k}^{1}l_{\tau}^{2}}\|u_{2}\|_{Y^{s-1}}
\leq C\|\mathscr{F}u_{1}\|_{l_{k}^{2}L_{\tau}^{2}}\|u_{2}\|_{Y^{s}}
\leq C\prod_{j=1}^{2}\|u_{j}\|_{Z^{s}}.
\end{eqnarray*}
When (b) occurs:
 by using Lemma 2.3 and the Young
 inequality as well as  the Cauchy-Schwartz inequality,
 since $-j+\frac{3}{2}+j\epsilon\leq s\leq1-\frac{j}{2}-j\epsilon$,
 we have that
\begin{eqnarray}
&&\left\|\Lambda^{-1}\partial_{x}(1-\partial_{x}^{2})^{-1}
\prod_{j=1}^{2}(\partial_{x}u_{j})\right\|_{Y^{s}}\leq C
\left\|\langle k\rangle ^{s+1}\langle \sigma\rangle ^{-\frac{1}{2j}}
\left[\mathscr{F}u_{1}*\mathscr{F}u_{2}\right]\right\|_{l_{k}^{2}L_{\tau}^{2}}\nonumber\\
&&\leq C\left\|\left[(|k|^{-\frac{1}{2j}}\langle \sigma
\rangle ^{\frac{1}{2j}}\mathscr{F}u_{1})*(\langle k\rangle ^{s}\mathscr{F}u_{2})\right]
\right\|_{l_{k}^{2}L_{\tau}^{2}}\nonumber\\
&&\leq C\left\||k|^{-\frac{1}{2j}}\langle \sigma \rangle ^{\frac{1}{2j}}\mathscr{F}u_{1}\right\|_{l_{k}^{1}l_{\tau}^{2}}\|u_{2}\|_{Y^{s}}\nonumber\\&&
\leq C\|\langle \sigma \rangle ^{\frac{1}{2j}}\mathscr{F}u_{1}\|_{l_{k}^{2}L_{\tau}^{2}}\|u_{2}\|_{Y^{s}}
\leq C\|u_{1}\|_{X_{0,\frac{1}{2j}}}\|u_{2}\|_{Y^{s}}
\leq C\prod_{j=1}^{2}\|u_{j}\|_{Z^{s}}.
\end{eqnarray}
When (c)  occurs:
by using Lemma 2.3
and the Young inequality as well as the Cauchy-Schwartz inequality,
since $-j+\frac{3}{2}+j\epsilon\leq s\leq1-\frac{j}{2}-j\epsilon$,
we have that
\begin{eqnarray*}
&&\left\|\Lambda^{-1}\partial_{x}(1-\partial_{x}^{2})^{-1}
\prod_{j=1}^{2}(\partial_{x}u_{j})\right\|_{Y^{s}}\leq C\left\|\Lambda^{-1}\partial_{x}(\prod_{j=1}^{2}u_{j})
\right\|_{X^{s}}\nonumber\\
&&\leq C\left\|\left[(|k|^{-\frac{1}{2j}}\mathscr{F}u_{1})
*(\langle k\rangle ^{s}\langle \sigma \rangle ^{\frac{1}{2j}}\mathscr{F}u_{2})\right]
\right\|_{l_{k}^{2}L_{\tau}^{2}}\nonumber\\
&&\leq C\left\||k|^{-\frac{1}{2j}}\mathscr{F}u_{1}
\right\|_{l_{k}^{1}l_{\tau}^{1}}\|u_{2}\|_{X_{s,\frac{1}{2j}}}\nonumber\\&&
\leq C\|\mathscr{F}u_{1}\|_{l_{k}^{2}L_{\tau}^{1}}\|u_{2}\|_{Y^{s}}
\leq C\|u_{1}\|_{Y^{s}}\|u_{2}\|_{X_{s,\frac{1}{2j}}}
\leq C\prod_{j=1}^{2}\|u_{j}\|_{Z^{s}}.
\end{eqnarray*}
(6)In region $\Omega_{6}$.
This case can be proved similarly to $\Omega_{4}$.

\noindent(7)In region $\Omega_{7}$.
This case can be proved similarly to $\Omega_{7}$.

\noindent (8)In region $\Omega_{8}$.
We consider  (a)-(c) of Lemma 2.5,  respectively.

\noindent
When (a) occurs:
  by using Lemma 2.5 and the Young inequality,
 since $-j+\frac{3}{2}+j\epsilon\leq s\leq1-\frac{j}{2}-j\epsilon$
 as well as the Cauchy-Schwartz inequality,
  we have that
\begin{eqnarray*}
&&\left\|\Lambda^{-1}\partial_{x}(1-\partial_{x}^{2})^{-1}
\prod_{j=1}^{2}(\partial_{x}u_{j})\right\|_{Y^{s}}\leq C
\left\|\langle k\rangle ^{s+1}\langle \sigma \rangle ^{-1}
\left(\mathscr{F}u_{1}*\mathscr{F}u_{2}\right) \right\|_{l_{k}^{2}L_{\tau}^{1}}
\nonumber\\&&\leq C\left\|\langle k\rangle ^{s-2j}
(\mathscr{F}u_{1}*\mathscr{F}u_{2})\right\|_{l_{k}^{2}L_{\tau}^{1}}\nonumber\\&&
\leq C\left\|\langle k\rangle ^{s}\mathscr{F}u_{1}\right\|_{l_{k}^{2}L_{\tau}^{1}}\|\langle k\rangle ^{-2j}\mathscr{F}u_{2}\|_{l_{k}^{1}L_{\tau}^{1}}\leq
C\|u_{1}\|_{Y^{s}}\|u_{2}\|_{Y^{s}}\leq C\prod_{j=1}^{2}\|u_{j}\|_{Z^{s}}.
\end{eqnarray*}
When (b)  occurs: we have that
 $\supp\mathscr{F}u_{1}\subset D_{3}$.

\noindent When  $\supp\mathscr{F}u_{2}
 \subset D_{1}\cup D_{2}$,     by using $X_{s,\frac{1}{2}+\epsilon}\hookrightarrow Y^{s}$,  since $-j+\frac{3}{2}+j\epsilon\leq s\leq1-\frac{j}{2}-j\epsilon$,
 we have that
\begin{eqnarray*}
&&\left\|\Lambda^{-1}\partial_{x}(1-\partial_{x}^{2})^{-1}
\prod_{j=1}^{2}(\partial_{x}u_{j})\right\|_{Y^{s}}\leq C
\left\|\Lambda^{-1}\partial_{x}(1-\partial_{x}^{2})^{-1}
\prod_{j=1}^{2}(\partial_{x}u_{j})\right\|_{Z^{s}}
\nonumber\\&&\leq C\left\|\Lambda^{-1}\partial_{x}(1-\partial_{x}^{2})^{-1}
\prod_{j=1}^{2}(\partial_{x}u_{j})\right\|_{X_{s,-\frac{1}{2}+\epsilon}}
\leq C\left\|\langle k\rangle ^{s+1}\langle \sigma \rangle ^{-\frac{1}{2}+\epsilon}\left(\mathscr{F}u_{1}*\mathscr{F}u_{2}\right)
\right\|_{l_{k}^{2}L_{\tau}^{2}}\nonumber\\
&&\leq C\left\|(J^{-\frac{s-1}{j}-1}\Lambda
 ^{\frac{s-1}{j}+1}u_{1})(J^{-s-(2j-3)}u_{2})
\right\|_{X_{0,-\frac{1}{2}+\epsilon}}\leq C
\|u_{1}\|_{X_{-\frac{s-1}{j}-1,\frac{s-1}{j}+1}}
\|u_{2}\|_{X_{s,\frac{1}{2j}}}\nonumber\\
&&\leq
C\prod_{j=1}^{2}\|u_{j}\|_{Z^{s}}.
\end{eqnarray*}
 When
 $\supp\mathscr{F}u_{j}\subset D_{3}$ with $j=1,2,$ without loss of generality, we can assume that $|\sigma|\leq C|k_{1}|^{2j+1}$ since case
 $|\sigma|>C|k_{1}|^{2j+1}$ can be easily proved.

\noindent By using $X_{s,\frac{1}{2}+\epsilon}\hookrightarrow Y^{s}$,
  since $-j+\frac{3}{2}+j\epsilon\leq s\leq1-\frac{j}{2}-j\epsilon$,
  we have that
\begin{eqnarray*}
&&\left\|\Lambda^{-1}\partial_{x}(1-\partial_{x}^{2})^{-1}
\prod_{j=1}^{2}(\partial_{x}u_{j})\right\|_{Y^{s}}
\leq C\left\|\partial_{x}(1-\partial_{x}^{2})^{-1}
\prod_{j=1}^{2}(\partial_{x}u_{j})\right\|_{X_{s,-\frac{1}{2}+\epsilon}}
\nonumber\\&&\leq C\left\|\langle k\rangle ^{s+1}\langle \sigma \rangle ^{-\frac{1}{2}+\epsilon}\left(\mathscr{F}u_{1}*\mathscr{F}u_{2}\right)
\right\|_{l_{k}^{2}L_{\tau}^{2}}\leq C\left\|\langle k\rangle ^{s+\frac{3}{2}+\epsilon}\langle \sigma \rangle ^{2\epsilon}\left(\mathscr{F}u_{1}*\mathscr{F}u_{2}\right)
\right\|_{l_{k}^{\infty}L_{\tau}^{\infty}}\nonumber\\&&\leq C\left\|\langle k\rangle^{-3s-4j+\frac{11}{2}+(4j+2)\epsilon}\right\|\prod_{j=1}^{2}
\|u_{j}\|_{X_{-\frac{s-1}{j}-1,\frac{s-1}{j}+1}}\leq
C\prod_{j=1}^{2}\|u_{j}\|_{Z^{s}}.
\end{eqnarray*}
Case (c)  can be proved similarly to Case (b).

The proof of Lemma 3.2 is completed.

\noindent {\bf Remark 5.} Regions $\Omega_{3}, \Omega_{4}$ are the most difficult to handle. Moreover, regions $\Omega_{3}, \Omega_{4}$ determine the indices $-j+\frac{3}{2}+j\epsilon\leq s\leq1-\frac{j}{2}-j\epsilon$.

\begin{Lemma}\label{Lemma3.3}
Let $j\geq 2$ and $-j+\frac{3}{2}+j\epsilon\leq s\leq1-\frac{j}{2}-j\epsilon$. Then, we have that
\begin{eqnarray}
      \left\|\Lambda^{-1}\partial_{x}(1-\partial_{x}^{2})^{-1}
      \prod_{j=1}^{2}(\partial_{x}u_{j})
      \right\|_{Z^{s}}
      \leq C\prod\limits_{j=1}^{2}\|u_{j}\|_{Z^{s}}.
        \label{3.03}
\end{eqnarray}
\end{Lemma}
{\bf Proof.} Combining the definition of $Z^{s}$ with Lemmas 3.1, 3.2, we have Lemma 3.3.

We have completed the proof of Lemma 3.3.

By using  a proof similar to Lemma 3.3,  we have Lemmas 3.3, 3.4.

\begin{Lemma}\label{Lemma3.4}
Let $j\geq 2$ and $-j+\frac{3}{2}+j\epsilon\leq s\leq1-\frac{j}{2}-j\epsilon$. Then, we have that
\begin{eqnarray}
      \left\|\Lambda^{-1}\partial_{x}\prod_{j=1}^{2}(u_{j})
      \right\|_{Z^{s}}
      \leq C\prod\limits_{j=1}^{2}\|u_{j}\|_{Z^{s}}.
        \label{3.04}
\end{eqnarray}
\end{Lemma}

\begin{Lemma}\label{Lemma3.5}
Let $j\geq 2$ and $-j+\frac{3}{2}+j\epsilon\leq s\leq1-\frac{j}{2}-j\epsilon$.  Then, we have that
\begin{eqnarray}
      \left\|\Lambda^{-1}\partial_{x}(1-\partial_{x}^{2})^{-1}\prod_{j=1}^{2}(u_{j})
      \right\|_{Z^{s}}
      \leq C\prod\limits_{j=1}^{2}\|u_{j}\|_{Z^{s}}.
        \label{3.05}
\end{eqnarray}
\end{Lemma}

\bigskip
\bigskip

\noindent {\large\bf 4. Proof of Theorem  1.1}

\setcounter{equation}{0}

 \setcounter{Theorem}{0}

\setcounter{Lemma}{0}

\setcounter{section}{4}

In this section, we present the proof of Theorem 1.1.

\begin{proof} We assume that $N\gg1$, $a\in \dot{Z}$ and
\begin{eqnarray*}
&&\mathscr{F}u_{1}(k,\tau)=\left(\chi_{(N)}(k)+\chi_{(N)}(-k)\right)\chi_{[-1,1]}(\tau+(-1)^{k}k^{2k+1}),
\nonumber\\&& \mathscr{F}u_{2}(k,\tau)=\left(\chi_{(1-N)}(k)+\chi_{(1-N)}(-k)\right)\chi_{[-1,1]}(\tau+(-1)^{k}k^{2k+1}),
\end{eqnarray*}
Where
\begin{eqnarray*}
\chi_{a}(k)=1 \quad {\rm if}\quad k=a, \chi_{a}(k)=0\quad {\rm if}\quad k\neq a,
\end{eqnarray*}
and
\begin{eqnarray*}
\chi_{[-1,1]}(\sigma)=1\quad {\rm if} \quad |\sigma|\leq 1, \chi_{[-1,1]}=0,\quad {\rm if}\quad |\sigma|>1.
\end{eqnarray*}
Obviously, by a direct computation, we have that
\begin{eqnarray*}
\|u_{j}\|_{W^{s}}\sim N^{s},j=1,2.
\end{eqnarray*}
Let
\begin{eqnarray*}
&&R_{1}(k_{1},k_{2})=\chi_{N}(k_{1})\chi_{(1-N)}(k_{2}), R_{2}(k_{1},k_{2})=\chi_{N}(k_{1})\chi_{(1-N)}(-k_{2}), \nonumber\\&&
R_{3}(k_{1},k_{2})=\chi_{N}(-k_{1})\chi_{(1-N)}(k_{2}), R_{4}(k_{1},k_{2})=\chi_{N}(-k_{1})\chi_{(1-N)}(-k_{2}).
\end{eqnarray*}
Then, we derive that
\begin{eqnarray*}
&&\left\|\mathscr{F}^{-1}\left[\langle \tau+(-1)^{j}k^{2j+1}\rangle^{-1}\mathscr{F}F(u_{1},u_{2})\right]\right\|_{W^{s}}\nonumber\\&&\hspace{-1cm}
=\left\|\sum_{j=1}^{4}\int_{\dot{Z}}\frac{|k|^{s+1}}{1+k^{2}}
\left[k^{2}+3+k_{1}k_{2}\right]R_{j}(k_{1},k_{2})\left(\int_{\SR}\langle\sigma\rangle^{-1/2}
\chi_{[-1,1]}(\sigma_{1})\chi_{[-1,1]}(\sigma_{2})d\sigma_{1}\right)dk_{1}\right\|_{l_{k}^{2}L_{\sigma}^{2}}
\end{eqnarray*}
By using  Lemma 2.7, we obtain that
\begin{eqnarray*}
\langle\sigma\rangle \sim |k_{\min}||k_{\max}|^{2j}
\end{eqnarray*}
since $|\sigma_{j}|\leq 1$ with $j=1,2.$ Thus, we have that
\begin{eqnarray*}
\int_{\SR}\langle \sigma\rangle^{-1/2}\chi_{[-1,1]}(\sigma_{1})\chi_{[-1,1]}(\sigma_{2})d\sigma_{1}\geq C|k_{\min}|^{-1/2}|k_{\max}|^{-j}.
\end{eqnarray*}
By using a direct computation, we obtain that
\begin{eqnarray*}
&&\left\|\mathscr{F}^{-1}\left[\langle \tau-k^{3}\rangle^{-1}\mathscr{F}F(u_{1},u_{2})\right]\right\|_{W^{s}}
\nonumber\\&&\geq C\left\|\sum_{j=1}^{4}\int_{\dot{Z}}\frac{|k|^{s+1}}{1+k^{2}}
\left[k^{2}+3+k_{1}k_{2}\right]R_{j}(k_{1},k_{2})|k_{\min}|^{-1/2}|k_{\max}|^{-j}dk_{1}\right\|_{l_{k}^{2}}\geq CN^{-j+2}.
\end{eqnarray*}
If (\ref{1.06}) is untrue, then we have that
\begin{eqnarray}
&&CN\leq \left\|\mathscr{F}^{-1}\left[\langle \tau+(-1)^{j}k^{2j+1}\rangle^{-1}\mathscr{F}F(u_{1},u_{2})\right]\right\|_{X_{s,\frac{1}{2}}}\nonumber\\&&\leq C
\left\|\mathscr{F}^{-1}\left[\langle \tau+(-1)^{j}k^{2j+1}\rangle^{-1}\mathscr{F}F(u_{1},u_{2})\right]\right\|_{W^{s}}\leq C\prod_{j=1}^{2}\|u_{j}\|_{W^{s}}\sim N^{2s}.\label{4.01}
\end{eqnarray}
Consequently, we obtain the contradiction since $s<-\frac{j}{2}+1.$
\end{proof}

We have completed the proof of Theorem 1.1.

\noindent {\large\bf 5. Proof of Theorem  1.2}

\setcounter{equation}{0}

 \setcounter{Theorem}{0}

\setcounter{Lemma}{0}

\setcounter{section}{5}

In this section, we prove Theorem 1.2.

Now we are in a position to prove Theorem 1.2.
 Let  $u^{\mu}:=\mu^{-2j}u^{\mu}(\mu^{-1} x,\mu^{-2j-1}t)$. Then, $u^{\mu}$  is the solution to
the following  problems
\begin{eqnarray}
&&u^{\mu}_{t}+\partial_{x}^{2j+1}u^{\mu}
     + \frac{1}{2}\partial_{x}((u^{\mu})^{2})+\partial_{x}(1-\mu^{2}\partial_{x}^{2})^{-1}\left[(u^{\mu})^2+\frac{1}{2}\mu^{2}(u_{x}^\mu)^{2}\right]
     = 0,\label{5.01}\\
    &&u^{\mu}(x,0)=\mu^{-2j}u_{0}(x/\mu):=u^{\mu}_{0}(x),\quad x\in \mathbf{T}=[0,2\pi\lambda\mu) \label{5.02}
\end{eqnarray}
if $u$ is the solution to (\ref{1.01})-(\ref{1.02}).
Let
\begin{eqnarray*}
&&F^{\mu}(t)=\frac{1}{2}\partial_{x}(u^{\mu})^{2}
+\partial_{x}(1-\mu^{2}\partial_{x}^{2})^{-1}\left[(u^{\mu})^{2}+\frac{1}{2}\mu^{2}(u_{x}^{\mu})^{2}\right].
\end{eqnarray*}
We define
\begin{eqnarray}
\Phi(u^{\mu})=\eta(t) S(t)u^{\mu}_{0}(x)-\frac{1}{2}\eta(t) \int_{0}^{t}S(t-t^{'})F^{\mu}(t^{'})dt^{'}.\label{5.03}
\end{eqnarray}
We claim that for $\|u_{0}^{\mu}\|_{H^{s}}\leq r$, there exists $u^{\mu}\in Z_{1}^{s}$ satisfying
\begin{eqnarray}
\Phi(u^{\mu})=u^{\mu}.\label{5.04}
\end{eqnarray}
By using   Lemmas 2.2,  2.3 , 3.3-3.5,  we have that
\begin{eqnarray*}
&&\left\|\Phi(u^{\mu})\right\|_{Z_{1}^{s}}\leq \left\|\eta(t) S(t)u_{0}^{\mu}\right\|_{Z_{1}^{s}}
+\left\|-\frac{1}{2}\eta(t) \int_{0}^{t}S(t-t^{'})
F^{\mu}(t^{'})dt^{'}\right\|_{Z_{1}^{s}}\nonumber\\&&\leq C_{1}\|u_{0}^{\mu}\|_{H^{s}(\mathbf{T})}
+C\left\|\partial_{x}((u^{\mu})^{2})\right\|_{Z_{1}^{s}}\leq C_{1}\left[\|u_{0}^{\mu}\|_{H^{s}(\mathbf{T})}+(1+\mu^{2})\|u^{\mu}\|_{Z_{1}^{s}}^{2}\right]
\end{eqnarray*}
and
\begin{eqnarray*}
&&\left\|\Phi(u^{\mu})-\Phi(v^{\mu})\right\|_{Z_{1}^{s}}\leq C_{1}(1+\mu^{2})\left\|u^{\mu}+v^{\mu}\right\|_{Z_{1}^{s}}\left\|u^{\mu}-v^{\mu}\right\|_{Z_{1}^{s}}\nonumber\\
&&\leq C_{1}(1+\mu^{2})\left[\left\|u^{\mu}\right\|_{Z_{1}^{s}}+\left\|v^{\mu}\right\|_{Z_{1}^{s}}\right]\left\|u^{\mu}-v^{\mu}\right\|_{Z_{1}^{s}}.
\end{eqnarray*}
Let
\begin{eqnarray*}
B=\left\{u\in Z_{1}^{s}: \quad \|u\|_{ Z_{1}^{s}}\leq (16C_{1})^{-1}\mu^{-\frac{j}{2}-1},j\geq2\right\}.
\end{eqnarray*}
Thus, if $\mu^{2}\geq \mu_{0}^{2}:=16C_{1}$ and
\begin{eqnarray*}
\|u_{0}^{\mu}\|_{H^{s}(\mathbf{T})}\leq (16C_{1})^{-2}\mu^{-j-\frac{1}{2}},
\end{eqnarray*}
 $\Phi$ will be a map from $B$  to itself and  $\Phi$ is a contraction map on $B.$ Thus,  the claim is valid.  Consequently, there exists a solution to (\ref{5.03}) for initial data  $\phi^{\mu}$
on the time interval $[-1,1].$ By using a similar manner, we can obtain the Lipschitz  continuity of the map $\Phi$. Next, we consider (\ref{1.01})-(\ref{1.02}) and $\|u_{0}\|_{H^{s}}\leq r$. If $r\leq (16C_{1})^{-2}\mu_{0}^{\frac{1}{2}}$, then we derive that
$\|u_{0}^{\mu_{0}}\|_{H^{s}}\leq \mu_{0}^{-j-1}\|u_{0}\|_{H^{s}}\leq (16C_{1})^{-2}\mu_{0}^{-j-\frac{1}{2}}$ and derive a solution $u^{\mu_{0}}$ to the $\mu_{0}$-rescaled
problem on $[0,1],$ thus derive a solution $u$ to (\ref{1.01}) with existence time $T=\mu ^{-2j-1}$.
If $(16C_{1})^{-2}\mu_{0}^{1/2}<r=:(16C_{1})^{-2}\mu^{1/2}$, by using the same way, we solve the $\mu(r)$-rescaled problem on $[-1,1]$ to obtain a solution
to (\ref{1.01}) with $T=\mu(r)^{-2j-1}.$

Now we prove the uniqueness of the solution.  Suppose that $u_{1}$ and $u_{2}$ are solutions to (\ref{1.01}) with the common data $u_{0}$ and the common existence time $T_{0}$ and $u$ and $v$ belong to $Z^{s}_{T_{0}}$. Then, $u_{k}^{\mu}(x,t)=\mu ^{-2j}u_{k}\left(\frac{x}{\mu^{2j+1}},\frac{t}{\mu}\right)$ with $k=1,2$ are solutions to (\ref{1.01})
corresponding to data $u_{0}^{\mu}=\mu^{-2j}u_{0}(\frac{x}{\mu^{2j+1}}).$ By using Lemmas 2.3, 2.4, 3.3, we have that
\begin{eqnarray}
&&\left\|u_{1}^{\mu}-u_{2}^{\mu}\right\|_{Z_{T}^{s}}\nonumber\\&&\leq C_{1}\left\|u_{1}^{\mu}+u_{2}^{\mu}\right\|_{Z_{T}^{s}}\left\|u_{1}^{\mu}-u_{2}^{\mu}\right\|_{Z_{T}^{s}}\nonumber\\
&&\leq C_{1}\left[\sum_{k=1}^{2}\left\|u_{k}^{\mu}\right\|_{Z_{T}^{s}}\right]\left\|u_{1}^{\mu}-u_{2}^{\mu}\right\|_{Z_{T}^{s}}\label{5.05}
\end{eqnarray}
for $0<T\leq {\rm min}\left\{1,\mu^{2j}T_{0}\right\}.$
From Lemma 2.3, we have that
\begin{eqnarray}
&&\|u_{k}^{\mu}\|_{Z_{T}^{s}}\leq \|u_{k}^{\mu}-e^{(-1)^{j}\partial_{x}^{2j+1}}u_{0}^{\mu}\|_{Z_{T}^{s}}+\|u_{0}^{\mu}\|_{Z_{T}^{s}}\nonumber\\&&\leq \|u_{k}^{\mu}-e^{(-1)^{j}\partial_{x}^{2j+1}}u_{0}^{\mu}\|_{Z_{T}^{s}}+\mu^{-2j}\|u_{0} \|_{H^{s}}\label{5.06}
\end{eqnarray}
with $k=1,2.$
Combining $(u_{k}^{\mu}-e^{(-1)^{j}\partial_{x}^{2j+1}}u_{0}^{\mu})_{t=0}$ with Lemma 2.6, we have that
$\|u^{\mu}_{k}-e^{(-1)^{j}\partial_{x}^{2j+1}}u_{0} ^{\mu}\|_{Z_{T}^{s}}\longrightarrow0$ as $T\longrightarrow0.$ We can choose a sufficiently large $\mu=\mu (\|u_{0}\|_{H^{s}})$,  such that
\begin{eqnarray}
2C_{1}C\mu^{-2j}\|u_{0}^{\mu}\|_{Z_{T}^{s}}\leq \frac{1}{4}\label{5.07}
\end{eqnarray}
and a sufficiently small $T=T(\mu,u_{1},u_{2})$ such that
\begin{eqnarray}
C_{1}\left[\sum_{k=1}^{2}\|u_{k}^{\mu}-e^{(-1)^{j}\partial_{x}^{2j+1}}u_{0}^{\mu}\|_{Z_{T}^{s}}\right]\leq \frac{1}{4}.\label{5.08}
\end{eqnarray}
Combining (\ref{5.05})-(\ref{5.06}) with (\ref{5.07})-(\ref{5.08}), we have that
\begin{eqnarray}
\left\|u_{1}^{\mu}-u_{2}^{\mu}\right\|_{Z_{T}^{s}}\leq \frac{1}{2}\left\|u_{1}^{\mu}-u_{2}^{\mu}\right\|_{Z_{T}^{s}}\label{5.09}.
\end{eqnarray}
Consequently, we derive that $u_{1}=u_{2}$ for $-\mu^{-2j}T\leq t\leq \mu^{-2j}T.$ If $\mu^{-2j}T=T_{0}$, the conclusion is valid. If $\mu^{-2j}T<T_{0}$, by using a continuity argument which can be seen in \cite{MT}, we can obtain the conclusion.

We have completed the proof of Theorem 1.2.

\leftline{\large \bf Acknowledgments}

\bigskip

\noindent

 This work is supported by the Natural Science Foundation of China
 under grant numbers 11171116 and 11401180. The first author is also
 supported in part by the Fundamental Research Funds for the
 Central Universities of China under the grant number 2012ZZ0072.
 The second author is  supported by the
 NSF of China (No.11371367) and Fundamental
 research program of NUDT(JC12-02-03).


\bigskip

  \bigskip


\leftline{\large\bf  References}



\begin{thebibliography}{99}

 \bibitem{BT}
 I.  Bejenaru, T. Tao,  Sharp well-posedness and ill-posedness results for a
quadratic non-linear Schr\"odinger equation, {\it J. Funct. Anal.} 233(2006), 228-259.



\bibitem{Bourgain93}
J. Bourgain, Fourier transform restriction phenomena for certain lattice subsets
and applications to nonlinear evolution equations,
part II: The KdV equation,  {\it Geom.   Funct. Anal.} 3(1993), 209-262.

\bibitem{Bourgain97}
 J. Bourgain, Periodic Korteweg  de vries equation with measures as initial data,
 {\it  Sel. Math.} 3(1997), 115-159.



\bibitem{BCARMA}
 A. Bressan, A. Constantin, Global  conservative solutions of the Camassa-Holm equation, {\it Arch. Ration. Mech. Anal.} 183(2007), 215-239.


\bibitem{CH} R. Camassa, D. Holm, An integrable shallow water equation with peaked solutions,
 {\it Phys. Rev. Lett.}  71(1993), 1661-1664.
 \bibitem{CHH} R. Camassa, D. Holm, J. Hyman, A new integrable shallow water equation, {\it Adv. Appl. Mech.} 31(1994), 1-33.

\bibitem{C} A. Constantin,  The Hamiltonian structure of the Camassa-Holm equation, {\it  Exposition. Math.} 15(1997),  53-85.
\bibitem{Con2000} A. Constantin, Existence of permanent and breaking waves for a shallow water equation: a geometric approach,
{\it Ann. Inst. Fourier (Grenoble)} 50(2000), 321-362.

\bibitem{C2000}A. Constantin, On the blow-up of solutions of a periodic shallow water equation, {\it  J. Nonlinear Sci.} 10(2000),  391-399.

\bibitem{C2001}A.  Constantin,  On the scattering problem for the Camassa-Holm equation, {\it Proc. R. Soc. Lond. Proc. Ser. A } 457(2001), 953-970.

\bibitem{C2006} A. Constantin, The trajectories of particles in Stokes waves, {\it Invent. Math.} 166(2006), 523-535.

\bibitem{CKSTT}  J. Colliander, M. Keel, G. Staffilani,  H. Takaoka,  T. Tao,
Sharp global well-posedness for KdV and modified KdV on $R$  and $T$,   {\it J. Amer. Math. Soc.}
16(2003),  705-749.
\bibitem{CE} A. Constantin, J. Escher,  Wave breaking for nonlinear nonlocal shallow water equations, {\it  Acta Math.} 181(1998),  229-243.

\bibitem{CECPAM}A. Constantin, J. Escher,  Well-posedness, global existence, and
blowup phenomena for a periodic quasi-linear hyperbolic equation, {\it  Comm. Pure Appl. Math.} 51(1998),  475-504.
\bibitem {CS} A. Constantin, W.  Strauss,  Stability of peakons, {\it  Comm. Pure Appl. Math.} 53(2000),  603-610.
\bibitem{EY}J. Escher, Z. Y. Yin, Initial boundary value problems of the
Camassa-Holm equation, {\it Comm. Partial   Diff. Eqns.} 33(2008),  377-395.

\bibitem{EYJFA} J. Escher, Z. Y. Yin, Initial boundary value problems
for nonlinear dispersive wave equations, {\it J. Funct. Anal.} 256(2009), 479-508.


\bibitem{FF}A. Fokas, B. Fuchssteiner, Symplectic structures, their
$B\ddot{a}$klund transformations and hereditary
symmetries, {\it  Phys. D.}  71(1981),  47-66.




\bibitem{Go} J. Gorsky, On the Cauchy problem for a KdV-type equation
 on the circle, dissertation, University of Notre Dame, 2004.

\bibitem{G} Z. H.  Guo,  Global well-posedness of Korteweg-de Vries equation in $H^{-3,4}(\R),$
{\it Journal de Math$\acute{e}$matiques Pures et Appliqu$\acute{e}$es}, 91(2009), 583-597.



\bibitem{H}
H. Hirayama, LocaL well-posedness for the  periodic higher order KdV type equations,
 {\it Nonlinear Differential equations and applications,}
19(2012), 677-693.

\bibitem{HM1998} A. A. Himonas, G. Misiolek, The  Cauchy   problem for
a shallow water type equation, {\it  Comm. Partial   Diff. Eqns.}
23(1998), 123-139.




\bibitem{HM2000} A. A. Himonas, G. Misiolek,  Well-posedness of the Cauchy
problem for a shallow water equation on the circle,
{\it J. Diff. Eqns.} 161(2000), 479-495.

\bibitem{HM}A. A. Himonas, G. Misiolek, The initial value problem for a
fifth order shallow water in: Analysis, Geometry, Number Theory:
The Mathematics of Leon Ehrenpreis, in: Contemp. Math. Vol. 251, Amer.
Math. Soc., Providence, RI, 2000, pp. 309-320.




\bibitem{HMPZ} A. A. Himonas,G.  Misiolek, G.  Ponce, Y. Zhou,  Persistence properties and unique continuation of solutions of the Camassa-Holm equation,
 {\it Comm. Math. Phys.} 271(2007), 511-522.



\bibitem{KT2006} T. Kappeler and P. Topalov, Global wellposedness of KdV in
$H^{-1}(T,R)$, {\it Duke Math. J.}  135(2006),  327-360.

\bibitem{TKato} T. K. Kato, Local well-posedness for Kawahara equation,
{Adv. Diff. Equ.} 16(2011), 257-287.

\bibitem{Kato}
T. Kato, Low regularity well-posedness for the periodic Kawahara  equation,
{\it Diff. Int. Eqns.} 25(2012), no. 11-12, 1011-1036.





\bibitem{IK}
A. D. Ionescu, C. E. Kenig, Global well-posedness of the Benjamin-Ono equation
in low-regularity spaces, {\it J. Amer. Math. Soc.} 20(2007), 753-798.

\bibitem{IKT}
A. D. Ionescu, C. E. Kenig, D. Tataru,  Global well-posedness of the KP-I
initial-value problem in the energy space, {\it Invent. Math.} 173(2008), 265-304.







\bibitem{KPV1996}
C. E. Kenig,  G. Ponce, L. Vega, A bilinear estimate with applications to the KdV equation,
 {\it J. Amer. Math. Soc.} 9(1996), 573-603.

\bibitem{KPV2001}
C. E. Kenig,  G. Ponce, L. Vega, On the ill-posedness of some canonical dispersive equations,
{\it Duke Math. J.} 106(2001), 617-633.




\bibitem{Kis}N. Kishimoto,
 Well-posedness of the Cauchy problem for the Korteweg-de
Vries  equation at the critical regularity, {\it Diff. Int. Eqns.}  22(2009), 447-464.

\bibitem{KJDE} N. Kishimoto, Low-regularity bilinear estimates for a quadratic nonlinear Schr\"odinger equation,
{\it J.  Diff. Eqns.} 247(2009),  1397-1439.

\bibitem{KJDE2013}N. Kishimoto,
Sharp local well-posedness for the ``good" Boussinesq equation,
{\it J. Diff. Eqns.} 254(2013), 2393-2433.




\bibitem{LIMRN} J. Lenells,  Stability of periodic peakons, {\it Int. Math. Res. Not.} 10(2004),  485-499.


\bibitem{LYLH}S. M. Li,  W. Yan, Y, S. Li, J. H. Huang, The Cauchy problem for   a higher order shallow water  type equation on the circle,
    {\it J. Diff. Eqns.}  259(2015),   4863-4896.


\bibitem{LYY}Y. S. Li,  W. Yan,   X. Y.  Yang,
Well-posedness of a higher order modified Camassa-Holm equation in spaces of low regularity,
{\it J. Evol. Eqns.} 10(2010),  465-486.


\bibitem{L} B. P. Liu, A priori bounds for KdV equation below $H^{-3/4}$,
{\it J. Funct. Anal.} 268(2015), 501-554.




\bibitem{Molinet} L. Molinet,  Sharp ill-posedness results for the KdV and mKdV
equations on the torus, {Advances in Mathematics,}  230(2012), 1895-1930.
\bibitem{MT} T. Muramatu, S. Taoka,
The initial data value problem for the 1-D semilinear Schr$\ddot{o}$dinger
equation in the Besov space,
{\it J. Math. Soc. Japan.} 56(2004), 853-888.

\bibitem{O} E. A. Olson,  Well posedness for a higher order modified Camassa-Holm
equation, {\it J. Diff. Eqns.}  246(2009), 4154-4172.
\bibitem{T}  T. Tao, Multilinear weighted convolution of $L^2$ functions, and
applicationsto non-linear dispersive equations, {\it Amer. J. Math. } 123(2001), 839-908.

\bibitem{WC}
H. Wang, S. B. Cui, Global well-posedness of the Cauchy problem of the fifth order
shallow water equation, {\it J. Diff. Eqns.} 230(2006), 600-613.

\bibitem{XZCPAM} Z. P. Xin, P. Zhang, On the weak solutions to a shallow water equation,
 {\it Comm. Pure Appl. Math.} 53(2000), 1411-1433.

\bibitem{XZCPDE} Z. P. Xin, P. Zhang,  On the uniqueness and
large time behavior of the weak solutions to a shallow water equation.
{\it Comm. Partial Diff. Eqns.} 27(2002),  1815-1844.



\bibitem{YJLH}  Wei Yan, Minjie Jiang, Yongsheng Li, Jianhua Huang,
 Sharp well-posedness and ill-posedness of the Cauchy problem for the higher-order KdV,
arXiv:1511.02430

\bibitem{YLZZ} Wei Yan, Yongsheng LI, Xiaoping Zhai, Yimin Zhang,   The Cauchy problem for the shallow water typ equations in low regularity spaces on the circle,
arXiv:1602.04533.








\end{thebibliography}
\end{document}